\documentclass[a4paper,10pt]{article}
\usepackage[T1]{fontenc}
\title{Heat Generation using Lorentzian Nanoparticles:\\ Estimation via Time-Domain Techniques}
\usepackage{amssymb}
\usepackage{graphicx}
\usepackage[cmtip,all]{xy}
\usepackage[dvipsnames]{xcolor}
\usepackage{amsmath,bm}
\usepackage{mathtools}
\usepackage[colorlinks,allcolors=blue]{hyperref} 
\usepackage[noabbrev,capitalize,nameinlink]{cleveref}
\newtheorem{theorem}{Theorem}[section]
\newtheorem{lemma}[theorem]{Lemma}
\newtheorem{corollary}{Corollary}[theorem]
\newtheorem{proposition}{Proposition}[theorem]

\numberwithin{equation}{section}
\usepackage[frak=esstix]{mathalpha}
\usepackage[margin=25mm]{geometry}
\usepackage{amsmath}
\usepackage{mathrsfs}
\DeclareRobustCommand{\rchi}{{\mathpalette\irchi\relax}}
\newcommand{\irchi}[2]{\raisebox{\depth}{$#1\chi$}}
\author{Arpan Mukherjee\footnote{Radon Institute (RICAM), Austrian Academy of
		Sciences, Altenberger Straße 69, 4040 Linz, Austria.\quad Email: arpan.mukherjee@oeaw.ac.at. This author is supported by the Austrian Science Fund (FWF): P32660.} \ and Mourad Sini\footnote{Radon Institute (RICAM), Austrian Academy of Sciences, Altenberger Straße 69, 4040 Linz, Austria. Email: mourad.sini@oeaw.ac.at. This author is partially supported by the Austrian Science Fund (FWF): P32660.}}
\begin{document}
	\maketitle
	\begin{abstract}
		\noindent
		We analyze the mathematical model that describes the heat generated by  electromagnetic nanoparticles. We use the known optical properties of the nanoparticles to control the support and amount of the heat needed around a nanoparticle. Precisely, we show that the dominant part of the heat around the nanoparticle is the electric field multiplied by a constant dependent, explicitly and only, on the permittivity and quantities related to the eigenvalues and eigenfunctions of the Magnetization (or the Newtonian) operator, defined on the nanoparticle, and inversely proportional to the distance to the nanoparticle. 
		\vskip 0.01in
		\noindent  
		The nanoparticles are described via the Lorentz model. If the used incident frequency is chosen related to the plasmonic frequency $\omega_p$ (via the Magnetization operator) then the nanoparticle behaves as a plasmonic one while if it is chosen related to the undamped resonance frequency $\omega_0$ (via the Newtonian operator), then it behaves as a dielectric one. The two regimes exhibit different optical behaviors. In both cases, we estimate the generated heat and discuss advantages of each incident frequency regime. 
		\vskip 0.01in
		\noindent  
		The analysis is based on time-domain integral equation techniques avoiding the use of (formal) Fourier type transformations.  
	\end{abstract}
	\section{\Large\textbf{Problem Formulation and the Main Results}}
	
	\subsection{General introduction}
	Heat generation is used in many applications including medical imaging and therapy \cite{A-B-B, B-C-R-2020, baffou, B-B-C, D-H-2013, G-S-K-S,  Maldovan}. For instance the photo-acoustic imaging modality is based on using the acoustic pressure fluctuations, collected in an accessible part of the region to image, to recover few optical, and eventually acoustical, properties of the tissue. This pressure is generated by the heat created after exciting locally the tissue with electric laser fields \cite{B-C-R-2020, baffou, Ahcene, Triki}. This heat generation phenomenon has been also proposed as a thermal therapy to cure anomalies (as tumors) by injecting nanoparticles in the region occupied by the anomaly \cite{ProfHabib, baffou}. The main principle behind this phenomenon can be described as follows. It is known that an electric laser field excites surface plasmons, on metallic nanoparticles, at optical frequencies. In turn, these nanoparticles produce heat from the absorbed energy that diffuses away from them to raise the temperature of the surrounding medium. 
	In this work, we mainly focus on this therapy application. For this application, the question raised is how to control the generated heat so that it is enough to clean the anomaly but not so high to harm the surrounding tissue.
	\vskip 0.1in
	\noindent
	To describe the mathematical model behind this heat generation using electromagnetic nanoparticles, we split it into parts.  In the first part, we describe electromagnetic wave propagation generated by injected nanoparticles and in the second one, we describe the related heat model.
	\vskip 0.1in
	\noindent
	Let us then recall the electromagnetic problem described by the time-dependent Maxwell equation
	\begin{equation}\label{eq:TM Maxwell}
		\begin{cases}   \nabla \times \textbf{E} = -\mu\frac{\partial}{\partial t}\textbf{H}\\[10pt]
			\nabla \times \textbf{H} = \varepsilon\frac{\partial}{\partial t}\textbf{E},
		\end{cases}
	\end{equation}
	where $\textbf{E}$ and $\textbf{H}$ are the total electric and magnetic field respectively. In addition, the coefficients $\mu$ and $\varepsilon$ are the magnetic permeability and electric permittivity, respectively. 
	\vskip 0.1in
	\noindent
	We assume that the nanoparticle occupy a bounded domain $\Omega$ in $\mathbb{R}^{2}$. In addition, denote $\Omega = \delta B + \mathrm{z}$, where $\delta$ denotes the size of the nanoparticle, $B$ is centred at origin and $\mathrm{z}$ represent the position of the nanoparticle.
	\vskip 0.1in
	\noindent
	Moreover, consider the electric permittivity and magnetic permeability, respectively, of the form $\varepsilon = {\varepsilon_{\mathrm{p}}}\rchi_{\Omega} + \varepsilon_{\mathrm{m}}\rchi_{\mathbb{R}^{2}\setminus\overline{\Omega}}$  $\mu = \mu_{\mathrm{p}}\rchi_{\Omega} + \mu_{\mathrm{m}}\rchi_{\mathbb{R}^{2}\setminus\overline{\Omega}}$. We denote by $\varepsilon_\mathrm{m} = \varepsilon_\infty \varepsilon'_\mathrm{m}$ to be the relative permittivity and $\mu_\mathrm{m} = \mu_{\infty}\mu'_\mathrm{m}$ the relative permeability of the host medium, respectively. Both are assumed to be constant and independent of the frequency $\omega$ of the incident wave (i.e. the host medium is non dispersive). Here, $\varepsilon_\infty$ and $\mu_{\infty}$ are the electric permittivity and magnetic permeability of the free space respectively. Next, we assume the nanoparticle to be nonmagnetic, i.e. $\mu_{\mathrm{p}} = \mu_{\infty}\mu'_{m}$. But its permittivity $\varepsilon_{\mathrm{p}}$ is given by the so-called Lorentz model which can be described as follows
	\begin{align}\label{Lorenz-model}
		\varepsilon_\mathrm{p}(\omega,\gamma) = \varepsilon_\infty\Bigg[1+\dfrac{\omega_\mathrm{p}^2}{\omega_0^2-\omega^2 - i\gamma\omega} \Bigg]
	\end{align}
	where $\omega^2_\mathrm{p}$ is the electric plasma frequency, $\omega^2_0$ is the undamped resonance frequency and $\gamma$ is the electric damping parameter, see Section \ref{Lorentz}.
	\vskip 0.1in
	\noindent
	We deal with two types of nanoparticles:
	\begin{enumerate}
		\item {\it{Plasmonic nanoparticles}}. These are characterized by the conditions that  $\boldsymbol{\Re}(\varepsilon_{\mathrm{p}})< 0$ and $\boldsymbol{\Im}(\varepsilon_{\mathrm{p}})>0$ and both of them have moderate amplitudes. The second condition represents the ability of the biological tissues to absorb electromagnetic radiation and hence create the heat while the first one is responsible for creating resonances (namely the plasmonic ones) that allow for controlling the heat.
		\item {\it{Dielectric nanoparticles}}. These are characterized by the conditions 
		$\boldsymbol{\Re}(\varepsilon_{\mathrm{p}})\gg 1$, and positive, and eventually $\boldsymbol{\Im}(\varepsilon_{\mathrm{p}})\gg 1$  but with a small ratio $\frac{\boldsymbol{\Im}(\varepsilon_{\mathrm{p}})}{\boldsymbol{\Re}(\varepsilon_{\mathrm{p}})}\ll 1$. Similarly as above, the second condition allows for heat generation while the second is responsible for creating the Dielectric resonances. With these resonances, we enhance the generated heat while with smallness of the above ratio, we control the upper bound of this heat.   
	\end{enumerate}
	\noindent
	As we can see it later, based on the Lorentz model (\ref{Lorenz-model}), the nanoparticle behaves as a plasmonic or as a dielectric one according to the regimes we choose for the incident frequency $\omega$. Precisely, if the incident frequency $\omega$ is chosen related to the the plasmonic frequency $\omega_p$ then it behaves as a plasmonic one while if it is chosen related to the undamped frequency $\omega_0$ then it behaves as a dielectric one. In both cases, one needs appropriate choices of the damping frequency $\gamma$ (but with small values), see Section \ref{plasmonics} and Section \ref{dielectric1}, respectively, for more details.
	
	\subsection{The electromagnetic models}
	In this work, we mainly focus on the $2\mathrm{D}$ model describing the time harmonic regime in TM-polarization and TE-polarization respectively. 
	\bigskip
	
	1.  First, in the TM-regime, we use incident waves of the form $\mathbb E^{in}:=e^{i\omega \sqrt{\mu_m \epsilon_m}(x \cdot \theta-t)}\; (1, 1, 0)$ and then the total field of the form
	$\textbf{{E}} = \mathbb{E}\; e^{-i\omega t}$ and $\textbf{H} =\mathbb{H}\; e^{-i\omega t}$, where $\mathbb{E}=(\mathrm{E}(\mathrm{x}_1,\mathrm{x}_2), 0):=(\mathrm{E}_1(\mathrm{x}_1,\mathrm{x}_2), \mathrm{E}_2(\mathrm{x}_1,\mathrm{x}_2), 0)$ and $\mathbb{H} = \Big(0,0,\mathrm{H}(\mathrm{x}_1,\mathrm{x}_2)\Big)$. From (\ref{eq:TM Maxwell}), we obtain the following equation
	\begin{equation}\label{eq:helmholtz}
		\nabla \cdot \left[ \dfrac{1}{\varepsilon} \nabla \mathrm{H} \right] +\omega^{2}\mu_{m} \mathrm{H} =0 \ \text{in} \ \mathbb{R}^{2} 
	\end{equation}
	where $\mathrm{H}=\mathrm{H}^{\textbf{in}}+\mathrm{H}^{\mathrm{s}}$ with $\mathrm{H}^{\textbf{in}}$ as the incident wave and $\mathrm{H}^{\mathrm{s}}$ as the scattered wave satisfying the Sommerfield radiation condition $\lim_{|\mathrm{x}|\to \infty}|\mathrm{x}|\Big(\frac{\partial\mathrm{H}^\mathrm{s}}{\partial|\mathrm{x}|}-\mathrm{i}\omega \sqrt{\varepsilon_\mathrm{m}\mu_\mathrm{m}}\ \mathrm{H}^\mathrm{s}\Big) = 0$. The incident magnetic field satisfies the following equation $
	\Delta \mathrm{H}^{\textbf{in}} + \omega^2\mu_{\mathrm{m}}\varepsilon_{\mathrm{m}}\mathrm{H}^{\textbf{in}} = 0 \quad \text{in} \quad \mathbb{R}^2$.
	Additionally, the electric field $\mathrm{E}$ can be calculated from the magnetic field using the following relation
	\begin{align}
		\mathrm{E}(\mathrm{x}) = \begin{pmatrix}
			\partial_{\mathrm{x_2}} \mathrm{H}(\mathrm{x}) \\
			-\partial_{\mathrm{x_1}} \mathrm{H}(\mathrm{x})
		\end{pmatrix}.
	\end{align}
	Therefore, it can be readily seen that $|\mathrm{E}(\mathrm{x})|^2 = |\nabla \mathrm{H}(\mathrm{x})|^2.$
	\vskip 0.1in
	
	2. Next, in the TE-regime, we use incident waves of the form $\mathbb E^{in}:=e^{i\omega \sqrt{\mu_m \epsilon_m}(x \cdot \theta-t)}\; (0, 0, 1)$ and then the total fields $\textbf{{E}} = \mathbb{E}\; e^{-i\omega t}$, where $\mathbb{E}=(0, 0,  \mathrm{E}_3(\mathrm{x}_1,\mathrm{x}_2))$, satisfies the following equation
	\begin{align}\label{el3}
		\Delta \mathrm{E} + \omega^2\mu_\mathrm{m}\varepsilon \mathrm{E} = 0 \quad \text{in} \quad \mathbb{R}^2,
	\end{align}
	and the corresponding scattered wave satisfies the Sommerfield radiation condition as $|\mathrm{x}|\to \infty$.
	\vskip 0.1in
	\noindent
	In the coming analysis, we use equation (\ref{el3}) to describe the electric field when it interacts with dielectric nanoparticles, while the model equation (\ref{eq:helmholtz}) will be used for plasmonic nanoparticles.
	
	\subsection{The heat generation model}
	As described above, laser electric fields can generate surface plasmons at certain ranges of frequencies so that the energy absorbed by them generates heat. The temperature at the surrounding medium created when the heat diffuses from the nanoparticle, is governed by the following model \cite{ProfHabib}, \cite{baffou}
	\begin{equation}\label{eq:heat}
		\begin{cases}   \rho \mathrm{c}\frac{\partial u}{\partial t} - \nabla.\ \gamma\nabla u = \frac{\omega}{2\pi} \boldsymbol{\Im}(\varepsilon)|\mathrm{E}|^{2} \ \ in \ \ (\mathbb{R}^{2}\setminus \partial\Omega) \times (0,\mathrm{T}),\\
			\gamma^{\textbf{int}}_{0}u \ - \ \gamma^{\textbf{ext}}_{0}u=0    \ \ \ \ \ \ \ \ \ \ \ \ \ \ \ \ \ on \ \partial \Omega,\\
			\gamma_{\mathrm{\mathrm{p}}}\gamma^{\textbf{int}}_{1}u\ - \ \gamma_{\mathrm{m}} \gamma^{\textbf{ext}}_{1}u= 0 \ \ \ \ \ \ \ \ \ \ on \ \partial \Omega, \\
			u(\mathrm{x},0) = 0 \ \ \ \ \ \ \ \ \ \ \ \ \ \ \ \ \ \ \ \ \ \ \ \ \ \ \ \ \text{for} \ \mathrm{x} \in \mathbb{R}^2,
		\end{cases}
	\end{equation}
	where $\rho = \rho_{\mathrm{p}}\rchi_{\Omega} + \rho_{m}\rchi_{\mathbb{R}^{2}\setminus\overline{\Omega}}$ is the mass density; $\mathrm{c} = \mathrm{c}_{\mathrm{p}}\rchi_{\Omega} + \mathrm{c}_{\mathrm{m}}\rchi_{\mathbb{R}^{2}\setminus\overline{\Omega}}$ is the thermal capacity; $\gamma = \gamma_{\mathrm{p}}\rchi_{\Omega} + \gamma_{m}\rchi_{\mathbb{R}^{2}\setminus\overline{\Omega}}$ is the thermal conductivity and we recall that $\varepsilon = \varepsilon_\mathrm{p} \rchi_{\Omega} + \varepsilon_\mathrm{m}\rchi_{\mathbb{R}^{2}\setminus\overline{\Omega}}$ is the electric permittivity respectively. Here, $\mathrm{T}\in \mathbb{R}$ is the final time of measurement.\\
	In the sequel, we use the notations for both an interior Dirichlet and Neumann trace, whenever they make sense,
	\begin{equation*}
		\gamma^{\textbf{int}}_{0} \mathrm{u}(\mathrm{x},\mathrm{t}) := \lim\limits_{\Omega \ni \Tilde{\mathrm{x}}\to \mathrm{x}\in \partial\Omega}\mathrm{u}(\mathrm{\Tilde{\mathrm{x}}},\mathrm{t}) \ \text{for} \ (\mathrm{x},\mathrm{t}) \in \partial\Omega \times \mathbb{R},
	\end{equation*}
	and 
	\begin{equation*}
		\gamma^{\textbf{int}}_{1} \mathrm{u}(\mathrm{x},\mathrm{t}) := \lim\limits_{\Omega \ni \Tilde{\mathrm{x}}\to \mathrm{x}\in \partial\Omega}\nu_{\mathrm{x}} \cdot \nabla_{\Tilde{\mathrm{x}}}\mathrm{u}(\mathrm{\Tilde{\mathrm{x}}},\mathrm{t}) \ \text{for} \ (\mathrm{x},\mathrm{t}) \in \partial\Omega \times \mathbb{R},\  \text{respectively}.
	\end{equation*}
	We use similar notations for the exterior Dirichlet and Neumann traces $\gamma^{\textbf{ext}}_{0} \ \text{and} \ \gamma^{\textbf{ext}}_{1}$ respectively.
	\bigbreak
	\noindent
	Furthermore, with $\mathrm{T}_{0}$ fixed, considering the fact that $u=0$ for $t<0$, we let $\mathrm{U}$ to be the solution of the following problem, stated in the whole time domain $\mathbb{R}$
	\begin{equation}\label{eq:heat}
		\begin{cases}   \rho \mathrm{c}\frac{\partial \mathrm{U}}{\partial t} - \nabla.\ \gamma\nabla \mathrm{U} = \frac{\omega}{2\pi} \boldsymbol{\Im}(\varepsilon)|\mathrm{E}|^{2}\rchi_{(0,T_{0})} \ \ in \ \ (\mathbb{R}^{2}\setminus \partial\Omega) \times \mathbb{R},\\
			\gamma^{\textbf{int}}_{0}\mathrm{U} \ - \ \gamma^{\textbf{ext}}_{0}\mathrm{U}=0    \ \ \ \ \ \ \ \ \ \ \ \ \ \ \ \ \ on \ \partial \Omega \times \mathbb{R},\\
			\gamma_{\mathrm{p}}\gamma^{\textbf{int}}_{1}\mathrm{U}\ - \ \gamma_{\mathrm{m}} \gamma^{\textbf{ext}}_{1}\mathrm{U}= 0 \ \ \ \ \ \ \ \ \ \ on \ \partial \Omega \times \mathbb{R}.
		\end{cases}
	\end{equation}
	As we have $\mathrm{U} = u$ on $\mathbb{R}^2\times (-\infty,T_{0})$ then to analyse $u \ \text{in}\ (0,T_{0})$, it is enough to study $\mathrm{U}.$
	\bigbreak
	\noindent
	We further assume, for the purpose of simplicity of the analysis, that $\rho_{\mathrm{p}},\rho_{m},\mathrm{c}_{\mathrm{p}},\mathrm{c}_{\mathrm{m}},\gamma_{\mathrm{p}},\gamma_{m}$ are positive constants. We also note that the heat equation becomes a homogeneous equation outside of $\Omega$ as $\boldsymbol{\Im}(\varepsilon) = 0$ in $(\mathbb{R}^{2}\setminus \overline{\Omega})$. Consequently, we can rewrite the governing heat equations and transmissions as follows
	\begin{align}\label{eq:heat1}
		\begin{cases}   \frac{\rho_{\mathrm{p}}\mathrm{c}_{\mathrm{p}}}{\gamma_{\mathrm{p}}}\frac{\partial \mathrm{U}_{\mathrm{i}}}{\partial t} - \Delta \mathrm{U}_{\mathrm{i}}= \frac{\omega}{2\pi\gamma_{\mathrm{p}}}\boldsymbol{\Im}(\varepsilon_\mathrm{p})|\mathrm{E}|^{2}\rchi_{(0,\mathrm{T_0})} \ \ in \ \ \Omega \times \mathbb{R}\\
			\frac{\rho_{\mathrm{m}}\mathrm{c}_{\mathrm{m}}}{\gamma_{\mathrm{m}}}\frac{\partial \mathrm{U}_{\mathrm{e}}}{\partial t} - \Delta \mathrm{U}_{\mathrm{e}} = 0 \ \ \ \ \ \ \ \ \ \ \ \ \ \ in \ \mathbb{R}^{2}\setminus \overline\Omega \times \mathbb{R}\\
			\gamma^{\textbf{int}}_{0}\mathrm{U}_{\mathrm{i}} \ - \ \gamma^{\textbf{ext}}_{0}\mathrm{U}_{\mathrm{e}}=0    \ \ \ \ \ \ \ \ \ \ \ \ \  on \ \partial \Omega \times \mathbb{R},\\
			\gamma_{\mathrm{p}}\gamma^{\textbf{int}}_{1}\mathrm{U}_{\mathrm{i}}\ - \ \gamma_{m}\gamma^{\textbf{ext}}_{1}\mathrm{U}_{\mathrm{e}}= 0 \ \ \ \ \ \ \ on \ \partial \Omega \times \mathbb{R}.
		\end{cases}
	\end{align}
	For ease and clarity of notation, throughout this work we denote the diffusion constants by $\alpha_{\mathrm{p}} := \frac{\rho_{\mathrm{p}}\mathrm{c}_{\mathrm{p}}}{\gamma_{\mathrm{p}}}$ and $\alpha_{m} := \frac{\rho_{\mathrm{m}}\mathrm{c}_{\mathrm{m}}}{\gamma_{\mathrm{m}}}$ respectively.
	\vskip 0.01in
	\noindent
	We assume that the nanoparticle has the following scales regarding the heat-related coefficients
	\begin{equation}\label{Assumption-heat-coefficients-section 1.3}
		\gamma_p \sim \delta^{-2} ~~\mbox{ and }~~ \rho_{\mathrm{p}}\mathrm{c}_{\mathrm{p}} \sim 1, ~~ \delta \ll 1.
	\end{equation}
	More general scales could be considered, namely $\gamma_p \sim \delta^{-m} ~~\mbox{ and }~~ \rho_{\mathrm{p}}\mathrm{c}_{\mathrm{p}} \sim \delta^{-n}, ~~ \delta \ll 1$, with $m$ and $n$ non-negative such that $m-n-1\geq 0$. \footnote{This last condition ensures that the diffusion coefficient $\alpha_p$ is small so that the, corresponding, last part of (\ref{Assumption-heat-coefficients}) in Theorem \ref{th2} is valid, see precisely (\ref{condition-apprio-estimate}).} To avoid introducing more parameters, we stick to the case $m=2, n=0$. Materials enjoying such scales can be found in the applied literature as \cite{biotissue, G-S-K-S} for instance.  
	
	\subsection{Statement of the results}
	\subsubsection{The electromagnetic field generated by plasmonic nanoparticles}
	We recall the fundamental solution of the Helmholtz equation in dimension two given by 
	\begin{equation*}
		\mathbb{G}^{(\mathrm{k})}(\mathrm{x},\mathrm{y}) = \dfrac{i}{4}\mathcal{H}^{(1)}_{0}(\mathrm{k}|\mathrm{x}-\mathrm{y}|), \quad \mathrm{x}\neq \mathrm{y},
	\end{equation*}
	where $\mathcal{H}^{(1)}_{0}$ is the Hankel function of the first kind of order zero.
	\hfill \break
	Moreover, we also introduce the Magnetic operator $\mathbb{M}: \nabla H_{arm}\rightarrow \nabla H_{arm},$ \ defined as follows
	\begin{equation}\label{magnetic operator}
		\mathbb{M}\Big[\nabla \mathrm{H}\Big](\mathrm{x}) = \nabla \int_{\Omega}  \nabla\mathbb{G}^{(0)}(\mathrm{x},\mathrm{y})\cdot\nabla \mathrm{H}(\mathrm{y})d\mathrm{y}.
	\end{equation}
	We also need to introduce the following decomposition of the space $\mathbb{L}^{2}(\Omega)$  into the following three sub-spaces as a direct sum as following, see \cite{raveski} for more details
	\begin{equation*}
		\mathbb{L}^{2} = \textbf{H}_{0}(\textbf{div},0) \oplus\textbf{H}_{0}(\textbf{curl},0)\oplus \nabla \mathrm{H}_{\textbf{arm}},
	\end{equation*}
	where we define these three sub-spaces as follows:
	\begin{equation}{\label{eq:subspaces}}
		\begin{cases}
			\textbf{H}_{0}(\textbf{div},0) = \left\{ u\in \mathbb{L}^{2}(\Omega): \nabla \cdot u = 0 \ \text{and} \ u\cdot \nu = 0 \right\},\\
			\textbf{H}_{0}(\textbf{curl},0) = \left\{ u\in \mathbb{L}^{2}(\Omega): \nabla \times u = 0 \ \text{and} \ u\times \nu = 0 \right\}, \\
			\nabla \mathrm{H}_{\textbf{arm}} = \left\{ u \in \mathbb{L}^{2}(\Omega): \exists \  \varphi \ \text{s.t.} \ u = \nabla \varphi \ \text{and} \ \Delta \varphi = 0 \right\}. 
		\end{cases}
	\end{equation}
	It is well known, see for instance \cite{friedmanI}, that the  Magnetization operator $\mathbb{M}: \nabla \mathrm{H}_{\text{arm}}\rightarrow \nabla \mathrm{H}_{\text{arm}}$ induces a complete orthonormal basis namely $\big(\lambda^{(3)}_{\mathrm{n}},\mathrm{e}^{(3)}_{\mathrm{n}}\big)_{\mathrm{n} \in \mathbb{N}}$. We are now able to present the first results of this research.
	\begin{theorem}\label{th1}
		Let a plasmonic nanoparticle occupying a domain $\Omega = \mathrm{z} + \delta B$ which is of class $\mathcal{C}^2$. We choose the incident frequency of the form
		
		\begin{equation}\label{choice-fre-plasmonic}
			\omega^2 = \omega_0^2 + \omega_\mathrm{p}^2\dfrac{\varepsilon_\mathrm{m}+\lambda_\mathrm{n_0}^{(3)}}{\lambda_\mathrm{n_0}^{(3)}(1-\varepsilon_\infty^{-1}\varepsilon_\mathrm{m})+\varepsilon_\mathrm{m}} + \mathcal{O}(\delta^\mathrm{h}) \mbox{ and } \gamma\omega \sim \delta^h.
		\end{equation}
		\noindent
		Then we have the following approximation of the electric field with $|\mathrm{E}|^2=|\nabla \mathrm{H}|^2$ , with $\mathrm{H}$ as the solution to (\ref{eq:helmholtz}), as $\delta\to 0$,
		\begin{align}
			\int_{\Omega}|\mathrm{E}|^2(\mathrm{y})d\mathrm{y} = \dfrac{1}{ |1 - \alpha\lambda^{(3)}_{\mathrm{n}_{0}}|^2}\Bigg[|\mathrm{E}^{\textbf{in}}|^2(\mathrm{z})\Big(\int_\mathrm{\Omega}\mathrm{e}^{(3)}_{\mathrm{n}_{0}}(\mathrm{x})d\mathrm{x}\Big)^2 + \mathcal{O}\Big(\delta^{3}\Big)\Bigg] \quad \text{with} \quad \mathrm{h}<1,
		\end{align}
		where $\alpha = \dfrac{1}{\varepsilon_{\mathrm{p}}(\omega)} - \dfrac{1}{\varepsilon_{m}}.$ In addition, we have $|1 - \alpha\lambda^{(3)}_{\mathrm{n}_{0}}| \sim \delta^h$ and $\Big(\int_\mathrm{\Omega}\mathrm{e}^{(3)}_{\mathrm{n}_{0}}(\mathrm{x})d\mathrm{x}\Big)^2\sim\delta^2$.
	\end{theorem}
	
	
	\subsubsection{The electromagnetic field generated by dielectric nanoparticles}
	Let us recall that the volumetric Logarithmic potential operator $\displaystyle \int_\Omega -\frac{1}{2\pi}\log|\mathrm{x}-\mathrm{y}|\mathrm{E}(\mathrm{y})d\mathrm{y}$ has a countable sequence of eigenvalues $\lambda_\mathrm{n}^{(\boldsymbol{\ell})}$ and the corresponding eigen-functions $\mathrm{e}_\mathrm{n}^{(\boldsymbol{\ell})}$ that form an orthonormal basis in $L^2(\Omega)$. In the sequel, we will need the following properties: $\displaystyle\int_\Omega\mathrm{e}_\mathrm{n}^{(\boldsymbol{\ell})}(\mathrm{x})d\mathrm{x} \neq 0$ and $\lambda_\mathrm{n}^{(\boldsymbol{\ell})}\sim \delta^2|\log\delta|.$ These properties are shown to be true for the first eigenvalue/eigen-function $n=1$ for $\Omega$ being the disc of radius $\delta$, see \cite{Ahcene}. We state the following theorem which was first derived in \cite{Ahcene} but we extend it to the case of the Lorentz model for the  permittivity.
	\begin{theorem}\label{1.2}
		Let a dielectric nanoparticle occupying a domain $\Omega = \mathrm{z} + \delta B$ which is of class $\mathcal{C}^2$. \footnote{Contrary to the plasmonic case, here we can reduce the regularity of $\Omega$ to Lipschitz.}  We choose the incident frequency of the form
		\begin{equation}\label{choice-fre-dielectric}
			\omega^2 -\omega_0^2\sim  -\delta^2|\log\delta|\Big(\overline{\lambda}_\mathrm{n_0}^{(\boldsymbol{\ell})}\mu_\mathrm{m}\omega_0^2\Big) \mbox{ and } \gamma\omega \sim \delta^2|\log\delta|^{1-\mathrm{h}-\mathrm{s}}\Big(\overline{\lambda}_\mathrm{n_0}^{(\boldsymbol{\ell})}\mu_\mathrm{m}\omega_0^2\Big)^2.
		\end{equation}
		\noindent
		Then an approximate representation of the electric field with $\mathrm{E}$ as a solution to (\ref{el3}) is as follows for $\delta\to 0$
		\begin{align}
			\displaystyle\int_\Omega|\mathrm{E}|^2(\mathrm{y})d\mathrm{y} = \dfrac{1}{|1-\omega^2\mu_\mathrm{m}\varepsilon_\mathrm{p}\lambda_\mathrm{n_0}^{(\boldsymbol{\ell})}|^2} \Big[|\mathrm{E}^{\textbf{in}}|^2(\mathrm{z)}\Big(\displaystyle\int_\Omega \mathrm{e}_\mathrm{n_0}^{(\boldsymbol{\ell})}(\mathrm{y})d\mathrm{y}\Big)^2 \; + \mathcal{O}\Big(\delta^2|\log\delta|^{\mathrm{h}-1}\Big)\Big] \quad \text{with} \quad \mathrm{h}<1.
		\end{align}
		In addition, we have $|1-\omega^2\mu_\mathrm{m}\varepsilon_\mathrm{p}\lambda_\mathrm{n_0}^{(\boldsymbol{\ell})}|\sim (\log \delta)^{-h}$ and $\Big(\int_\mathrm{\Omega}\mathrm{e}^{(\boldsymbol{\ell})}_{\mathrm{n}_{0}}(\mathrm{x})d\mathrm{x}\Big)^2\sim \delta^2$.
	\end{theorem}

	\subsubsection{The heat generated by the plasmonic/dielectric nanoparticles}
	
	The fundamental solution of the heat operator $\alpha\; \partial_t-\Delta$ for the two dimensional spatial space is given by the expression
	\begin{equation}
		\Phi(\mathrm{x},t;\mathrm{y},\tau):=\  \begin{cases}
			\frac{\alpha}{4\pi(t-\tau)}\textbf{exp}\Bigg(-\frac{\alpha |\mathrm{x}-\mathrm{y}|^2}{4(t-\tau)}\Bigg), \ \ \ t > \tau \\
			0 ,\quad \text{otherwise}
		\end{cases}
	\end{equation}
	The fundamental solutions for the interior and exterior heat equation (\ref{eq:heat1}) are $\Phi(\mathrm{x},t;\mathrm{y},\tau)$ and $\Phi^\mathrm{e}(\mathrm{x},t;\mathrm{y},\tau)$ respectively, which depend on the variables $\alpha_{\mathrm{p}}$ and $\alpha_\mathrm{m}$.
	\bigbreak
	\noindent
	We state the main result of our work.
	\begin{theorem}\label{th2}
		Let a Lorentzian nanoparticle, occupying a domain $\Omega = \mathrm{z} + \delta B$ which is of class $\mathcal{C}^2$, be such that its heat coefficients $(\rho_p, C_c, \gamma_c)$ satisfy the conditions
		\begin{equation}\label{Assumption-heat-coefficients}
			\gamma_p = \overline{\gamma}_p\; \delta^{-2} ~~\mbox{ and }~~ \rho_{\mathrm{p}}\mathrm{c}_{\mathrm{p}} \sim 1, ~~\mbox{such that}~~ \gamma_m <\sqrt{\overline{\gamma}_p\; \rho_{\mathrm{p}}\mathrm{c}_{\mathrm{p}} }, ~~ \delta \ll 1.
		\end{equation}
		
		Let $\xi$ such that $\textbf{dist}(\xi ,\Omega)\sim \delta^\mathrm{p}$ $\Big(|\xi -z|\sim \delta^\mathrm{p}+\delta\Big)$. 
		
		\begin{enumerate}
			
			\item If we use incident frequency $\omega$ satisfying (\ref{choice-fre-plasmonic}), then for $\mathrm{r}<\frac{1}{2}$, if $\frac{1+2\mathrm{p}(1-\mathrm{r})}{2} < \mathrm{h}<1 $, the heat conducted by the \underline{plasmonic nanoparticle}, as a solution to (\ref{eq:heat1}), is given by, as $\delta\to 0$,
			\begin{align}
				\mathrm{U}_{\mathrm{e}}(\xi,t) &= \frac{\gamma_{\mathrm{p}}}{\gamma_{m}}\frac{1}{\alpha_\mathrm{m}}\Bigg[\frac{\omega \cdot \boldsymbol{\Im}(\varepsilon_\mathrm{p})}{2\pi\gamma_{\mathrm{p}}}\int_{0}^{t} \Phi^{\textbf{e}}(\xi,t;z,\tau)d\tau\int_{\Omega}|\mathrm{E}|^{2}(\mathrm{y}) d\mathrm{y} + \mathcal{O}\Bigg(\frac{\omega \cdot \boldsymbol{\Im}(\varepsilon_\mathrm{p})}{2\pi}\delta^{3-2\mathrm{p}(1-\mathrm{r})}\sqrt{\mathcal{K}^{(\mathrm{T_0})}_{\mathrm{r}}} \Bigg)\Bigg].
			\end{align}
			
			\item If we use incident frequency $\omega$ satisfying (\ref{choice-fre-dielectric}), then for $2\mathrm{p}(1-\mathrm{r})<1$, the heat conducted by the \underline{dielectric nanoparticle}, as a solution to (\ref{eq:heat1}), is given by, as $\delta\to 0$,
			\begin{align}
				\mathrm{U}_{\mathrm{e}}(\xi,t) &= \frac{\gamma_{\mathrm{p}}}{\gamma_{\mathrm{m}}}\frac{1}{\alpha_\mathrm{m}}\Bigg[\frac{\omega \cdot \boldsymbol{\Im}(\varepsilon_\mathrm{p})}{2\pi\gamma_{\mathrm{p}}}\int_{0}^{t} \Phi^{\textbf{e}}(\xi,t;\mathrm{z},\tau)d\tau\int_{\Omega}|\mathrm{E}|^{2}(\mathrm{y}) d\mathrm{y} + \mathcal{O}\Bigg(\frac{\omega \cdot \boldsymbol{\Im}(\varepsilon_\mathrm{p})}{2\pi}\delta^{5-2\mathrm{p}(1-\mathrm{r})}|\log\delta|^{\frac{3\mathrm{h}}{2}}\sqrt{\mathcal{K}^{(\mathrm{T_0})}_{\mathrm{r}}}\Bigg)\Bigg].
			\end{align}
		\end{enumerate}
		
		Here, $\displaystyle\mathcal{K}^{(\mathrm{T_0})}_{\mathrm{r}}:=\sup_{t\in (0, T_0)}\int^{T_0}_0\frac{1}{(t-\tau)^{2r}}d\tau$ and it makes sense if $\mathrm{r}< \frac{1}{2}$.
	\end{theorem}
	\noindent
	The above approximations can be further detailed by the following observations. Let $\mathrm{s} := \frac{|\xi-\mathrm{z}|}{2\sqrt{t-\tau}}$. Then, as $\displaystyle\int_{0}^{t} \Phi^{\textbf{e}}(\xi,t;z,\tau)d\tau$ is independent of the space variable, and $t\le \mathrm{T}_0$, it can be shown that 
	\begin{align}
		\int_{0}^{t} \Phi^{\textbf{e}}(\xi,t;z,\tau)d\tau = \frac{1}{2} \Large{\Gamma}\Big(0,\frac{|\xi-\mathrm{z}|^2}{4t}\Big).
	\end{align}
	Furthermore, for $|\xi-\mathrm{z}|<<t$, it is well known that 
	\begin{align}
		\Large{\Gamma}\Big(0,\frac{|\xi-\mathrm{z}|^2}{4\mathrm{t}}\Big) = \mathcal{E}_1\Big(\frac{|\xi-\mathrm{z}|^2}{4\mathrm{t}}\Big),
	\end{align} 
	where $\mathcal{E}_1$ is the exponential integral. Now, the exponential integral $\mathcal{E}_1$ can be expanded as follows
	\begin{align}
		\mathcal{E}_1\Big(\frac{|\xi-\mathrm{z}|^2}{4\mathrm{t}}\Big) = -\gamma - \ln \frac{|\xi-\mathrm{z}|^2}{4\mathrm{t}} + \mathrm{e}_{1}\Big(\frac{|\xi-\mathrm{z}|^2}{4\mathrm{t}}\Big),
	\end{align}
	where $\gamma$ is the Euler–Mascheroni constant and $\mathrm{e_1}$ is a smooth function.
	\bigbreak
	\noindent
	In addition, let us recall that $|1 - \alpha\lambda^{(3)}_{\mathrm{n}_{0}}|=\mathcal{O}(\delta^\mathrm{h})$. Therefore, we can derive a more precise dominating terms as follows
	\begin{corollary} Under the assumptions of Theorem \ref{th2} and with the combination of Theorem \ref{th1}, we have the following approximation for a \underline{plasmonic nanoparticle} occupying a domain $\Omega = \mathrm{z} + \delta B$, as $\delta \to 0$,
		\begin{align}\label{Precise-heat-plasmonics}
			\nonumber
			\mathrm{U}_{\mathrm{e}}(\xi,t) &= \frac{\omega \cdot \boldsymbol{\Im}(\varepsilon_\mathrm{p})}{2\pi\rho_\mathrm{m} \mathrm{c}_\mathrm{m} }\;\Bigg[ \vert \mathrm{E}^{\textbf{in}} (\mathrm{z})\vert^2 \Big(\int_{\Omega}{e}^{(3)}_{\mathrm{n}_{0}}(\mathrm{x})d\mathrm{x} \Big)^2\; \log|\xi-\mathrm{z}|^{-1}\delta^{-2h}+ \mathcal{O}\Big(\delta^{3-2\mathrm{h}}\log|\xi-\mathrm{z}|^{-1}\Big) \\ &+ \mathcal{O}\Big(\delta^{1-2\mathrm{p}(1-\mathrm{r})}\sqrt{\mathcal{K}^{(\mathrm{T_0})}_{\mathrm{r}}} \Big)\Bigg].
		\end{align}
	\end{corollary}
	Similarly, let us recall that $|1-\omega^2\mu_\mathrm{m}\varepsilon_\mathrm{p}\lambda^{(\boldsymbol{\ell})}_\mathrm{n_0}| = \mathcal{O}(|\log\delta|^{-\mathrm{h}})$, where the electric permittivity $\varepsilon_\mathrm{p}$ derived from the Lorentz model for a dielectric nanoparticle. We state the following corollary.
	\begin{corollary}
		Under the assumptions of Theorem \ref{th2} and with the combination of Theorem \ref{1.2}, we have the following approximation for a \underline{dielectric nanoparticle} occupying a domain $\Omega = \mathrm{z} + \delta B$, as $\delta \to 0$,
		\begin{align}\label{Precise-heat-dielectrics}
			\nonumber
			\mathrm{U}_{\mathrm{e}}(\xi,t) &= \frac{\omega \cdot \boldsymbol{\Im}(\varepsilon_\mathrm{p})}{2\pi\rho_\mathrm{m} \mathrm{c}_\mathrm{m} } \Bigg[\vert \mathrm{E}^{\textbf{in}} (\mathrm{z})\vert^2 \Big(\int_{\Omega}\mathrm{e}^{(\boldsymbol{\ell})}_{\mathrm{n}_{0}}(\mathrm{x})d\mathrm{x}\Big)^2  \log|\xi-\mathrm{z}|^{-1}|\log\delta|^{2\mathrm{h}}+ \mathcal{O}\Big(\delta^2|\log\delta|^{3\mathrm{h}-2}\log|\xi-\mathrm{z}|^{-1}\Big) \\ &+ \mathcal{O}\Big(\delta^{3-2\mathrm{p}(1-\mathrm{r})}|\log\delta|^{\frac{3\mathrm{h}}{2}}\sqrt{\mathcal{K}^{(\mathrm{T_0})}_{\mathrm{r}}} \Big)\Bigg].
		\end{align}
	\end{corollary}
	
	\subsubsection{Discussion about the obtained results}
	
	We make the following two observations based on the two formulas (\ref{Precise-heat-plasmonics}) and (\ref{Precise-heat-dielectrics}).
	\begin{enumerate}
		\item Using plasmonics. In this case, according to the choice of the incident frequency, from the Lorentz model, we have $\boldsymbol{\Im}(\varepsilon_\mathrm{p}) \sim \delta^h$. As \ $\ln{\vert \xi-z\vert}\sim  p\; \vert\ln \delta\vert$, then from (\ref{Precise-heat-plasmonics}), recalling that $\displaystyle\Big(\int_{\Omega}\mathrm{e}^{(3)}_{\mathrm{n}_{0}}(\mathrm{x})d\mathrm{x}\Big)^2 \sim \delta^2$, we have
		\begin{equation}\label{exact-amount-plasmonics}
			\mathrm{U}_{\mathrm{e}}(\xi,t)\sim p\; \vert \ln \delta\vert \delta^{2-h}(1+o(1)),\; \delta \ll 1.
		\end{equation}
		We observe that the amount of heat can be enhanced till the order near to $\delta$.
		
		\item Using dielectrics. In this case, we have $\boldsymbol{\Im}(\varepsilon_\mathrm{p}) \sim \delta^{-2} (\vert \ln \delta \vert)^{-1-\mathrm{h}-\mathrm{s}}$. Therefore from (\ref{Precise-heat-dielectrics}), recalling that  $\displaystyle \Big(\int_{\Omega}\mathrm{e}^{(\boldsymbol{\ell})}_{\mathrm{n}_{0}}(\mathrm{x})d\mathrm{x}\Big)^2 \sim \delta^2$, we have 
		\begin{equation}\label{exact-amount-dielectrics}
			\mathrm{U}_{\mathrm{e}}(\xi,t)\sim p\; \vert \ln \delta\vert^{h-s}(1+o(1)),\; \delta \ll 1.
		\end{equation}
		Choosing $\mathrm{s}=\mathrm{h}$, we have $\mathrm{U}_{\mathrm{e}}(\xi,t)\sim p$. Choosing $\mathrm{s}<\mathrm{h}$ we get too much heat and for $s>h$, we get less heat. 
	\end{enumerate}
	\bigbreak
	\noindent
	Recall that $dist(\xi, \Omega)\sim \delta^\mathrm{p}$. In both the formulas (\ref{exact-amount-plasmonics}) and (\ref{exact-amount-dielectrics}), the parameter $p$ indicates that the heat decreases when $dist(\xi, \Omega)$ increases. For instance, the heat at the distance $dist(\xi, \Omega)=\delta^{\frac{1}{10}}$ is of the order $\frac{1}{5}$ the heat generated at a distance $dist(\xi, \Omega)=\delta^{\frac{1}{2}}$.
	\bigbreak
	\noindent
	In both cases, the multiplicative terms in the dominant parts are computable as the two respective quantities  $\displaystyle\frac{\omega}{2\pi\rho_\mathrm{m} \mathrm{c}_\mathrm{m} }\;\vert \mathrm{E}^{\textbf{in}} (\mathrm{z})\vert^2 \Big\vert\int_{\Omega}e^{(3)}_{\mathrm{n}_{0}}(\mathrm{x})d\mathrm{x} \Big\vert^2$ and $\displaystyle\frac{\omega}{2\pi\rho_\mathrm{m} \mathrm{c}_\mathrm{m} }\vert \mathrm{E}^{\textbf{in}} (\mathrm{z})\vert^2 \Big(\int_{\Omega}\mathrm{e}^{(\boldsymbol{\ell})}_{\mathrm{n}_{0}}(\mathrm{x})d\mathrm{x}\Big)^2$, for the plasmonics and the dielectrics, and they are at our disposal.
	\bigbreak
	\noindent
	It is worth mentioning that according to the expansions above, dielectric nanoparticles provide more pronounce amount of heat around its surrounding than the plasmonic ones. Therefore, the choice of the incident frequency $\omega$ makes a difference in generating the heat. Recall that, if it is chosen close to the undamped resonance frequency $\omega_0$, then the nanoparticle behaves as a dielectric one while if it is chosen via the electric plasma frequency $\omega_p$, then nanoparticle  behaves a a plasmonic one.
	\bigbreak
	\noindent
	The use of the heat generated by nanoparticles in different applications, as imaging, therapy etc, is well documented in the engineering literature, see \cite{A-B-B, B-C-R-2020, baffou, B-B-C, D-H-2013, G-S-K-S,  Maldovan} for instance. However, for our best knowledge, regarding its actual mathematical estimation, there is only the work \cite{ProfHabib} by Ammari et al. where the mathematical model was stated and the estimation where given on the surface of the nanoparticles. 
	To justify these expansions, in their work, they use semi-formal computations to characterize the dominant term. Precisely, they first use Fourier-Laplace transformation, then derive the small-volume approximations of in the Fourier-Laplace domain and then come back by inverse Fourier-Laplace transformation. However, in the Fourier-Laplace domain, the approximations of the fields are possible only at the 'low frequency' range and we have no access to the approximation in the full frequency range. In the current work, we use time-domain techniques, instead, to provide the actual dominant part of the heat. In addition to its mathematical interests, it allows us to deal with more general sources of heat and not only locally supported (in time) ones (as the initial or instantaneous sources for instance).
	\bigbreak
	\noindent
	The remaining part of the paper is organized as follows. In Section \ref{se2}, we state the equivalent boundary integral equations for the heat model after recalling the needed boundary integral operators and volume potentials. In Section \ref{se3}, we introduce the needed functional spaces i.e. the anisotropic Sobolev spaces, with their used properties. We also recall some basic mapping properties of the layer operators and volume potential for the heat equation on those spaces with their proper scales. The main steps of the proof of Theorem \ref{th2} are covered in Section \ref{sec4} while those of Theorem \ref{th1} and Theorem \ref{1.2} are provided in Section \ref{sec5} and Section \ref{dielectric} respectively. In order to simplify the exposition, some needed technical results and steps are moved to Appendix \ref{sec6}.
	
	\section{\Large\textbf{Representation Formula, Integral Operators and Boundary Integral Equations for the Heat Equation}}\label{se2}
	\noindent
	It is known that equations of boundary integral type are required to rephrase the transmission boundary value problem (\ref{eq:heat1}). With the use of both the interior and exterior Dirichlet problems between two media with different physical parameters, the boundary integral equation system can now be connected with the transmission condition. The following representation formula can therefore be used to determine the solution to the interior heat equation (\ref{eq:heat1}) using the direct approach, see \cite{costabel} and  \cite{dohr},
	\begin{align}\label{rep formula}
		\nonumber \mathrm{U}_{\mathrm{i}}(\mathrm{x},t) &=  \frac{1}{\alpha_{\mathrm{p}}}  \displaystyle\int_{-\infty}^{t}\int_{\partial\Omega} \bigg(\Phi(\mathrm{x},t;\mathrm{y},\tau)\gamma^{\textbf{int}}_{1}\mathrm{U}_{\mathrm{i}}(y,\tau) - \gamma^{\textbf{int}}_{0}\mathrm{U}_{\mathrm{i}}(y,\tau)\gamma^{\textbf{int}}_{1,\mathrm{y}}\Phi(\mathrm{x},t;\mathrm{y},\tau)\bigg)d\sigma_\mathrm{y}d\tau  \\ &+ \frac{1}{\alpha_\mathrm{p}}\int_{-\infty}^t\int_{\Omega}\Phi(\mathrm{x},\mathrm{t}; \mathrm{y},\mathrm{\tau})\frac{\omega \cdot \boldsymbol{\Im}(\varepsilon_\mathrm{p})}{2\pi\gamma_{\mathrm{p}}}|\mathrm{E}|^{2}(\mathrm{y})\rchi_{(0,T_0)} d\mathrm{y}d\tau.
	\end{align}
	\noindent
	Before proceeding, we need to define the following boundary integral operators $\text{for} \ (\mathrm{x},t) \in \partial \Omega\times \mathbb{R}$, namely the classical single-layer, double-layer and its spatial adjoint, and the hyper-singular heat operator as follows
	$$ \mathcal{S}\Big[\gamma^{\textbf{int}}_{1}\mathrm{U}_{\mathrm{i}}\Big](\mathrm{x},t) := \frac{1}{\alpha}  \displaystyle\int_\mathbb{R}\int_{\partial\Omega} \Phi(\mathrm{x},t;\mathrm{y},\tau)\gamma^{\textbf{int}}_{1}\mathrm{U}_{\mathrm{i}}(y,\tau) d\sigma_{\mathrm{y}}d\tau,
	$$
	$$ \mathcal{K}\Big[\gamma^{\textbf{int}}_{0}\mathrm{U}_{\mathrm{i}}\Big](\mathrm{x},t) := \frac{1}{\alpha}  \displaystyle\int_\mathbb{R}\int_{\partial\Omega} \gamma^{\textbf{int}}_{1,\mathrm{y}}\Phi(\mathrm{x},t;\mathrm{y},\tau)\gamma^{\textbf{int}}_{0}\mathrm{U}_{\mathrm{i}}(y,\tau) d\sigma_{\mathrm{y}}d\tau,
	$$
	$$ \mathcal{K^*}\Big[\gamma^{\textbf{int}}_{1}\mathrm{U}_{\mathrm{i}}\Big](\mathrm{x},t) := \frac{1}{\alpha}  \displaystyle\int_\mathbb{R}\int_{\partial\Omega} \gamma^{\textbf{int}}_{1,\mathrm{x}}\Phi(\mathrm{x},t;\mathrm{y},\tau)\gamma^{\textbf{int}}_{1}\mathrm{U}_{\mathrm{i}}(y,\tau) d\sigma_{\mathrm{y}}d\tau,
	$$
	$$
	\mathbb{H}\Big[\gamma^{\textbf{int}}_{0}\mathrm{U}_{\mathrm{i}}\Big](\mathrm{x},\mathrm{t}) := -\frac{1}{\alpha}\gamma^{\textbf{int}}_{1,\mathrm{x}}\displaystyle\int_\mathbb{R}\int_{\partial\Omega} \gamma^{\textbf{int}}_{1,\mathrm{y}}\Phi(\mathrm{x},t;\mathrm{y},\tau) \gamma^{\textbf{int}}_{0}\mathrm{U}_{\mathrm{i}}(y,\tau)d\sigma_\mathrm{y}d\tau
	$$
	respectively. Furthermore, we refer to the Newtonian heat potential associated with the source term $f := \frac{\omega \cdot  \boldsymbol{\Im}(\varepsilon_\mathrm{p})}{2\pi\gamma_{\mathrm{p}}} |\mathrm{E}|^{2}\rchi_{(0,\mathrm{T_0})}$ as
	$$ \mathcal{V}\Big[f\Big](\mathrm{x},t) :=  \displaystyle\int_{-\infty}^{t}\int_{\Omega}\Phi(\mathrm{x},t;\mathrm{y},\tau)f(\mathrm{y},\tau) d\mathrm{y}d\tau. 
	$$
	As a result, from the interior problem, we derive the following boundary integral equations:
	\begin{align}\label{intdrichlet1}
		\Big(\frac{1}{2}I_{d} + \mathcal{K}_{\alpha_{\mathrm{p}}}\Big)\Big[\gamma^{\textbf{int}}_{0}\mathrm{U}_{\mathrm{i}}\Big](\mathrm{x},t) = \mathcal{S}_{\alpha_{\mathrm{p}}}\Big[\gamma^{\textbf{int}}_{1}\mathrm{U}_{\mathrm{i}}\Big](\mathrm{x},t) + \gamma^{\textbf{int}}_{0}\mathcal{V}\Big[ f \Big](\mathrm{x},t), \ \text{for} \ \mathrm{x} \in \partial\Omega
	\end{align}
	and
	\begin{align}\label{intdrichlet2}
		\Big(\frac{1}{2}I_{d} - \mathcal{K}_{\alpha_{\mathrm{p}}}^{*}\Big)\Big[\gamma^{\textbf{int}}_{1}\mathrm{U}_{\mathrm{i}}\Big](\mathrm{x},t) =  \mathbb{H}_{\alpha_{\mathrm{p}}}\Big[\gamma^{\textbf{int}}_{0}\mathrm{U}_{\mathrm{i}}\Big](\mathrm{x},t) + \gamma^{\textbf{int}}_{1,\mathrm{x}}\mathcal{V}\Big[ f \Big](\mathrm{x},t), \ \text{for} \ \mathrm{x} \in \partial\Omega.
	\end{align}
	These equations together yield the Calderón system of boundary integral equations, \cite{dohr}
	\begin{align}\label{int}
		\begin{pmatrix}
			\gamma^{\textbf{int}}_{1}\mathrm{U}_{\mathrm{i}} \\
			\gamma^{\textbf{int}}_{0}\mathrm{U}_{\mathrm{i}}
		\end{pmatrix} = \underbrace{\begin{pmatrix}
				\frac{1}{2} - \mathcal{K}_{\alpha_{\mathrm{p}}} &  \mathcal{S}_{\alpha_{\mathrm{p}}} \\
				\mathbb{H}_{\alpha_{\mathrm{p}}} & \frac{1}{2} + \mathcal{K}_{\alpha_{\mathrm{p}}}^{*}
			\end{pmatrix}
		}_{=:\ \mathcal{C}}\begin{pmatrix}
			\gamma^{\textbf{int}}_{1}\mathrm{U}_{\mathrm{i}} \\
			\gamma^{\textbf{int}}_{0}\mathrm{U}_{\mathrm{i}}
		\end{pmatrix} +  \begin{pmatrix}
			\gamma^{\textbf{int}}_{0}\mathcal{V}[f] \\
			\gamma^{\textbf{int}}_{1,\mathrm{x}}\mathcal{V}[ f ]
		\end{pmatrix}.
	\end{align}
	The operator $\mathcal{C}$ is called the Calderón projection operator. A similar representation is also valid for the external heat problem in (\ref{eq:heat1})
	\begin{align}\label{ext}
		\begin{pmatrix}
			\gamma^{\textbf{ext}}_{1}\mathrm{U}_{\mathrm{i}} \\
			\gamma^{\textbf{ext}}_{0}\mathrm{U}_{\mathrm{i}}
		\end{pmatrix} = \underbrace{\begin{pmatrix}
				\frac{1}{2} + \mathcal{K}_{\alpha_{\mathrm{m}}} &  -\mathcal{S}_{\alpha_{\mathrm{m}}} \\
				-\mathbb{H}_{\alpha_{\mathrm{m}}} & \frac{1}{2} - \mathcal{K}_{\alpha_{\mathrm{m}}}^{*}
			\end{pmatrix}
		}_{=:\ \mathcal{C}}\begin{pmatrix}
			\gamma^{\textbf{ext}}_{1}\mathrm{U}_{\mathrm{i}} \\
			\gamma^{\textbf{ext}}_{0}\mathrm{U}_{\mathrm{i}}
		\end{pmatrix}
	\end{align}
	It is well known that $\mathcal{C}$ is a projection, i.e. $\mathcal{C} = \mathcal{C}^2,$ see \cite{costabel} and \cite{sao}. Therefore, we have the following identity
	\begin{align}\label{calderon}
		\mathcal{S}\mathbb{H} = (\frac{1}{2}I_{d} - \mathcal{K}) (\frac{1}{2}I_{d} + \mathcal{K}).
	\end{align}
	The above identity will be useful to scale the hyper-singular heat operator and hence the Steklov-Poincare operator $\mathbb{A}^{\textbf{ext}}$. Now together with the transmission condition in (\ref{eq:heat1}) and the boundary integral equation described in (\ref{int}), (\ref{ext}), we end up with the following system of boundary integral equation, 
	\begin{align}\label{IEfinal}
		\begin{bmatrix}
			\Big(\frac{1}{2}I_{d} - \mathcal{K}_{\alpha_{\mathrm{p}}}^{*}\Big)\Big[\gamma^{\textbf{int}}_{1}\mathrm{U}_{\mathrm{i}}\Big](\mathrm{x},t) = \mathbb{H}_{\alpha_\mathrm{p}}\Big[\gamma^{\textbf{int}}_{0}\mathrm{U}_{\mathrm{i}}\Big](\mathrm{x},t) + \gamma^{\textbf{int}}_{1,\mathrm{x}}\mathcal{V}\Big[ f \Big](\mathrm{x},t) \\
			\gamma^{\textbf{int}}_{0}\mathrm{U}_{\mathrm{i}}(\mathrm{x},t) = -\frac{\gamma_{\mathrm{m}}}{\gamma_{\mathrm{p}}} \Big(\frac{1}{2}I_{d} + \mathcal{K}_{\alpha_{\mathrm{p}}}\Big)^{-1}\mathcal{S}_{\alpha_{\mathrm{p}}} \mathbb{A}^{\textbf{ext}}\Big[\gamma^{\textbf{int}}_{0}\mathrm{U}_{\mathrm{i}}\Big](\mathrm{x},t) + \Big(\frac{1}{2}I_{d} + \mathcal{K}_{\alpha_{\mathrm{p}}}\Big)^{-1}\gamma^{\textbf{int}}_{0}\mathcal{V}\Big[ f\Big](\mathrm{x},t)
		\end{bmatrix},
	\end{align}
	associated with the source term $f := \frac{\omega \cdot  \boldsymbol{\Im}(\varepsilon_\mathrm{p})}{2\pi\gamma_{\mathrm{p}}} |\mathrm{E}|^{2}\rchi_{(0,\mathrm{T_0})}.$
	Further details on these derivations can be found in Appendix \ref{sec6}.
	
	\section{\Large\textbf{Function Spaces, Properties and Scales of Integral Operators }}\label{se3}
	We start this section with defining the anisotropic Sobolev spaces on space-time boundaries $\partial\Omega \times \mathbb{R}$ on which the above system of integral equations (\ref{IEfinal}) makes sense and can be inverted. 
	\subsection{Function Spaces}
	Let $\mathrm{H}^{\mathrm{r},\mathrm{s}}\Big(\partial \Omega \times \mathbb{R}\Big), \Big[\mathrm{r},\mathrm{s} \in (0,1)\Big]$, denote the anisotropic Sobolev spaces with respect to the following norms, see \cite[Chapter 4]{lions2} and \cite{Noonthe},
	\begin{equation}\label{def1}
		\Big\Vert u \Big\Vert_{\mathrm{H}^{\mathrm{r},\mathrm{s}}\Big(\partial\Omega \times \mathbb{R}\Big)}^2 := \Big\Vert u \Big\Vert_{\mathrm{L}^2\Big(\partial\Omega \times \mathbb{R}\Big)}^2 + \Big|u\Big|_{\mathrm{L}^2\Big(\mathbb{R};\mathrm{H}^{\mathrm{r}}(\partial\Omega)\Big)}^2 + \Big|u\Big|_{\mathrm{H}^{\mathrm{s}}\Big(\mathbb{R};\mathrm{L}^2(\partial\Omega)\Big)}^2,
	\end{equation}
	where
	\begin{equation}\label{def2}
		\Big|u\Big|_{\mathrm{L}^2\Big(\mathbb{R};\mathrm{H}^{\mathrm{r}}(\partial\Omega)\Big)}^2 := \int_{\partial\Omega}\int_{\partial\Omega}\dfrac{\Vert u(\mathrm{x},\cdot)-u(\mathrm{y},\cdot)\Vert_{\mathrm{L}^2(\mathbb{R})}^2}{|\mathrm{x}-\mathrm{y}|^{n-1+2\mathrm{r}}}d\sigma_\mathrm{x}d\sigma_\mathrm{y},
	\end{equation}
	and
	\begin{equation}\label{def3}
		\Big|u\Big|_{\mathrm{H}^{\mathrm{s}}\Big(\mathbb{R};\mathrm{L}^2(\partial\Omega)\Big)}^2 := \int_\mathbb{R}\int_\mathbb{R}\dfrac{\Vert u(\cdot,t)-u(\cdot,\tau)\Vert_{\mathrm{L}^2(\partial\Omega)}^2}{|t-\tau|^{1+2\mathrm{s}}}dt d\tau.
	\end{equation}
	Then, we define the negative-ordered anisotropic Sobolev space as the dual space of $\mathrm{H}^{-\mathrm{r},-\mathrm{s}}\Big(\partial\Omega \times \mathbb{R}\Big)$ when $\mathrm{r},\mathrm{s}<0$. We denote it by $\Big[\mathrm{H}^{-\mathrm{r},-\mathrm{s}}\Big(\partial\Omega \times \mathbb{R}\Big)\Big]'$ and equip it with the norm
	\begin{align}\label{def4}
		\Big\Vert u \Big\Vert_{\mathrm{H}^{\mathrm{r},\mathrm{s}}\Big(\partial\Omega \times \mathbb{R}\Big)} = \sup_{0\neq \varphi \in \mathrm{H}^{-\mathrm{r},-\mathrm{s}}\Big(\partial\Omega \times \mathbb{R}\Big)} \dfrac{\Big|\Big\langle u,\varphi \Big\rangle_{\partial \Omega\times \mathbb{R}}\Big|}{\Big\Vert \varphi \Big\Vert_{\mathrm{H}^{-\mathrm{r},-\mathrm{s}}\Big(\partial\Omega \times \mathbb{R}\Big)}}.
	\end{align}
	Furthermore, we introduce the following anisotropic Sobolev space $\mathrm{H}^{1,\frac{1}{2}}\Big(\partial\Omega \times \mathbb{R}\Big)$ with respect to the following norm
	\begin{align}\label{halfnorm}
		\Big\Vert u \Big\Vert_{\mathrm{H}^{1,\frac{1}{2}}\Big(\partial\Omega \times \mathbb{R}\Big)}^2 = \int_\mathbb{R}\int_{\partial\Omega} \Big[|\nabla_{\textbf{tan}}u(\mathrm{x},t)|^2 + u^2(\mathrm{x},t) + \Big(\partial_{t}^{\frac{1}{2}}u(\mathrm{x},t)\Big)^2\Big] d\sigma_{\mathrm{x}}d\mathrm{t}, 
	\end{align}
	where the tangential gradient $\nabla_{\textbf{tan}}u$ on $\partial\Omega$ is described as $\nabla_{\textbf{tan}}u := \nabla u -(\partial_{\nu}u)\nu$ and we denote by $\partial_{t}^{\frac{1}{2}}$ the fractional time derivative of order $\frac{1}{2}$. Next, due to the interpolation theory, we can state the following result regarding the anisotropic space $\mathrm{H}^{\frac{1}{2},\frac{1}{4}}(\Omega)$, see \cite[Chapter 2]{lions2} for detailed description. We have for $\mathrm{r},\mathrm{s}\ge 0$ and $\theta \in (0,1)$ 
	\begin{align}
		\Bigg[\mathrm{H}^{1,\frac{1}{2}}\Big(\partial \Omega\times \mathbb{R}\Big), \mathrm{L}^2\Big((\partial\Omega)_T\Big)\Bigg]_{\theta = \frac{1}{2}} = \mathrm{H}^{\frac{1}{2},\frac{1}{4}}\Big((\partial\Omega)_T\Big).
	\end{align}
	In addition, we point out the following interpolation inequality
	\begin{align}
		\Big\Vert u\Big\Vert_{\mathrm{H}^{\frac{1}{2},\frac{1}{4}}\Big((\partial\Omega)_T\Big)} \lesssim \Big\Vert u\Big\Vert^{\frac{1}{2}}_{\mathrm{H}^{1,\frac{1}{2}}\Big((\partial\Omega)_T\Big)} \Big\Vert u \Big\Vert^{\frac{1}{2}}_{\mathrm{L}^2\Big((\partial\Omega)_T\Big)}.
	\end{align}
	In the following paragraphs, we review the mapping properties of the boundary integral operators mentioned above. It has been well established that all boundary integral operators exhibit a mapping property, see, for instance \cite{baffou}, \cite{costabel}, and \cite{fabes}. In our analysis, we need the $\mathcal{C}^2$-regularity of the boundary $\partial\Omega$ to study the magnetization operator, otherwise only the Lipschitz regularity is enough. 

	\subsection{Mapping Properties of the Integral Operators}
	\begin{lemma}\label{3.1}\cite[Section 5]{Noonthe}
		On smooth $\mathcal{C}^2$-boundary $\partial\Omega$ and for all $\mathrm{r}\ge 0$ the single layer heat operator $\mathcal{S}$ maps $H^{-\frac{1}{2}+\mathrm{r},(-\frac{1}{2}+\mathrm{r})/2}\Big(\partial \Omega \times \mathbb{R}\Big)$ to $ H^{\frac{1}{2}+\mathrm{r},(\frac{1}{2}+\mathrm{r})/2}\Big(\partial \Omega \times \mathbb{R}\Big)$ isomorphically.
	\end{lemma} 
	\noindent
	This lemma still holds for a Lipschitz boundary, if $\mathrm{r}$ is in a certain finite range. Next, we state the following result.
	\begin{lemma}\label{3.2}\cite{hofman} 
		Let $\Omega$ be the domain above the graph of a Lipschitz function $\varphi:\mathbb{R}^{\mathrm{n}-1} \rightarrow \mathbb{R}\ (\mathrm{n}\ge 2).$ Then the operator
		\begin{align}
			\frac{1}{2} I + \mathcal{K}^* : L^\mathrm{2}\Big(\partial\Omega \times \mathbb{R}\Big) \rightarrow L^\mathrm{2}\Big(\partial\Omega \times \mathbb{R}\Big)
		\end{align}
		is invertible.
		\vskip 0.1in
		\noindent
		Moreover, the operators
		\begin{align}
			\frac{1}{2} I + \mathcal{K} : \mathrm{L}^2\Big(\partial\Omega \times \mathbb{R}\Big) \rightarrow \mathrm{L}^2\Big(\partial\Omega \times \mathbb{R}\Big),
		\end{align}
		\begin{align}
			\frac{1}{2} I + \mathcal{K} : H^{1,\frac{1}{2}}\Big(\partial\Omega \times \mathbb{R}\Big) \rightarrow H^{1,\frac{1}{2}}\Big(\partial\Omega \times \mathbb{R}\Big)
		\end{align}
		are invertible.
	\end{lemma} 
	\noindent
	Consequently, with the help of interpolation theory, we can demonstrate the following corollary to the above lemma:
	\begin{corollary}
		The operators
		\begin{align}
			\frac{1}{2} I + \mathcal{K} : \mathrm{H}^{\frac{1}{2},\frac{1}{4}}\Big(\partial\Omega \times \mathbb{R}\Big) \rightarrow \mathrm{H}^{\frac{1}{2},\frac{1}{4}}\Big(\partial\Omega \times \mathbb{R}\Big),
		\end{align}
		\begin{align}
			\frac{1}{2} I + \mathcal{K}^* : \mathrm{H}^{-\frac{1}{2},-\frac{1}{4}}\Big(\partial\Omega \times \mathbb{R}\Big) \rightarrow \mathrm{H}^{-\frac{1}{2},-\frac{1}{4}}\Big(\partial\Omega \times \mathbb{R}\Big)
		\end{align}
		are invertible.
	\end{corollary}
	Then, we state the mapping property of Newtonian heat potential, see for instance \cite{costabel}.
	\begin{lemma}\label{3.3}
		The convolution between the fundamental solution $\Phi(\mathrm{x},\mathrm{t};\mathrm{y},\mathrm{\tau})$ and density i.e. the Newtonian heat potential
		\begin{align}
			\mathcal{V} : H^{r,\frac{r}{2}}_{\text{comp}}\Big(\mathbb{R}^{\mathrm{n}}\times \mathbb{R}\Big) \rightarrow H^{r+2,\frac{r}{2}+1}_{\text{loc}}\Big(\mathbb{R}^{\mathrm{n}}\times \mathbb{R}\Big)
		\end{align}
		is linear and bounded for any $\mathrm{r} \in \mathbb{R}.$ 
		\vskip 0.1in
		\noindent
		Here, we denote by $H^{\mathrm{r},\frac{\mathrm{r}}{2}}_{\text{comp}}\Big(\mathbb{R}^{\mathrm{n}}\times (0,\mathrm{T})\Big)$ set of functions which has compact support in space variable and by $H^{\mathrm{r}+2,\frac{\mathrm{r}}{2}+1}_{\text{loc}}\Big(\mathbb{R}^{\mathrm{n}}\times (0,\mathrm{T})\Big)$, we describe the local behaviour in the space variable. Now, it is clear that for $\mathrm{r}=-1$ and by restriction of the domain, we find the Newtonian heat potential 
		\begin{align}
			\mathcal{V} : \Big[H^{1,\frac{1}{2}}\Big(\Omega \times \mathbb{R}\Big)\Big]' \rightarrow \Big[H^{1,\frac{1}{2}}\Big(\Omega \times \mathbb{R}\Big)\Big]
		\end{align}
		is linear and bounded.
	\end{lemma}
	\noindent
	The next step is to state the trace theorem. See, for example, \cite[Theorem 2.1]{lions2}, \cite[Chapter 2]{Noonthe}.
	\begin{lemma}\label{3.4}
		(Trace Theorem) The interior trace operator 
		\begin{align}
			\gamma^{\textbf{int}}_{0}:  H^{1,\frac{1}{2}}\Big(\Omega \times \mathbb{R}\Big) \rightarrow \mathrm{H}^{\frac{1}{2},\frac{1}{4}}\Big(\partial\Omega \times \mathbb{R}\Big) 
		\end{align}
		is linear and bounded.
	\end{lemma}
	\noindent
	Now, regarding the Newtonian heat potential, we restrict the domain to $\mathrm{L}^2(\Omega)$ as we consider the corresponding source term belonging to that space. Consequently, from Lemma \ref{3.3} and trace theorem we deduce the following corollary
	\begin{corollary}\label{cc3}
		The following operator
		\begin{align}
			\gamma^{\textbf{int}}_{0}\mathcal{V} : \mathrm{L}^2(\Omega \times \mathbb{R}) \rightarrow \mathrm{H}^{\frac{1}{2},\frac{1}{4}}\Big(\partial \Omega \times \mathbb{R}\Big)
		\end{align}
		is linear and continuous.
	\end{corollary}
	As we proceed, we introduce an anisotropic Sobolev space that is required to describe the mapping properties of the interior and exterior Neumann trace operators, as well as their application to the Newtonian heat operator.
	\begin{align}
		H^{1,\frac{1}{2}}\Big(\Omega \times \mathbb{R},\mathcal{L}\Big) := \Big\{ \mathrm{u} \in H^{1,\frac{1}{2}}\Big(\Omega \times \mathbb{R}\Big) : \mathcal{L}u \in \mathrm{L}^2\Big(\Omega \times \mathbb{R}\Big)\Big\},
	\end{align}
	where, $\mathcal{L} := \alpha\partial_{t}-\Delta$ represents the corresponding heat differential operator. Now we will state the following lemma, see \cite[Proposition 2.18]{costabel}, for details.
	\begin{lemma}\label{3.5}
		The interior Neumann Trace Operator
		\begin{align}
			\gamma^{\textbf{int}}_{1} : H^{1,\frac{1}{2}}\Big(\Omega \times \mathbb{R},\mathcal{L}\Big) \rightarrow \mathrm{H}^{-\frac{1}{2},-\frac{1}{4}}\Big(\partial \Omega \times \mathbb{R}\Big) 
		\end{align}
		defines a linear bounded operator.
	\end{lemma}
	\noindent
	Furthermore, from Lemma \ref{3.3}, after restricting the domain to $\mathrm{L}^2(\Omega)$ and considering that the Newtonian heat potential satisfies the corresponding heat equation, we obtain the following corollary
	\begin{corollary}\label{c4}
		The following operator
		\begin{align}
			\gamma^{\textbf{int}}_{1}\mathcal{V} : \mathrm{L}^2(\Omega \times \mathbb{R}) \rightarrow \mathrm{H}^{-\frac{1}{2},-\frac{1}{4}}\Big(\partial \Omega \times \mathbb{R}\Big)
		\end{align}
		is linear and continuous.
	\end{corollary}
	Also, let us define the initial heat potential for $f \in \mathrm{L}^2\Big(\Omega \Big)$ as follows
	\begin{align}
		\mathbb{I}[f](\mathrm{x},\mathrm{t}) = \int_{\Omega}\Phi(\mathrm{x},\mathrm{t};\mathrm{y})f(\mathrm{y})d\mathrm{y},
	\end{align}
	which enjoys the following mapping property. For more details we refer to \cite[Theorem 4.4]{dohr} and \cite[Section 7]{Noonthe}.
	\begin{lemma}\label{4.3}
		The initial heat operator $\mathbb{I}: \mathrm{L}^2\Big(\mathrm{B} \Big) \rightarrow \mathrm{H}^{1,\frac{1}{2}}\Big(\mathrm{B} \times \mathbb{R}_+\Big) $ is linear and bounded.
	\end{lemma}
	\subsection{Scales for the Function Spaces}
	Using the same notation as in \cite{prof sini}, we consider a nanoparticle occupying a domain $\Omega = \delta \mathrm{B} + \mathrm{z}$, where $\mathrm{B}$ is centered at the origin and $|\mathrm{B}| \sim 1.$ Moreover, we use the notation below in defining functions  $\varphi$ and $\psi$ on $\partial\Omega\times \mathbb{R}$ and $\partial \mathrm{B}\times \mathbb{R}$, respectively,
	\begin{equation}
		\hat{\varphi}( \eta, \Tilde{\tau}) = \varphi^{\Lambda}( \eta, \Tilde{\tau}) := \varphi(\delta\eta + \mathrm{z}, \alpha\delta^2\Tilde{\tau}), \quad \quad \quad  \check{\psi}(\mathrm{x}, t) = \psi^\vee(\mathrm{x}, t) := \psi\Bigg(\frac{\mathrm{x}-\mathrm{z}}{\delta}, \frac{t}{\alpha\delta^2}\Bigg)    
	\end{equation}
	for $(\mathrm{x},t) \in \partial\Omega \times\mathbb{R} $ and $(\eta,\Tilde{\tau}) \in \partial B\times \mathbb{R}$ respectively. Our next step is to derive the following lemma.
	\begin{lemma}\label{3.6}
		Suppose $0<\delta\leq 1$, $\Omega:=\delta B + \mathrm{z}$ and $t := \alpha \delta^2\Tilde{t} $. Then for $\varphi \in \mathrm{H}^{\frac{1}{2},\frac{1}{4}}\Big(\partial\Omega\times \mathbb{R}\Big) $  and $\psi \in \mathrm{H}^{-\frac{1}{2},-\frac{1}{4}}\Big(\partial\Omega\times \mathbb{R}\Big)$, we have the following scales
		\begin{align}\label{l1}
			\alpha \delta^{3}\Big\Vert\hat{\varphi}\Big\Vert_{\mathrm{H}^{\frac{1}{2},\frac{1}{4}}\Big(\partial B\times \mathbb{R}\Big)}^2 \leq \Big\Vert\varphi\Big\Vert_{\mathrm{H}^{\frac{1}{2},\frac{1}{4}}\Big(\partial\Omega\times \mathbb{R}\Big)}^2 \leq \alpha^{\frac{1}{2}}\delta^2\Big\Vert\hat{\varphi}\Big\Vert_{\mathrm{H}^{\frac{1}{2},\frac{1}{4}}\Big(\partial B\times \mathbb{R}\Big)}^2
		\end{align}
		and
		\begin{align}\label{l2}
			\alpha^{\frac{3}{4}}\delta^2\Big\Vert\hat{\psi}\Big\Vert_{\mathrm{H}^{-\frac{1}{2},-\frac{1}{4}}\Big(\partial B\times \mathbb{R}\Big) }\leq \Big\Vert\varphi\Big\Vert_{\mathrm{H}^{-\frac{1}{2},-\frac{1}{4}}\Big(\partial\Omega\times \mathbb{R}\Big)} \leq \alpha^{\frac{1}{2}}\delta^{\frac{3}{2}}\Big\Vert\hat{\psi}\Big\Vert_{\mathrm{H}^{-\frac{1}{2},-\frac{1}{4}}\Big(\partial B\times \mathbb{R}\Big)}.
		\end{align}
	\end{lemma}
	\textbf{Proof.} 
	Let us set $\mathrm{x} := \delta \xi + \mathrm{z}$, $\mathrm{y}:= \delta \eta + \mathrm{z}$, $t:= \alpha \delta^2 \Tilde{t}$ and $\tau = \alpha \delta^2 \Tilde{\tau}$. As a first step, we need to scale $\Big\Vert \varphi \Big\Vert_{\mathrm{L}^2\Big(\partial \Omega \times \mathbb{R}\Big)}^2$. Indeed,
	\begin{align}\label{e1}
		\Big\Vert \varphi \Big\Vert_{\mathrm{L}^2\Big(\partial \Omega \times \mathbb{R}\Big)}^2 &= \int_{-\infty}^T\int_{\partial \Omega}|\varphi|^2(\mathrm{x},t) d\sigma_\mathrm{x} dt =   \alpha \delta^3 \int_{-\infty}^{T_{\delta}}\int_{\partial B}|\varphi|^2(\delta \eta + \mathrm{z},\alpha \delta^2 \Tilde{\tau}) d\sigma_\eta d\Tilde{\tau}
		= \alpha \delta^3 \Big\Vert \hat{\varphi} \Big\Vert_{\mathrm{L}^2\Big(\partial B \times \mathbb{R}\Big)}^2
	\end{align}
	Now, in order to scale the norms (\ref{def2}) and (\ref{def3}), we need to estimate $\Big\Vert \varphi \Big\Vert_{\mathrm{L}^2(\mathbb{R})}^2$ and $ \Big\Vert \varphi \Big\Vert_{\mathrm{L}^2(\partial \Omega)}^2$ respectively.
	We have
	\begin{align}\label{e2}
		\Big\Vert \varphi \Big\Vert_{\mathrm{L}^2(\mathbb{R})}^2 = \int_{-\infty}^T|\varphi|^2(\mathrm{x},t) dt =   \alpha \delta^2 \int_{-\infty}^{T_{\delta}}|\varphi|^2(\delta \eta + \mathrm{z},\alpha \delta^2 \Tilde{\tau}) d\Tilde{\tau} = \alpha \delta^2 \Big\Vert \hat{\varphi} \Big\Vert_{\mathrm{L}^2(\mathbb{R})}^2
	\end{align}
	and
	\begin{align}\label{e3}
		\Big\Vert \varphi \Big\Vert_{\mathrm{L}^2(\partial \Omega)}^2 = \int_{\partial \Omega}|\varphi|^2(\mathrm{x},t)  d\sigma(\mathrm{x}) = \delta \int_{\partial B}|\varphi|^2(\delta \eta + \mathrm{z},\alpha \delta^2 \Tilde{\tau})  d\sigma_\eta = \delta \Big\Vert \hat{\varphi} \Big\Vert_{\mathrm{L}^2(\partial B)}^2.
	\end{align}
	Therefore, using (\ref{e2}) and (\ref{e3}), we estimate the following norms
	\begin{align}\label{a1}
		\nonumber
		\Big|\varphi\Big|_{\mathrm{L}^2\Big(\mathbb{R};H^{\frac{1}{2}}(\partial\Omega)\Big)}^2 &= \int_{\partial\Omega}\int_{\partial\Omega}\dfrac{\Vert \varphi(x,\cdot) - \varphi(y,\cdot)\Vert_{\mathrm{L}^2(\mathbb{R})}^2}{|x-y|^{2}}d\sigma_\mathrm{x}d\sigma_\mathrm{y}\\ \nonumber &= \delta^2 \int_{\partial B}\int_{\partial B}\dfrac{\alpha \delta^2\Vert \hat{\varphi}(\xi,\cdot)-\hat{\varphi}(\eta,\cdot)\Vert_{\mathrm{L}^2(\mathbb{R})}^2}{\delta^2|\xi-\eta|^{2}}d\sigma_\xi d\sigma_\eta \\ &= \alpha \delta^2 \Big|\hat{\varphi}\Big|_{\mathrm{L}^2\Big(\mathbb{R};H^{\frac{1}{2}}(\partial B)\Big)}^2.
	\end{align}
	Again, we have
	\begin{align}\label{a2}
		\nonumber
		\Big|\varphi\Big|_{H^{\frac{1}{4}}\Big(\mathbb{R};\mathrm{L}^2(\partial\Omega)\Big)}^2 &= \int_{-\infty}^T\int_{-\infty}^T\dfrac{\Vert \varphi(\cdot,t)-\varphi(\cdot,\tau)\Vert_{\mathrm{L}^2(\partial \Omega)}^2}{|t-\tau|^{\frac{3}{2}}}dtd\tau\\ \nonumber &= \alpha^2\delta^4 \int_{-\infty}^{T_{\delta}}\int_{-\infty}^{T_{\delta}}\dfrac{ \delta\Vert \varphi(\cdot, \Tilde{t})-\varphi(\cdot, \Tilde{\tau})\Vert_{\mathrm{L}^2(\partial B)}^2}{\alpha^{\frac{3}{2}}\delta^3|\Tilde{t}-\Tilde{\tau}|^{\frac{3}{2}}}d\Tilde{t}d\Tilde{\tau} \\ &= \alpha^{\frac{1}{2}} \delta^2 \Big|\hat{\varphi}\Big|_{H^{\frac{1}{4}}\Big(\mathbb{R};\mathrm{L}^2(\partial B)\Big)}^2 
	\end{align}
	Finally, considering the norm (\ref{e1}), (\ref{a1}), and (\ref{a2}), noticing the fact {\color{Black}{$\alpha\delta^3 \leq \alpha\delta^2 \leq \alpha^{\frac{1}{2}} \delta^2$}}, we prove the first double inequality,
	\begin{align}\label{l1}
		\alpha\delta^{3}\Big\Vert \hat{\varphi}\Big\Vert_{\mathrm{H}^{\frac{1}{2},\frac{1}{4}}\Big(\partial B \times \mathbb{R}\Big)}^2 \leq \Big\Vert \varphi\Big\Vert_{\mathrm{H}^{\frac{1}{2},\frac{1}{4}}\Big(\partial \Omega \times \mathbb{R}\Big)}^2 \leq  \alpha^{\frac{1}{2}} \delta^2\Big\Vert \hat{\varphi}\Big\Vert_{\mathrm{H}^{\frac{1}{2},\frac{1}{4}}\Big(\partial B \times \mathbb{R}\Big)}^2.
	\end{align}
	In order to derive the second inequality, we need to scale $\Big\langle u,\psi \Big\rangle_{\partial \Omega \times \mathbb{R}}$.
	\begin{align}
		\nonumber
		\Big\langle u,\psi \Big\rangle_{\partial \Omega \times \mathbb{R}} &= \int_{-\infty}^{T}\int_{\partial \Omega} u(y,\tau)\psi(y,\tau)d\mathrm{y}d\tau\\ \nonumber &= \alpha \delta^3 \int_{-\infty}^{T_{\delta}}\int_{\partial B} u(\delta \eta + \mathrm{z},\alpha \delta^2 \Tilde{\tau})\psi(\delta \eta + \mathrm{z},\alpha \delta^2 \Tilde{\tau})d\sigma_\eta d\Tilde{\tau} \\ &= \alpha \delta^3 \Big\langle \hat{u},\hat{\psi} \Big\rangle_{\partial B \times \mathbb{R}}
	\end{align}
	Hence, from the definition of norm (\ref{def4}), we get our desired second inequality
	\begin{align}\label{le2}
		\alpha^{\frac{3}{4}}\delta^2\Big\Vert \hat{\psi}\Big\Vert_{\mathrm{H}^{-\frac{1}{2},-\frac{1}{4}}\Big(\partial B \times \mathbb{R}\Big)} \leq \Big\Vert \psi\Big\Vert_{\mathrm{H}^{-\frac{1}{2},-\frac{1}{4}}\Big(\partial \Omega \times \mathbb{R}\Big)} \leq  \alpha^{\frac{1}{2}}\delta^{\frac{3}{2}}\Big\Vert \hat{\psi}\Big\Vert_{\mathrm{H}^{-\frac{1}{2},-\frac{1}{4}}\Big(\partial B \times \mathbb{R}\Big)}.
	\end{align}
	Hence it completes the proof.
	\begin{lemma}
		For $\partial_t\varphi \in \mathrm{H}^{-\frac{1}{2},-\frac{1}{4}}\Big(\partial \Omega\times\mathbb{R}\Big)$ we have the following scales
		\begin{align}\label{le2}
			\alpha^{-\frac{1}{4}}\Big\Vert \partial_{\Tilde{t}}\hat{\varphi}\Big\Vert_{\mathrm{H}^{-\frac{1}{2},-\frac{1}{4}}\Big(\partial B\times \mathbb{R}\Big)} \leq \Big\Vert \partial_t\varphi\Big\Vert_{\mathrm{H}^{-\frac{1}{2},-\frac{1}{4}}\Big(\partial \Omega\times \mathbb{R}\Big)} \leq  \alpha^{-\frac{1}{2}}\delta^{-\frac{1}{2}}\Big\Vert \partial_{\Tilde{t}}\hat{\varphi}\Big\Vert_{\mathrm{H}^{-\frac{1}{2},-\frac{1}{4}}\Big(\partial B\times \mathbb{R}\Big)}.
		\end{align}
	\end{lemma}
	\noindent
	\textbf{Proof.} 
	As a first step, we will do the scaling of  $\Big\Vert \partial_t\varphi \Big\Vert_{\mathrm{L}^2\Big(\partial \Omega\times\mathbb{R}\Big)}^2$. Indeed,
	\begin{align}\label{e1}
		\nonumber
		\Big\Vert \partial_t\varphi \Big\Vert_{\mathrm{L}^2\Big(\partial \Omega \times \mathbb{R}\Big)}^2 &= \int_{-\infty}^T\int_{\partial \Omega}|\partial_t\varphi(\mathrm{x},t)|^2 d\sigma_\mathrm{x} dt =   \alpha \delta^3 \int_{-\infty}^{T_{\delta}}\int_{\partial B}\alpha^{-2}\delta^{-4}|\partial_{\Tilde{t}}\varphi(\delta \eta + \mathrm{z},\alpha \delta^2 \Tilde{\tau})|^2 d\sigma_\eta d\Tilde{\tau}
		\\ &= \frac{1}{\alpha\delta} \Big\Vert \partial_{\Tilde{t}}\hat{\varphi} \Big\Vert_{\mathrm{L}^2\Big(\partial B \times \mathbb{R}\Big)}^2.
	\end{align}
	Next, we calculate the following product.
	\begin{align}
		\nonumber
		\Big\langle u, \partial_t\varphi \Big\rangle_{\partial \Omega\times \mathbb{R}} &= \int_{0}^T\int_{\partial \Omega}u(\mathrm{y},\tau)\partial_t\varphi(\mathrm{x},t) d\sigma_\mathrm{y} d\tau \\ \nonumber&= \alpha \delta^3 \int_{0}^{T_{\delta}}\int_{\partial B}\alpha^{-1}\delta^{-2}\partial_{\Tilde{t}}\varphi(\delta \eta + \mathrm{z},\alpha \delta^2 \Tilde{\tau})u(\delta \eta + \mathrm{z},\alpha \delta^2 \Tilde{\tau}) d\sigma_\eta d\Tilde{\tau} \\ \nonumber
		&= \delta \int_{0}^{T_{\delta}}\int_{\partial B}\partial_{\Tilde{t}}\hat{\varphi}(\eta,\Tilde{\tau})\hat{u}(\eta, \Tilde{\tau}) d\sigma_\eta d\Tilde{\tau}
		\\ &= \delta \Big\langle \hat{u}, \partial_{\Tilde{t}}\hat{\varphi} \Big\rangle_{\partial B\times\mathbb{R}}
	\end{align}
	Then by Lemma \ref{3.6} and from the previous estimate we obtain the desired inequality
	\begin{align}\label{le2}
		\alpha^{-\frac{1}{4}}\Big\Vert \partial_{\Tilde{t}}\hat{\varphi}\Big\Vert_{\mathrm{H}^{-\frac{1}{2},-\frac{1}{4}}\Big(\partial B \times \mathbb{R}\Big)} \leq \Big\Vert \partial_t\varphi\Big\Vert_{\mathrm{H}^{-\frac{1}{2},-\frac{1}{4}}\Big(\partial \Omega \times \mathbb{R}\Big)} \leq  \alpha^{-\frac{1}{2}}\delta^{-\frac{1}{2}}\Big\Vert \partial_{\Tilde{t}}\hat{\varphi}\Big\Vert_{\mathrm{H}^{-\frac{1}{2},-\frac{1}{4}}\Big(\partial B \times \mathbb{R}\Big)}.
	\end{align}
	\begin{lemma}
		For $\varphi \in H^{1,\frac{1}{2}}\Big(\Omega \times\mathbb{R}_+\Big) $ and $\psi \in H^{-1,-\frac{1}{2}}\Big(\Omega\times\mathbb{R}_+\Big) $ we have the following scales
		\begin{align}\label{l1}
			\alpha^{\frac{1}{2}}\delta^{2}\Big\Vert \hat{\varphi}\Big\Vert_{\mathrm{H}^{1,\frac{1}{2}}\Big(B\times\mathbb{R}_+\Big)} \leq \Big\Vert \varphi\Big\Vert_{\mathrm{H}^{1,\frac{1}{2}}\Big(\Omega\times\mathbb{R}_+\Big)} \leq \delta\Big\Vert \hat{\varphi}\Big\Vert_{\mathrm{H}^{1,\frac{1}{2}}\Big( B\times\mathbb{R}_+\Big)}
		\end{align}
		and
		\begin{align}\label{le2}
			\alpha\delta^3\Big\Vert \hat{\psi}\Big\Vert_{H^{-1,-\frac{1}{2}}\Big( B\times\mathbb{R}_+\Big)} \leq \Big\Vert \psi\Big\Vert_{H^{-1,-\frac{1}{2}}\Big(\Omega\times\mathbb{R}_+\Big)} \leq  \alpha^{\frac{1}{2}}\delta^2\Big\Vert \hat{\psi}\Big\Vert_{H^{-1,-\frac{1}{2}}\Big(B\times\mathbb{R}_+\Big)}.
		\end{align}
	\end{lemma}
	\textbf{Proof.} For $\varphi \in H^{1,\frac{1}{2}}\Big(\Omega\times\mathbb{R_+}\Big)$ we have
	\begin{align}
		\nonumber
		\Big\Vert \psi\Big\Vert_{\mathrm{H}^{1,\frac{1}{2}}\Big( \Omega\times\mathbb{R_+}\Big)}^2 &= \int_{0}^T\int_{\Omega}\Big[|\nabla_{\textbf{tan}}\varphi(\mathrm{x},t)|^2 + |\varphi(\mathrm{x},t)|^2 + \Big(\partial_{t}^{\frac{1}{2}}\varphi(\mathrm{x},t)\Big)^2\Big] d\mathrm{x}dt
		\\ \nonumber &= \alpha\delta^4\int_{0}^{T_\delta}\int_{B}\Big[\delta^{-2}|\nabla_{\textbf{tan}}\varphi(\delta\xi+\mathrm{z},\alpha\delta^2\Tilde{t})|^2 + |\varphi(\delta\xi+\mathrm{z},\alpha\delta^2\Tilde{t})|^2 \Big] d\xi d\Tilde{t}
		\\ &+ \int_{0}^{T_\delta}\int_{B}\alpha^{-1}\delta^{-2}\Big(\partial_{t}^{\frac{1}{2}}\varphi(\delta\xi+\mathrm{z},\alpha\delta^2\Tilde{t})\Big)^2 d\xi d\Tilde{t}
	\end{align}
	Now, observing the fact that $\alpha\delta^4 \le \delta^2$, we obtain the desired first inequality
	\begin{align}\label{l1}
		\alpha^{\frac{1}{2}}\delta^{2}\Big\Vert \hat{\varphi}\Big\Vert_{\mathrm{H}^{1,\frac{1}{2}}\Big(B\times\mathbb{R_+}\Big)} \leq \Big\Vert \varphi\Big\Vert_{\mathrm{H}^{1,\frac{1}{2}}\Big(\Omega\times\mathbb{R_+}\Big)} \leq \delta\Big\Vert \hat{\varphi}\Big\Vert_{\mathrm{H}^{1,\frac{1}{2}}\Big(B\times\mathbb{R_+}\Big)}.
	\end{align}
	Next, from the definition of norm for the dual space we know that
	\begin{align}\label{def4}
		\Big\Vert \psi \Big\Vert_{H^{-1,-\frac{1}{2}}\Big(\Omega\times\mathbb{R_+}\Big)} = \sup_{0\neq u \in H^{1,\frac{1}{2}}\Big( \Omega\times\mathbb{R_+}\Big)} \dfrac{\Big|\Big\langle \psi,u \Big\rangle_{ \Omega\times\mathbb{R_+}}\Big|}{\Big\Vert u \Big\Vert_{\mathrm{H}^{1,\frac{1}{2}}\Big(\Omega\times\mathbb{R_+}\Big)}}.
	\end{align}
	Now, to get the desired inequality we need to scale the scalar product $\Big\langle \psi,u \Big\rangle_{(\partial \Omega)_{T}}$.
	\begin{align}
		\nonumber
		\Big\langle \psi,u \Big\rangle_{\Omega\times\mathbb{R_+}} &= \int_{0}^{T}\int_{\Omega} u(\mathrm{y},\tau)\psi(\mathrm{y},\tau)d\mathrm{y}d\tau \\ \nonumber &= \alpha \delta^4 \int_{0}^{T_{\delta}}\int_{B} \hat{u}(\eta,\Tilde{\tau})\hat{\psi}(\eta,\Tilde{\tau})d\eta d\Tilde{\tau}\\ &= \alpha \delta^4 \Big\langle \hat{\psi},\hat{u}\Big\rangle_{B\times\mathbb{R_+}}
	\end{align}
	Consequently, this leads to
	\begin{align}\label{le2}
		\alpha\delta^3\Big\Vert \hat{\psi}\Big\Vert_{H^{-1,-\frac{1}{2}}\Big(B\times\mathbb{R_+}\Big)} \leq \Big\Vert \psi\Big\Vert_{H^{-1,-\frac{1}{2}}\Big( \Omega\times\mathbb{R_+}\Big)} \leq  \alpha^{\frac{1}{2}}\delta^2\Big\Vert \hat{\psi}\Big\Vert_{H^{-1,-\frac{1}{2}}\Big(B\times\mathbb{R_+}\Big)}.
	\end{align}
	\bigbreak
	\noindent
	The next step is to state and derive the following lemmas, which will be useful in estimating the hyper-singular and Steklov-Poincare operators, as well as providing an exact dominating term for the estimated heat potential.
	\subsection{Scales for the Integral Operators}
	\begin{lemma}\label{3.7}
		For $\varphi \in \mathrm{H}^{\frac{1}{2},\frac{1}{4}}\Big(\partial \Omega\times \mathbb{R}\Big)$ and $\psi \in \mathrm{H}^{-\frac{1}{2},-\frac{1}{4}}\Big(\partial \Omega\times \mathbb{R}\Big)$, we have the following estimate
		\begin{align}
			\mathcal{S}_{\partial \Omega\times \mathbb{R}}\Big[\psi\Big](\mathrm{x},t) =  \delta \Big(\overline{\mathcal{S}}_{\partial B\times \mathbb{R}}\Big[\hat{\psi}\Big]\Big)^\vee \ \text{and} \
			\mathcal{S}^{-1}_{\partial \Omega\times \mathbb{R}}\Big[\varphi\Big](\mathrm{x},t) =  \delta^{-1} \Big(\overline{\mathcal{S}}^{-1}_{\partial B\times \mathbb{R}}\Big[\hat{\varphi}\Big]\Big)^\vee.
		\end{align}
		Based on the above estimates, we obtain the following one
		\begin{align}\label{operatoresti}
			\Big\Vert \mathcal{S}^{-1}_{\partial \Omega\times \mathbb{R}} \Big\Vert_{\mathcal{L}{\Big( \mathrm{H}^{\frac{1}{2},\frac{1}{4}}\Big(\partial \Omega\times \mathbb{R}\Big), \mathrm{H}^{-\frac{1}{2},-\frac{1}{4}}\Big(\partial \Omega\times \mathbb{R}\Big)}\Big) } \le \delta^{-1} \Big\Vert \overline{\mathcal{S}}^{-1}_{\partial B\times \mathbb{R}} \Big\Vert_{\mathcal{L}{\Big( \mathrm{H}^{\frac{1}{2},\frac{1}{4}}\Big(\partial B\times \mathbb{R}\Big), \mathrm{H}^{-\frac{1}{2},-\frac{1}{4}}\Big(\partial B\times \mathbb{R}}\Big)\Big)},
		\end{align}
		where, $\overline{\mathcal{S}}$ represents the single layer operator corresponding to the fundamental solution with $\alpha:=1$.
	\end{lemma}
	\noindent
	\textbf{Proof.} As usual, we set $\mathrm{x} := \delta \xi + \mathrm{z}$, $\mathrm{y}:= \delta \eta + \mathrm{z}$, $t := \alpha \delta^2 \Tilde{t}$ and $\tau := \alpha \delta^2 \Tilde{\tau}$. We write 
	\begin{align}
		\nonumber
		\mathcal{S}_{\partial \Omega\times \mathbb{R}}\Big[\psi\Big](\mathrm{x},t) &=  \frac{1}{\alpha}\displaystyle\int_\mathbb{R}\int_{\partial\Omega} \Phi(\mathrm{x},t;\mathrm{y},\tau)\psi(\mathrm{y},\tau)d\sigma_{\mathrm{y}}d\tau
		\\ \nonumber &= \frac{1}{\alpha} \displaystyle\int_{-\infty}^{\Tilde{t}}\int_{\partial B} \dfrac{\alpha}{4\pi\alpha\delta^2(\Tilde{t}-\Tilde{\tau})} e^{-\dfrac{\alpha \delta^2|\xi-\eta|^2}{4\alpha\delta^2(\Tilde{t}-\Tilde{\tau})}}\psi(\delta \eta + \mathrm{z}, \delta^2 \Tilde{\tau})\delta^3 d\sigma_{\eta}d\Tilde{\tau}
		\\ \nonumber &= \delta \displaystyle\int_{-\infty}^{\Tilde{t}}\int_{\partial B} \dfrac{1}{4\pi(\Tilde{t}-\Tilde{\tau})} e^{-\dfrac{|\xi-\eta|^2}{4(\Tilde{t}-\Tilde{\tau})}}\hat{\psi}(\eta,\Tilde{\tau})d\sigma_{\eta}d\hat{\tau}  \\ &=  \delta \overline{\mathcal{S}}_{\partial B \times \mathbb{R}}\Big[\hat{\psi}\Big](\xi,\Tilde{t})
	\end{align}
	which gives the first identity. The second identity can be derived from the following
	$$ \mathcal{S}_{\partial \Omega \times \mathbb{R}}\Big[\Big(\overline{\mathcal{S}}^{-1}_{\partial B \times \mathbb{R}}\Big[\hat{\varphi}\Big]\Big)^\vee\Big] = \delta \Bigg(\overline{\mathcal{S}}_{\partial B \times \mathbb{R}} \Big(\overline{\mathcal{S}}^{-1}_{\partial B \times \mathbb{R}}\Big[\hat{\varphi}\Big]\Big)\Bigg)^\vee = \delta(\hat{\varphi})^\vee = \delta \varphi.$$
	Using the definition of the operator norm, to derive the estimate (\ref{operatoresti}), we obtain
	\begin{align}
		\nonumber
		\Big\Vert \mathcal{S}^{-1}_{\partial \Omega \times \mathbb{R}} \Big\Vert_{\mathcal{L}{\Big( \mathrm{H}^{\frac{1}{2},\frac{1}{4}}\Big(\partial \Omega \times \mathbb{R}\Big), \mathrm{H}^{-\frac{1}{2},-\frac{1}{4}}\Big(\partial \Omega \times \mathbb{R}\Big)}\Big) } &:= \sup_{0\neq \varphi \in \mathrm{H}^{\frac{1}{2},\frac{1}{4}}\Big(\partial \Omega \times \mathbb{R}\Big)} \dfrac{\Big\Vert \mathcal{S}^{-1}_{\partial \Omega \times \mathbb{R}}[\varphi] \Big\Vert_{\mathrm{H}^{-\frac{1}{2},-\frac{1}{4}}\Big(\partial \Omega \times \mathbb{R}\Big)}}{\Big\Vert \varphi \Big\Vert_{\mathrm{H}^{\frac{1}{2},\frac{1}{4}}\Big(\partial \Omega \times \mathbb{R}\Big)}} \\ \nonumber &\le \sup_{0\neq \varphi \in \mathrm{H}^{\frac{1}{2},\frac{1}{4}}\Big(\partial \Omega \times \mathbb{R}\Big)} \dfrac{\alpha\delta^{\frac{3}{2}}\Big\Vert \Big(\mathcal{S}^{-1}_{\partial \Omega \times \mathbb{R}}[\varphi] \Big)^\wedge\Big\Vert_{\mathrm{H}^{-\frac{1}{2},-\frac{1}{4}}\Big(\partial \Omega \times \mathbb{R}\Big)}}{\alpha \delta^{\frac{3}{2}}\Big\Vert \hat{\varphi} \Big\Vert_{\mathrm{H}^{\frac{1}{2},\frac{1}{4}}\Big(\partial B\times \mathbb{R}\Big)}} \\ \nonumber &= \sup_{0\neq \hat{\varphi} \in \mathrm{H}^{\frac{1}{2},\frac{1}{4}}\Big(\partial B\times \mathbb{R}\Big)} \dfrac{\Big\Vert \Big(\delta^{-1} \Big(\overline{\mathcal{S}}^{-1}_{\partial B\times \mathbb{R}}\Big[\hat{\varphi}\Big]\Big)^\vee \Big)^\wedge\Big\Vert_{\mathrm{H}^{-\frac{1}{2},-\frac{1}{4}}\Big(\partial \Omega \times \mathbb{R}\Big)}}{\Big\Vert \hat{\varphi} \Big\Vert_{\mathrm{H}^{\frac{1}{2},\frac{1}{4}}\Big(\partial B\times \mathbb{R}\Big)}} \\ \nonumber 
		&= \sup_{0\neq \hat{\varphi} \in \mathrm{H}^{\frac{1}{2},\frac{1}{4}}\Big(\partial B\times \mathbb{R}\Big)} \dfrac{\delta^{-1}\Big\Vert\overline{\mathcal{S}}^{-1}_{\partial B\times \mathbb{R}}\Big[\hat{\varphi}\Big]\Big\Vert_{\mathrm{H}^{-\frac{1}{2},-\frac{1}{4}}\Big(\partial B\times \mathbb{R}\Big)}}{\Big\Vert \hat{\varphi} \Big\Vert_{\mathrm{H}^{\frac{1}{2},\frac{1}{4}}\Big(\partial B\times \mathbb{R}\Big)}} \\
		&= \delta^{-1} \Big\Vert \overline{\mathcal{S}}^{-1}_{\partial B \times \mathbb{R}} \Big\Vert_{\mathcal{L}{\Big( \mathrm{H}^{\frac{1}{2},\frac{1}{4}}\Big(\partial B\times \mathbb{R}\Big), \mathrm{H}^{-\frac{1}{2},-\frac{1}{4}}\Big(\partial B\times \mathbb{R}\Big)}\Big)}.
	\end{align}
	The proof is thus complete.
	\begin{lemma}\label{4.7}
		Suppose $0<\delta\leq 1$ and $\Omega=\delta B + \mathrm{z}$ and $ t = \alpha \delta^2\Tilde{t} $. Then for $\psi \in \mathrm{H}^{-\frac{1}{2},-\frac{1}{4}}\Big(\partial \Omega \times \mathbb{R}\Big) $ we have
		\begin{align}\label{l1}
			\alpha \delta^5\Big\Vert \overline{\mathcal{S}}\Big[\hat{\psi}\Big]\Big\Vert_{\mathrm{H}^{\frac{1}{2},\frac{1}{4}}\Big(\partial B\times \mathbb{R}\Big)}^2 \leq \Big\Vert \mathcal{S}\Big[\psi\Big]\Big\Vert_{\mathrm{H}^{\frac{1}{2},\frac{1}{4}}\Big(\partial \Omega\times \mathbb{R}\Big)}^2 \leq  \alpha^{\frac{1}{^2}}\delta^4 \Big\Vert \overline{\mathcal{S}}\Big[\hat{\psi}\Big]\Big\Vert_{\mathrm{H}^{\frac{1}{2},\frac{1}{4}}\Big(\partial B\times \mathbb{R}\Big)}^2.
		\end{align}
	\end{lemma}
	\textbf{Proof.} Based on the proof of Lemma \ref{3.7}, we conclude that
	$$ 
	\mathcal{S}_{\partial \Omega\times\mathbb{R}}\Big[\psi\Big](\mathrm{x},t) =  \delta \Big(\overline{\mathcal{S}}_{\partial B\times \mathbb{R}}\Big[\hat{\psi}\Big]\Big)^\vee,
	$$
	where, $\overline{\mathcal{S}}$ is the single layer operator independent of the diffusion constant $\alpha$ is set to be $1$. In order to scale $\Big\Vert \mathcal{S}\Big[\psi\Big]\Big\Vert_{\mathrm{H}^{\frac{1}{2},\frac{1}{4}}\Big(\partial \Omega \times \mathbb{R}\Big)}^2$, from the definition of norm (\ref{def1}), we need first to consider $\Big\Vert \mathcal{S}\Big[\psi\Big]\Big\Vert_{\mathrm{L}^2\Big(\partial \Omega \times \mathbb{R}\Big)}^2.$\\
	Indeed,
	\begin{align}\label{estil1}
		\nonumber
		\Big\Vert\mathcal{S}\Big[\Psi\Big] \Big\Vert_{\mathrm{L}^2\Big(\partial \Omega \times \mathbb{R}\Big)}^2 &= \int_{-\infty}^T \int_{\partial\Omega} \Bigg| \frac{1}{\alpha}\displaystyle\int_{-\infty}^{T}\int_{\partial\Omega} \Phi(\mathrm{x},t;\mathrm{y},\tau)\psi(y,\tau)d\sigma_{\mathrm{y}}d\tau \Bigg|^2 d\sigma_{\mathrm{x}} dt \\ \nonumber&= \alpha \delta^5 \int_{-\infty}^{\mathrm{T}_\delta} \int_{\partial B} \Bigg| \displaystyle\int_{-\infty}^{\mathrm{T}_\delta}\int_{\partial B} \Phi(\xi,\Tilde{t};\eta,\Tilde{\tau})\hat{\psi}(\eta,\Tilde{\tau})d\sigma_\eta d\Tilde{\tau} \Bigg|^2 d\sigma_\xi d\Tilde{t} \\ &= \alpha \delta^5 \Big\Vert\overline{\mathcal{S}}\Big[\hat{\psi}\Big] \Big\Vert_{\mathrm{L}^2\Big(\partial B \times \mathbb{R}\Big)}^2.
	\end{align}
	Now, in order to estimate (\ref{def2}) and (\ref{def3}), we need to scale  $\Big\Vert\mathcal{S}\Big[\psi\Big] \Big\Vert_{\mathrm{L}^2(\mathbb{R})}^2$ and $ \Big\Vert \mathcal{S}\Big[\psi\Big]\Big \Vert_{\mathrm{L}^2(\partial \Omega)}^2$ respectively.\\
	So,
	\begin{align}\label{estil2}
		\nonumber
		\Big\Vert\mathcal{S}\Big[\psi\Big] \Big\Vert_{\mathrm{L}^2(\mathbb{R})}^2 &= \int_{-\infty}^T \Bigg| \frac{1}{\alpha}\displaystyle\int_{-\infty}^{T}\int_{\partial\Omega} \Phi(\mathrm{x},t;\mathrm{y},\tau)\psi(y,\tau)d\sigma_\mathrm{y} d\tau \Bigg|^2 dt \\ \nonumber &= \alpha \delta^4 \int_{-\infty}^{T_{\delta}} \Bigg| \displaystyle\int_{-\infty}^{T_{\delta}}\int_{\partial B} \Phi(\xi,\Tilde{t};\eta,\Tilde{\tau})\hat{\psi}(\eta,\Tilde{\tau})d\sigma_\eta d\Tilde{\tau} \Bigg|^2 d\Tilde{t} \\ &= \alpha \delta^4 \Big\Vert\overline{\mathcal{S}}\Big[\hat{\psi}\Big] \Big\Vert_{\mathrm{L}^2(\mathbb{R})}^2.
	\end{align}
	Also,
	\begin{align}\label{estil3}
		\nonumber
		\Big\Vert\mathcal{S}\Big[\psi\Big] \Big\Vert_{\mathrm{L}^2(\partial \Omega)}^2 &= \int_{\partial\Omega} \Bigg| \frac{1}{\alpha}\displaystyle\int_{-\infty}^{T}\int_{\partial\Omega} \Phi(\mathrm{x},t;\mathrm{y},\tau)\psi(y,\tau)d\sigma_\mathrm{y} d\tau \Bigg|^2 d\sigma_\mathrm{x}\\ \nonumber&= \delta^3 \int_{\partial B} \Bigg| \displaystyle\int_{-\infty}^{T_{\delta}}\int_{\partial B} \Phi(\xi,\Tilde{t};\eta,\Tilde{\tau})\hat{\psi}(\eta,\Tilde{\tau})d\sigma_\eta d\Tilde{\tau} \Bigg|^2 d\sigma_\xi \\ &= \delta^3 \Big\Vert\overline{\mathcal{S}}\Big[\hat{\psi}\Big] \Big\Vert_{\mathrm{L}^2(\partial B)}^2.
	\end{align}
	So, using (\ref{estil2}) and (\ref{estil3}), we calculate the following norms
	\begin{align}\label{estil4}
		\nonumber
		\Big|\mathcal{S}\Big[\psi\Big]\Big|_{\mathrm{L}^2\Big(\mathbb{R};H^{\frac{1}{2}}(\partial\Omega)\Big)}^2 &= \int_{\partial\Omega}\int_{\partial\Omega}\dfrac{\Big\Vert \mathcal{S}\Big[\psi\Big](x,\cdot)-\mathcal{S}\Big[\psi\Big](y,\cdot)\Big\Vert_{\mathrm{L}^2(\mathbb{R})}^2}{|x-y|^{2}}d\sigma_\mathrm{x}d\sigma_\mathrm{y}\\ \nonumber &= \delta^2 \int_{\partial B}\int_{\partial B}\dfrac{\alpha \delta^4 \Big\Vert \overline{\mathcal{S}}\Big[\hat{\psi}\Big](\xi,\cdot) - \overline{\mathcal{S}}\Big[\hat{\psi}\Big](\eta,\cdot)\Big\Vert_{\mathrm{L}^2(\mathbb{R})}^2}{\delta^2|\xi - \eta |^{2}}d\sigma_\xi d\sigma_\eta \\ &= \alpha \delta^4 \Big|\overline{\mathcal{S}}\Big[\hat{\psi}\Big]\Big|_{\mathrm{L}^2\Big(\mathbb{R};H^{\frac{1}{2}}(\partial B)\Big)}^2 
	\end{align}
	Again,
	\begin{align}\label{estil5}
		\nonumber
		\Big|\mathcal{S}\Big[\psi\Big]\Big|_{H^{\frac{1}{4}}\Big(\mathbb{R};\mathrm{L}^2(\partial\Omega)\Big)}^2 &= \int_{-\infty}^T\int_{-\infty}^T\dfrac{\Big\Vert \mathcal{S}\Big[\psi\Big](\cdot,t)-\mathcal{S}\Big[\psi\Big](\cdot,\tau)\Big\Vert_{\mathrm{L}^2(\partial \Omega)}^2}{|t-\tau|^{\frac{3}{2}}}dt d\tau\\ \nonumber &= \alpha^2\delta^4 \int_{-\infty}^{T_{\delta}}\int_{-\infty}^{T_{\delta}}\dfrac{ \delta^3\Big\Vert \overline{\mathcal{S}}\Big[\hat{\psi}\Big](\xi,\cdot) - \overline{\mathcal{S}}\Big[\hat{\psi}\Big](\eta,\cdot)\Big\Vert_{\mathrm{L}^2(\partial B)}^2}{\alpha^{\frac{3}{2}}\delta^3|\Tilde{t}-\Tilde{\tau}|^{\frac{3}{2}}}d\Tilde{t}d\Tilde{\tau} \\ &= \alpha^{\frac{1}{2}} \delta^4 \Big|\overline{\mathcal{S}}\Big[\hat{\psi}\Big]\Big|_{H^{\frac{1}{4}}\Big(\mathbb{R};\mathrm{L}^2(\partial B)\Big)}^2 
	\end{align}
	In this case, we obtain our desired result by estimating the norms (\ref{estil1}), (\ref{estil4}) and (\ref{estil5}), the norm (\ref{def1}) and, noticing the fact that ${\color{Black}{  \alpha \delta^5} \leq \alpha \delta^4\leq \alpha^{\frac{1}{2}} \delta^4}$,
	\begin{align}\label{l1}
		\alpha \delta^5 \Big\Vert \overline{\mathcal{S}}\Big[\hat{\psi}\Big]\Big\Vert_{\mathrm{H}^{\frac{1}{2},\frac{1}{4}}\Big(\partial B \times \mathbb{R}\Big)}^2 \leq \Big\Vert \mathcal{S}\Big[\psi\Big]\Big\Vert_{\mathrm{H}^{\frac{1}{2},\frac{1}{4}}\Big(\partial \Omega \times \mathbb{R}\Big)}^2 \leq \alpha^{\frac{1}{2}} \delta^4 \Big\Vert \overline{\mathcal{S}}\Big[\hat{\psi}\Big]\Big\Vert_{\mathrm{H}^{\frac{1}{2},\frac{1}{4}}\Big(\partial B \times \mathbb{R}\Big)}^2.
	\end{align}
	\begin{lemma}\label{4.8}
		For $\psi \in \mathrm{L}^2\Big(\Omega\times\mathbb{R}\Big) $ we have the following identities
		\begin{align}
			\mathcal{V}_{\Omega\times\mathbb{R}}\Big[\psi\Big](\mathrm{x},t) =  \delta^2 \Big(\overline{\mathcal{V}}_{B\times\mathbb{R}}\Big[\hat{\psi}\Big]\Big)^\vee
		\end{align}
		and
		\begin{align}\label{l1}
			\alpha \delta^7\Big\Vert \overline{\mathcal{V}}\Big[\hat{\psi}\Big]\Big\Vert_{\mathrm{H}^{\frac{1}{2},\frac{1}{4}}\Big(\partial B\times \mathbb{R}\Big)}^2 \leq \Big\Vert \mathcal{V}\Big[\psi\Big]\Big\Vert_{\mathrm{H}^{\frac{1}{2},\frac{1}{4}}\Big(\partial \Omega\times \mathbb{R}\Big)}^2 \leq  \alpha^{\frac{1}{^2}}\delta^6 \Big\Vert \overline{\mathcal{V}}\Big[\hat{\psi}\Big]\Big\Vert_{\mathrm{H}^{\frac{1}{2},\frac{1}{4}}\Big(\partial B\times \mathbb{R}\Big)}^2.
		\end{align}
	\end{lemma}
	\textbf{Proof.} Setting again $\mathrm{x} := \delta \xi + \mathrm{z}$, $\mathrm{y}:= \delta \eta + \mathrm{z}$, $t := \alpha \delta^2 \Tilde{t}$ and $\tau := \alpha \delta^2 \Tilde{\tau}$, we have
	\begin{align}
		\nonumber
		\mathcal{V}_{\Omega\times\mathbb{R}}\Big[\psi\Big](\mathrm{x},t) &=  \frac{1}{\alpha}\displaystyle\int_{-\infty}^{T}\int_{\Omega} \Phi(\mathrm{x},t;\mathrm{y},\tau)\psi(y,\tau)d\mathrm{y}d\tau
		\\ \nonumber &= \frac{1}{\alpha} \displaystyle\int_{-\infty}^{T_{\delta}}\int_{B} \dfrac{\alpha}{4\pi\alpha \delta^2(\Tilde{t}-\Tilde{\tau})} e^{-\dfrac{\alpha \delta^2|\xi-\eta|^2}{4\alpha\delta^2(\Tilde{t}-\Tilde{\tau})}}\hat{\psi}(\delta \eta + \mathrm{z},\alpha\delta^2\Tilde{\tau})\alpha \delta^4 d\eta d\Tilde{\tau}
		\\ \nonumber &= \delta^2  \displaystyle\int_{-\infty}^{T_{\delta}}\int_{B} \dfrac{1}{4\pi(\Tilde{t}-\Tilde{\tau})} e^{-\dfrac{|\xi-\eta|^2}{4(\Tilde{t}-\Tilde{\tau})}}\hat{\psi}(\eta,\Tilde{\tau})d\eta d\Tilde{\tau}  \\ &=  \delta^2 \overline{\mathcal{V}}_{B\times\mathbb{R}}\Big[\hat{\psi}\Big](\xi,\Tilde{t}).
	\end{align}
	where, $\overline{\mathcal{V}}$ is the Newtonian heat operator corresponding to $\alpha=1$ and to the source $\psi.$\\
	In order to scale $\Big\Vert \mathcal{V}\Big[\psi\Big]\Big\Vert_{\mathrm{H}^{\frac{1}{2},\frac{1}{4}}\Big(\partial \Omega \times \mathbb{R}\Big)}^2$, from the definition of norm (\ref{def1}), we need to first calculate $\Big\Vert \mathcal{V}\Big[\psi\Big]\Big\Vert_{\mathrm{L}^2\Big(\partial \Omega \times \mathbb{R}\Big)}^2.$\\
	Indeed,
	\begin{align}\label{estill1}
		\nonumber
		\Big\Vert\mathcal{V}\Big[\psi\Big] \Big\Vert_{\mathrm{L}^2\Big(\partial \Omega \times \mathbb{R}\Big)}^2 &= \int_{-\infty}^T \int_{\partial\Omega} \Bigg| \frac{1}{\alpha}\displaystyle\int_{-\infty}^{T}\int_{\Omega} \Phi(\mathrm{x},t;\mathrm{y},\tau)\Psi(y,\tau)d\mathrm{y}d\tau \Bigg|^2 d\sigma_\mathrm{x} dt \\ \nonumber&= \alpha \delta^7 \int_{-\infty}^{T_{\delta}} \int_{\partial B} \Bigg| \displaystyle\int_{-\infty}^{T_{\delta}}\int_{B} \Phi(\xi,\Tilde{t};\eta,\Tilde{\tau})\hat{\psi}(\eta,\Tilde{\tau})d\eta d\Tilde{\tau} \Bigg|^2 d\sigma_\xi d\Tilde{t} \\ &= \alpha \delta^7 \Big\Vert\overline{\mathcal{V}}\Big[\hat{\psi}\Big] \Big\Vert_{\mathrm{L}^2\Big(\partial B \times \mathbb{R}\Big)}^2.
	\end{align}
	Further, in order to calculate (\ref{def2}) and (\ref{def3}), we need to scale the following $\Big\Vert\mathcal{V}\Big[\psi\Big] \Big\Vert_{\mathrm{L}^2(\mathbb{R})}$ and $\Big\Vert \mathcal{V}\Big[\psi\Big]\Big \Vert_{\mathrm{L}^2(\partial \Omega)}^2$ respectively. Indeed,
	\begin{align}\label{estill2}
		\nonumber
		\Big\Vert\mathcal{V}\Big[\psi\Big] \Big\Vert_{\mathrm{L}^2(\mathbb{R})}^2 &= \int_{-\infty}^T \Bigg| \frac{1}{\alpha}\displaystyle\int_{-\infty}^{T}\int_{\Omega} \Phi(\mathrm{x},t;\mathrm{y},\tau)\Psi(y,\tau)d\mathrm{y} d\tau \Bigg|^2 dt \nonumber \\ \nonumber &= \alpha \delta^6 \int_{-\infty}^{T_{\delta}} \Bigg| \displaystyle\int_{-\infty}^{T_{\delta}}\int_{B} \Phi(\xi,\Tilde{t};\eta,\Tilde{\tau})\hat{\psi}(\eta,\Tilde{\tau})d\eta d\Tilde{\tau} \Bigg|^2 d\hat{t} \\  &= \alpha \delta^6 \Big\Vert\overline{\mathcal{V}}\Big[\hat{\psi}\Big] \Big\Vert_{\mathrm{L}^2(\mathbb{R})}^2
	\end{align}
	and
	\begin{align}\label{estill3}
		\nonumber
		\Big\Vert\mathcal{V}\Big[\psi\Big] \Big\Vert_{\mathrm{L}^2(\partial \Omega)}^2 &= \int_{\partial\Omega} \Bigg| \frac{1}{\alpha}\displaystyle\int_{-\infty}^{T}\int_{\Omega} \Phi(\mathrm{x},t;\mathrm{y},\tau)\Psi(y,\tau)d\mathrm{y} d\tau \Bigg|^2 d\sigma_\mathrm{x} \\ \nonumber &= \delta^5 \int_{\partial B} \Bigg| \displaystyle\int_{-\infty}^{T_{\delta}}\int_{B} \Phi(\xi,\Tilde{t};\eta,\Tilde{\tau})\hat{\psi}(\eta,\Tilde{\tau})d\eta d\Tilde{\tau} \Bigg|^2 d\sigma_\xi \\ &= \delta^5 \Big\Vert\overline{\mathcal{V}}\Big[\hat{\psi}\Big] \Big\Vert_{\mathrm{L}^2(\partial B)}^2.
	\end{align}
	Based on these two expressions, we can calculate the following norms
	\begin{align}\label{estill4}
		\nonumber
		\Big|\mathcal{V}\Big[\psi\Big]\Big|_{\mathrm{L}^2\Big(\mathbb{R};H^{\frac{1}{2}}(\partial\Omega)\Big)}^2 &= \int_{\partial\Omega}\int_{\partial\Omega}\dfrac{\Big\Vert \mathcal{V}\Big[\psi\Big](x,\cdot)-\mathcal{V}\Big[\psi\Big](y,\cdot)\Big\Vert_{\mathrm{L}^2(\mathbb{R})}^2}{|x-y|^{2}}d\sigma_\mathrm{x}d\sigma_\mathrm{y}\\ \nonumber &= \delta^2 \int_{\partial B}\int_{\partial B}\dfrac{\alpha \delta^6 \Big\Vert \overline{\mathcal{V}}\Big[\hat{\psi}\Big](\xi,\cdot)-\overline{\mathcal{V}}\Big[\hat{\psi}\Big](\eta,\cdot)\Big\Vert_{\mathrm{L}^2(\mathbb{R})}^2}{\delta^2|\xi - \eta|^{2}}d\sigma_\xi d\sigma_\eta \\ &= \alpha \delta^6 \Big|\overline{\mathcal{V}}\Big[\hat{\psi}\Big]\Big|_{\mathrm{L}^2\Big(\mathbb{R};H^{\frac{1}{2}}(\partial B)\Big)}^2 
	\end{align}
	and
	\begin{align}\label{estill5}
		\nonumber
		\Big|\mathcal{V}\Big[\psi\Big]\Big|_{H^{\frac{1}{4}}\Big(\mathbb{R};\mathrm{L}^2(\partial\Omega)\Big)}^2 &= \int_{-\infty}^T\int_{-\infty}^T\dfrac{\Big\Vert \mathcal{V}\Big[\psi\Big](\cdot,t)-\mathcal{V}\Big[\psi\Big](\cdot,\tau)\Big\Vert_{\mathrm{L}^2(\partial \Omega)}^2}{|t-\tau|^{\frac{3}{2}}}d\mathrm{t}d\tau\\ \nonumber &= \alpha^2\delta^4 \int_{-\infty}^{T_{\delta}}\int_{-\infty}^{T_{\delta}}\dfrac{ \delta^5\Big\Vert \overline{\mathcal{V}}\Big[\hat{\psi}\Big](\xi,\cdot)-\overline{\mathcal{V}}\Big[\hat{\psi}\Big](\eta,\cdot)\Big\Vert_{\mathrm{L}^2(\partial B)}^2}{\alpha^{\frac{3}{2}}\delta^3|\Tilde{t}-\Tilde{\tau}|^{\frac{3}{2}}}d\Tilde{t}d\Tilde{\tau} \\ &= \alpha^{\frac{1}{2}} \delta^6 \Big|\overline{\mathcal{V}}\Big[\hat{\Psi}\Big]\Big|_{H^{\frac{1}{4}}\Big(\mathbb{R};\mathrm{L}^2(\partial B)\Big)}^2 
	\end{align}
	Hence, considering the norm-calculation (\ref{estill1}), (\ref{estill4}) and (\ref{estill5}), from the norm definition (\ref{def1}), noticing the fact that ${\color{Black}{\alpha \delta^7} \leq \alpha \delta^6\leq \alpha^{\frac{1}{2}} \delta^6}$, we get our desired inequality
	\begin{align}\label{l1}
		\alpha \delta^7 \Big\Vert \noindent
		\overline{\mathcal{V}}\Big[\hat{\psi}\Big]\Big\Vert_{\mathrm{H}^{\frac{1}{2},\frac{1}{4}}\Big(\partial B \times \mathbb{R}\Big)}^2 \leq \Big\Vert \mathcal{V}\Big[\psi\Big]\Big\Vert_{\mathrm{H}^{\frac{1}{2},\frac{1}{4}}\Big(\partial \Omega \times \mathbb{R}\Big)}^2 \leq \alpha^{\frac{1}{2}} \delta^6 \Big\Vert \overline{\mathcal{V}}\Big[\hat{\psi}\Big]\Big\Vert_{\mathrm{H}^{\frac{1}{2},\frac{1}{4}}\Big(\partial B\times\mathbb{R}\Big)}^2
	\end{align}
	Hence it completes the proof.
	\begin{lemma}\label{3.11}
		Suppose $0<\delta\leq 1$, $\Omega:=\delta B + \mathrm{z}$ and $t := \alpha \delta^2\Tilde{t} $. Then for $ \psi \in L_{2}\Big(\Omega \times \mathbb{R}\Big) $ we have the following
		\begin{align}
			\gamma^{\textbf{int}}_{1,\mathrm{x}}\mathcal{V}_{\Omega \times\mathbb{R}}\Big[\psi\Big](\mathrm{x},t) =  \delta \Big(\gamma^{\textbf{int}}_{1,\xi}\overline{\mathcal{V}}_{B\times\mathbb{R}}\Big[\hat{\psi}\Big]\Big)^\vee
		\end{align}
		and
		\begin{align}\label{l1}
			\alpha^{\frac{3}{4}} \delta^3\Big\Vert \gamma^{\textbf{int}}_{1,\xi}\overline{\mathcal{V}}\Big[\hat{\psi}\Big]\Big\Vert_{\mathrm{H}^{-\frac{1}{2},-\frac{1}{4}}\Big(\partial B\times \mathbb{R}\Big)} \leq \Big\Vert \gamma^{\textbf{int}}_{1,\mathrm{x}}\mathcal{V}\Big[\psi\Big]\Big\Vert_{\mathrm{H}^{-\frac{1}{2},-\frac{1}{4}}\Big(\partial \Omega\times \mathbb{R}\Big)}\leq  \alpha^{\frac{1}{^2}}\delta^{\frac{5}{2}} \Big\Vert \gamma^{\textbf{int}}_{1,\xi}\overline{\mathcal{V}}\Big[\hat{\psi}\Big]\Big\Vert_{\mathrm{H}^{-\frac{1}{2},-\frac{1}{4}}\Big(\partial B\times \mathbb{R}\Big)}.
		\end{align}
	\end{lemma}
	\textbf{Proof.} Before proceeding to prove the first identity, we recall that on the boundary $(\partial\Omega)_{\mathrm{T}}$ i.e. $ \gamma^{\textbf{int}}_{1,\mathrm{x}} \mathrm{u} = \partial_{\nu}u,$ with $\nu$ being the outward unit normal vector to $\partial\Omega.$ 
	\begin{align}
		\nonumber
		\partial_{\nu_{\mathrm{x}}}\mathcal{V}_{\Omega\times\mathbb{R}}\Big[\psi\Big](\mathrm{x},t) &=  \frac{1}{\alpha}\partial_{\nu_{\mathrm{x}}}\displaystyle\int_{-\infty}^{T}\int_{\Omega} \Phi(\mathrm{x},t;\mathrm{y},\tau)\psi(y,\tau)d\mathrm{y}d\tau
		\\ \nonumber &=  \frac{1}{\alpha}\displaystyle\int_{-\infty}^{T}\int_{\Omega}-\dfrac{\alpha(x-y)\cdot \nu_{y}}{2(t-\tau)}  \Phi(\mathrm{x},t;\mathrm{y},\tau)\psi(y,\tau)d\mathrm{y}d\tau
		\\ \nonumber &= \frac{1}{\alpha} \displaystyle\int_{-\infty}^{T_{\delta}}\int_{B} -\dfrac{\alpha\delta(\xi-\eta)\cdot \nu_{\eta}}{2\alpha\delta^2(\Tilde{t}-\Tilde{\tau})}\cdot \dfrac{\alpha}{4\pi\alpha \delta^2(\Tilde{t}-\Tilde{\tau})} e^{-\dfrac{\alpha \delta^2|\xi-\eta|^2}{4\alpha\delta^2(\Tilde{t}-\Tilde{\tau})}}\hat{\psi}(\eta,\Tilde{\tau})\alpha \delta^4 d\eta d\Tilde{\tau}
		\\ \nonumber &= \delta  \partial_{\nu_{\xi}}\displaystyle\int_{-\infty}^{T_{\delta}}\int_{B} \Phi(\xi,\Tilde{t};\eta,\Tilde{\tau})\hat{\psi}(\eta,\Tilde{\tau})d\eta d\Tilde{\tau}  \\ &=  \delta \partial_{\nu_{\xi}}\overline{\mathcal{V}}_{B\times\mathbb{R}}\Big[\hat{\Psi}\Big](\hat{x},\hat{\tau}),
	\end{align}
	where $\partial_{\nu_{\xi}}\overline{\mathcal{V}}$ is related to the Newtonian heat operator with $\alpha=1$ and corresponding to the source $\psi$. Hence the first identity is proved.\\
	Now, to get the desired inequality, we need to scale the product $\Big\langle \varphi,\partial_{\nu_\mathrm{x}}\mathcal{V}\Big[\psi\Big] \Big\rangle_{\partial \Omega \times \mathbb{R}}$.
	\begin{align}
		\nonumber
		\Big\langle \varphi,\partial_{\nu_\mathrm{x}}\mathcal{V}\Big[\psi\Big] \Big\rangle_{\partial \Omega\times\mathbb{R}} &= \int_{-\infty}^{T}\int_{\partial \Omega} \varphi(\mathrm{x},t)\frac{1}{\alpha}\partial_{\nu_{\mathrm{x}}}\displaystyle\int_{-\infty}^{T}\int_{\Omega} \Phi(\mathrm{x},t;\mathrm{y},\tau)\psi(y,\tau)d\mathrm{y}d\tau d\mathrm{y}d\tau d\sigma_\mathrm{x}d\mathrm{t}\\ \nonumber &= \alpha \delta^4 \int_{-\infty}^{T_{\delta}}\int_{\partial \Omega} \varphi(\xi,\Tilde{t})\partial_{\nu_\xi}\displaystyle\int_{-\infty}^{T_{\delta}}\int_{\Omega} \Phi(\xi,\Tilde{t};\eta,\Tilde{\tau})\hat{\psi}(\xi,\Tilde{\tau})d\eta d\Tilde{\tau}d\sigma_\xi d\Tilde{t} \\ &= \alpha \delta^4 \Big\langle \hat{\varphi},\partial_{\nu_\xi}\overline{\mathcal{V}}\Big[\hat{\psi}\Big] \Big\rangle_{\partial B\times\mathbb{R}}.
	\end{align}
	Hence, from the definition of norm definition (\ref{def4}), we get our desired second inequality,
	\begin{align}\label{l2}
		\alpha^{\frac{3}{4}}\delta^3\Big\Vert \gamma^{\textbf{int}}_{1,\xi}\overline{\mathcal{V}}\Big[\hat{\psi}\Big]\Big\Vert_{\mathrm{H}^{-\frac{1}{2},-\frac{1}{4}}\Big(\partial B\times\mathbb{R}\Big)} \leq \Big\Vert \gamma^{\textbf{int}}_{1,\mathrm{x}}\mathcal{V}\Big[\psi\Big]\Big\Vert_{\mathrm{H}^{-\frac{1}{2},-\frac{1}{4}}\Big(\partial \Omega\times\mathbb{R}\Big)} \leq  \alpha^{\frac{1}{2}}\delta^{\frac{5}{2}}\Big\Vert \gamma^{\textbf{int}}_{1,\xi}\overline{\mathcal{V}}\Big[\hat{\psi}\Big]\Big\Vert_{\mathrm{H}^{-\frac{1}{2},-\frac{1}{4}}\Big(\partial B\times\mathbb{R}\Big)}.
	\end{align}
	Thus it completes the proof.
	\begin{lemma}
		For $\psi \in \mathrm{L}^2\Big(\Omega\Big) $ we have the following
		\begin{align}\label{l1}
			\alpha^{\frac{1}{2}} \delta^2\Big\Vert \overline{\mathbb{I}}\Big[\hat{\psi}\Big]\Big\Vert_{\mathrm{H}^{1,\frac{1}{2}}\Big(B\times\mathbb{R}_+\Big)} \leq \Big\Vert \mathbb{I}\Big[\psi\Big]\Big\Vert_{\mathrm{H}^{1,\frac{1}{2}}\Big(\Omega\times\mathbb{R}_+\Big)} \leq  \alpha^{\frac{1}{2}}\delta \Big\Vert \overline{\mathbb{I}}\Big[\hat{\psi}\Big]\Big\Vert_{\mathrm{H}^{1,\frac{1}{2}}\Big(B\times\mathbb{R}_+\Big)}.
		\end{align}
	\end{lemma} 
	\textbf{Proof.} We define the initial heat potential for $\psi \in \mathrm{L}^2\Big(\Omega \Big)$ as follows
	\begin{align}
		\nonumber
		\mathbb{I}\Big[\psi\Big](\mathrm{x},t) = \int_{\Omega}\Phi(x,t;y)\psi(y)d\mathrm{y}.
	\end{align}
	We also note that for $\psi \in \mathrm{L}^2\Big(\Omega\Big)$, we have the following two estimates
	\begin{align}\label{i1}
		\mathbb{I}_{(\Omega)}\Big[\psi\Big](\mathrm{x},t) = \Big(\overline{\mathbb{I}}_{(B)}\Big[\hat{\psi}\Big]\Big)^\vee,
	\end{align}
	and
	\begin{align}\label{i2}
		\nabla_{\mathrm{x}}\mathbb{I}_{(\Omega)}\Big[\psi\Big](\mathrm{x},t) =\frac{1}{\delta} \Big(\nabla_{\xi}\overline{\mathbb{I}}_{(B)}\Big[\hat{\psi}\Big]\Big)^\vee.
	\end{align}
	From the definition of initial heat potential we obtain
	\begin{align}
		\nonumber
		\mathbb{I}_{(\Omega)}\Big[\psi\Big](\mathrm{x},t) &= \int_{\Omega}\Phi(x,t;y)\psi(y)d\mathrm{y}
		\\ \nonumber &= \int_{\Omega}\frac{\alpha}{4\pi t}\textbf{exp}\Big(-\frac{\alpha|\mathrm{x}-\mathrm{y}|^2}{4t}\Big)\psi(y)d\mathrm{y}
		\\ \nonumber &= \delta^2 \int_B \frac{\alpha}{4\pi \alpha\delta^2\Tilde{t}}\textbf{exp}\Big(-\frac{|\xi-\eta|^2}{4\Tilde{t}}\Big)\psi(\delta\eta + \mathrm{z})d\eta
		\\ \nonumber &= \int_B \frac{1}{4\pi\Tilde{t}}\textbf{exp}\Big(-\frac{|\xi-\eta|^2}{4\Tilde{t}}\Big)\hat{\psi}(\eta )d\eta
		\\ \nonumber &= \Big(\overline{\mathbb{I}}_{(B)}\Big[\hat{\psi}\Big]\Big)^\vee.
	\end{align}
	Hence, it completes the first inequality (\ref{i1}). Next, we provide the proof of the second inequality (\ref{i2}). So,
	\begin{align}
		\nonumber
		\nabla_{\mathrm{x}}\mathbb{I}_{(\Omega)}\Big[\psi\Big](\mathrm{x},t) &= \nabla_{\mathrm{x}}\int_{\Omega}\Phi(x,t;y)\psi(y)d\mathrm{y}
		\\ \nonumber &= \int_{\Omega}\frac{\alpha^2|x-y|}{8\pi t^2}\textbf{exp}\Big(-\frac{\alpha|\mathrm{x}-\mathrm{y}|^2}{4t}\Big)\psi(y)d\mathrm{y}
		\\ \nonumber &= \delta^2 \int_B \frac{\alpha^2\delta|\xi-\eta|}{8\pi \alpha^2\delta^4\Tilde{t}}\textbf{exp}\Big(-\frac{|\xi-\eta|^2}{4\Tilde{t}}\Big)\psi(\delta\eta + \mathrm{z})d\eta
		\\ \nonumber &= \frac{1}{\delta}\int_B \frac{1}{8\pi\Tilde{t}}\textbf{exp}\Big(-\frac{|\xi-\eta|^2}{4\Tilde{t}}\Big)\hat{\psi}(\eta )d\eta
		\\ \nonumber &= \frac{1}{\delta}\Big(\nabla_{\xi}\overline{\mathbb{I}}_{(B)}\Big[\hat{\psi}\Big]\Big)^\vee,
	\end{align}
	which leads to (\ref{i2}).
	Consequently, from the above two estimates we deduce
	\begin{align}
		\nonumber
		\Big\Vert \mathbb{I}_{(\Omega)}\Big[\psi\Big] \Big\Vert_{\mathrm{L}^2(\Omega \times\mathbb{R}_+)}^2 &= \int_{0}^{T}\int_\Omega\Bigg|\int_{\Omega}\Phi(x,t;y)\psi(y)d\mathrm{y}\Bigg|^2 d\mathrm{x}dt
		\\ \nonumber &= \alpha\delta^4\int_{0}^{T_\delta}\int_B\Bigg|\int_{B}\hat{\Phi}(\xi,\Tilde{t};\eta)\hat{\psi}(\eta)d\eta\Bigg|^2 d\xi d\Tilde{t}
		\\  &= \alpha\delta^4 \Big\Vert \mathbb{I}_{(B)}\Big[\hat{\psi}\Big] \Big\Vert_{\mathrm{L}^2(B\times\mathbb{R_+})}^2.
	\end{align}
	In a similar manner, we deduce that 
	\begin{align}
		\Big\Vert \nabla_{\mathrm{x}}\mathbb{I}_{(\Omega)}\Big[\psi\Big] \Big\Vert_{\mathrm{L}^2(\Omega\times\mathbb{R_+})}^2 = \alpha\delta^2 \Big\Vert\nabla_{\xi} \mathbb{I}_{(B)}\Big[\hat{\psi}\Big] \Big\Vert_{\mathrm{L}^2(B\times\mathbb{R_+})}^2,
	\end{align}
	and
	\begin{align}
		\Big\Vert \partial_{t}^{\frac{1}{2}}\mathbb{I}_{(\Omega)}\Big[\psi\Big] \Big\Vert_{\mathrm{L}^2(\Omega\times\mathbb{R_+})}^2 = \alpha\delta^2 \Big\Vert\partial_{\Tilde{t}}^{\frac{1}{2}} \mathbb{I}_{(B)}\Big[\hat{\psi}\Big] \Big\Vert_{\mathrm{L}^2(B\times\mathbb{R_+})}^2.
	\end{align}
	Consequently, from the above three estimates, $H^{1,\frac{1}{2}}$-norm definition and noticing the fact that $\alpha\delta^4\le \alpha\delta^2$ we obtain the desired inequality
	\begin{align}
		\alpha^{\frac{1}{2}} \delta^2\Big\Vert \overline{\mathbb{I}}\Big[\hat{\psi}\Big]\Big\Vert_{\mathrm{H}^{1,\frac{1}{2}}\Big( B\times\mathbb{R_+}\Big)} \leq \Big\Vert \mathbb{I}\Big[\psi\Big]\Big\Vert_{\mathrm{H}^{1,\frac{1}{2}}\Big( \Omega\times\mathbb{R}_+\Big)} \leq  \alpha^{\frac{1}{2}}\delta \Big\Vert \overline{\mathbb{I}}\Big[\hat{\psi}\Big]\Big\Vert_{\mathrm{H}^{1,\frac{1}{2}}\Big( B\times\mathbb{R}_+\Big)}.
	\end{align}
	\section{\Large\textbf{Proof of Theorem \ref{th2}}}\label{sec4}
	In this section, we show the asymptotic analysis of the solution to (\ref{eq:heat1}) as $\delta\to 0$ when a plasmonic nanoparticle occupy a bounded domain $\Omega = \mathrm{z} + \delta \mathrm{B}.$ 
	\subsection{A Priori Estimates of the Heat Potential}
	We start by stating the following Proposition.
	\begin{proposition} \label{pr4}
		We have the following a priori estimate for the solution to the boundary integral equations (\ref{IEfinal})
		\begin{align}\label{gamma0}
			\Big\Vert \gamma^{\textbf{int}}_{0}\mathrm{U}_{\mathrm{i}} \Big\Vert_{\mathrm{H}^{\frac{1}{2},\frac{1}{4}}\Big(\partial \Omega\times \mathbb{R}\Big)} \lesssim \frac{\omega \cdot \boldsymbol{\Im}(\varepsilon_\mathrm{p})}{2\pi\gamma_{\mathrm{p}}}\alpha_{\mathrm{p}}^{\frac{1}{4}} \delta^2 \Big\Vert|\mathrm{E}|^{2}\Big\Vert_{\mathrm{L}^2\Big(\Omega \Big)}.
		\end{align}
		and
		\begin{align}\label{gamma1}
			\Big\Vert \gamma^{\textbf{int}}_{1}\mathrm{U}_{\mathrm{i}} \Big\Vert_{\mathrm{H}^{-\frac{1}{2},-\frac{1}{4}}\Big(\partial \Omega\times \mathbb{R}\Big)} \lesssim \frac{\omega \cdot \boldsymbol{\Im}(\varepsilon_\mathrm{p})}{2\pi\gamma_{\mathrm{p}}}\Big(\alpha^{\frac{1}{4}}_\mathrm{p}\delta + \alpha^{\frac{1}{2}}_\mathrm{p}\delta^{\frac{3}{2}}\Big) \Big\Vert|\mathrm{E}|^{2}\Big\Vert_{\mathrm{L}^2\Big(\Omega\Big)}.
		\end{align}
	\end{proposition}
	\noindent
	\textbf{Proof.} We start with the following equation for $(\mathrm{x},t) \in \partial \Omega\times \mathbb{R}$
	\begin{align}\label{aheat}
		\gamma^{\textbf{int}}_{0}\mathrm{U}_{\mathrm{i}}(\mathrm{x},t) = -\frac{\gamma_{\mathrm{m}}}{\gamma_{\mathrm{p}}} \Big(\frac{1}{2}I_{d} + \mathcal{K}_{\alpha_{\mathrm{p}}}\Big)^{-1}\mathcal{S}_{\alpha_{\mathrm{p}}} \mathbb{A}^{\textbf{ext}}\Big[\gamma^{\textbf{int}}_{0}\mathrm{U}_{\mathrm{i}}\Big](\mathrm{x},t) + \frac{\omega \cdot \boldsymbol{\Im}(\varepsilon_\mathrm{p})}{2\pi\gamma_{\mathrm{p}}}\Big(\frac{1}{2}I_{d} + \mathcal{K}_{\alpha_{\mathrm{p}}}\Big)^{-1}\gamma^{\textbf{int}}_{0,\mathrm{x}}\mathcal{V}\Big[ |\mathrm{E}|^{2}\rchi_{(0,\mathrm{T_0})} \Big](\mathrm{x},t)   
	\end{align}
	From Lemma \ref{3.1} and taking $\mathrm{H}^{\frac{1}{2},\frac{1}{4}}\Big((\partial\Omega \times \mathbb{R}\Big)$-norm on the both sides of the equation (\ref{aheat}) , we obtain
	\begin{align}\label{gamma}
		\nonumber
		\Big\Vert \gamma^{\textbf{int}}_{0}\mathrm{U}_{\mathrm{i}} \Big\Vert_{\mathrm{H}^{\frac{1}{2},\frac{1}{4}}\Big(\partial \Omega \times \mathbb{R}\Big)} &=  \frac{\gamma_{\mathrm{m}}}{\gamma_{\mathrm{p}}}\Big\Vert\Big[\frac{1}{2}I_{d} + \mathcal{K}_{\alpha_\mathrm{p}}\Big]^{-1}\Big\Vert_{\mathcal{L}\Big(\mathrm{H}^{\frac{1}{2},\frac{1}{4}}(\partial \Omega \times \mathbb{R})\Big)} \cdot \Big\Vert \mathcal{S}_{\partial \Omega\times \mathbb{R}} \Big\Vert_{\mathcal{L}{\Big( \mathrm{H}^{-\frac{1}{2},-\frac{1}{4}}\Big((\partial \Omega\times \mathbb{R})\Big), \mathrm{H}^{\frac{1}{2},\frac{1}{4}}\Big((\partial \Omega\times\mathbb{R)}\Big)}\Big) }
		\\ \nonumber & \cdot \Big\Vert \mathbb{A}^{\textbf{ext}}\Big[\gamma^{\textbf{int}}_{0}\mathrm{U}_{\mathrm{i}}\Big] \Big\Vert_{\mathrm{H}^{-\frac{1}{2},-\frac{1}{4}}\Big(\partial \Omega \times \mathbb{R}\Big)} \\  &+  \frac{\omega \cdot \boldsymbol{\Im}(\varepsilon_\mathrm{p})}{2\pi\gamma_{\mathrm{p}}}\Big\Vert\Big[\frac{1}{2}I_{d} + \mathcal{K}_{\alpha_\mathrm{p}}\Big]^{-1}\Big\Vert_{\mathcal{L}\Big(\mathrm{H}^{\frac{1}{2},\frac{1}{4}}(\partial \Omega\times \mathbb{R})\Big)} \Big\Vert \gamma^{\textbf{int}}_{0}\mathcal{V}\Big[ |\mathrm{E}|^{2}\rchi_{(0,\mathrm{T_0})} \Big] \Big\Vert_{\mathrm{H}^{\frac{1}{2},\frac{1}{4}}\Big((\partial \Omega \times \mathbb{R}\Big)}.
	\end{align}
	Furthermore, from Lemma \ref{3.6}, Lemma \ref{4.7} and Lemma \ref{3.1}, i.e. the continuity of the operator $\overline{\mathcal{S}}: \mathrm{H}^{-\frac{1}{2},-\frac{1}{4}}\Big(\partial B\times \mathbb{R}\Big) \to  \mathrm{H}^{\frac{1}{2},\frac{1}{4}}\Big(\partial B\times \mathbb{R}\Big)$, we obtain
	\begin{align}\label{singlelayeresti1}
		\nonumber
		\Big\Vert \mathcal{S}_{\alpha_\mathrm{p}}\Big[\varphi\Big]\Big\Vert_{\mathrm{H}^{\frac{1}{2},\frac{1}{4}}\Big(\partial \Omega\times \mathbb{R}\Big)} &\leq C\alpha_{\mathrm{p}}^{\frac{1}{4}} \delta^{2}   \Big\Vert \overline{\mathcal{S}}_{\alpha_\mathrm{p}}\Big[\hat{\varphi}\Big] \Big\Vert_{\mathrm{H}^{-\frac{1}{2},-\frac{1}{4}}\Big(\partial B\times \mathbb{R}\Big)} 
		\\&\leq \nonumber
		C\alpha_{\mathrm{p}}^{\frac{1}{4}} \delta^{2}   \Big\Vert\hat\varphi \Big\Vert_{\mathrm{H}^{-\frac{1}{2},-\frac{1}{4}}\Big(\partial B\times \mathbb{R}\Big)} 
		\\&\leq
		C\alpha_{\mathrm{p}}^{-\frac{1}{2}}  \Big\Vert\varphi \Big\Vert_{\mathrm{H}^{-\frac{1}{2},-\frac{1}{4}}\Big(\partial \Omega \times \mathbb{R}\Big)}.
	\end{align}
	Similar as in the proof of the third assertion of the Lemma \ref{3.7}, we conclude
	\begin{align}\label{operatoresti1}
		\Big\Vert \mathcal{S}_{\partial \Omega\times \mathbb{R}} \Big\Vert_{\mathcal{L}{\Big( \mathrm{H}^{-\frac{1}{2},-\frac{1}{4}}\Big(\partial \Omega\times \mathbb{R}\Big), \mathrm{H}^{\frac{1}{2},\frac{1}{4}}\Big(\partial \Omega\times \mathbb{R}\Big)}\Big) } \le \alpha_{\mathrm{p}}^{-\frac{1}{2}} \Big\Vert \overline{\mathcal{S}}_{\partial B\times \mathbb{R}} \Big\Vert_{\mathcal{L}{\Big( \mathrm{H}^{-\frac{1}{2},-\frac{1}{4}}\Big(\partial B\times \mathbb{R}\Big), \mathrm{H}^{\frac{1}{2},\frac{1}{4}}\Big(\partial B\times \mathbb{R}\Big)}\Big)}.
	\end{align}
	\vskip 0.1in
	\noindent
	Next, we use Lemma \ref{3.6}, Lemma \ref{4.8}, Corollary \ref{cc3}, i.e. the continuity of the operator $\gamma^{\textbf{int}}_{0}\overline{\mathcal{V}}: \mathrm{L}^2(\mathrm{B}\times\mathbb{R}) \to \mathrm{H}^{\frac{1}{2},\frac{1}{4}}\Big(\partial  B\times \mathbb{R}\Big)$, and due to the fact that $|\mathrm{E}|^2$ is independent of the time $t \in (0,T_0)$, we deduce
	\begin{align}\label{volesti}
		\nonumber
		\Big\Vert\gamma^{\textbf{int}}_{0}\mathcal{V}\Big[ |\mathrm{E}|^{2}\rchi_{(0,\mathrm{T_0})} \Big]\Big\Vert_{\mathrm{H}^{\frac{1}{2},\frac{1}{4}}\Big(\partial \Omega\times \mathbb{R}\Big)}
		&\leq 
		C\alpha_{\mathrm{p}}^{\frac{1}{4}} \delta^{3}\Big\Vert\gamma^{\textbf{int}}_{0}\overline{\mathcal{V}}\Big[ |\mathrm{E}|^{2}\rchi_{(0,\mathrm{T}_{0,\delta})} \Big]\Big\Vert_{\mathrm{H}^{\frac{1}{2},\frac{1}{4}}\Big(\partial  B\times \mathbb{R}\Big)}
		\\ \nonumber &\leq 
		C\alpha_{\mathrm{p}}^{\frac{1}{4}} \delta^{3}   \Big\Vert|\mathrm{E}|^{2}\rchi_{(0,\mathrm{T}_{0,\delta})}\Big\Vert_{\mathrm{L}^2\Big( B\times \mathbb{R}\Big)} \\&\leq 
		CT_0\alpha_{\mathrm{p}}^{\frac{1}{4}} \delta^2 \Big\Vert|\mathrm{E}|^{2}\Big\Vert_{\mathrm{L}^2\Big(\Omega\Big)}.
	\end{align}
	As we move on to the next step, it should be noted that the continuity of the hyper-singular operator $\mathbb{H}: \mathrm{H}^{\frac{1}{2},\frac{1}{4}}\Big(\partial \Omega\times \mathbb{R}\Big) \rightarrow \mathrm{H}^{-\frac{1}{2},-\frac{1}{4}}\Big(\partial \Omega\times \mathbb{R}\Big)$ is determined by the mapping property of the single and double layer operator, as well as the definition of Calderón projector. Now, if we recall the identity of the Calderón projector $\mathcal{C}$, we note that for $\varphi \in \mathrm{H}^{\frac{1}{2},\frac{1}{4}}\Big(\partial \Omega\times \mathbb{R}\Big)$, 
	\begin{align}\label{calderon}
		\mathcal{S}\mathbb{H}\Big[\varphi\Big] = (\frac{1}{2}I_{d} - \mathcal{K}) (\frac{1}{2}I_{d} + \mathcal{K})\Big[\varphi\Big].
	\end{align}
	Therefore, if we consider the appropriate norm, we deduce the following
	\begin{align}
		\Big\Vert \mathbb{H}[\varphi]\Big\Vert_{\mathrm{H}^{-\frac{1}{2},-\frac{1}{4}}\Big(\partial \Omega\times \mathbb{R}\Big)} &= \Big\Vert \mathcal{S}^{-1}_{\partial \Omega\times \mathbb{R}} \Big\Vert_{\mathcal{L}{\Big( \mathrm{H}^{\frac{1}{2},\frac{1}{4}}\Big(\partial \Omega\times \mathbb{R}\Big), \mathrm{H}^{-\frac{1}{2},-\frac{1}{4}}\Big(\partial \Omega\times \mathbb{R}\Big)}\Big)} \cdot \Big\Vert \frac{1}{2}I_{d} - \mathcal{K}\Big\Vert_{\mathcal{L}\Big(\mathrm{H}^{\frac{1}{2},\frac{1}{4}}(\partial \Omega\times \mathbb{R})\Big)} \\ &\cdot \Big\Vert\Big[\frac{1}{2}I_{d} + \mathcal{K}\Big](\varphi)\Big\Vert_{\mathrm{H}^{\frac{1}{2},\frac{1}{4}}\Big(\partial \Omega\times \mathbb{R}\Big)}.
	\end{align}
	\noindent 
	Now, Lemma \ref{3.7} and the non-scalability of  the double layer heat operator allows us to compute
	\begin{align}\label{hyperesti}
		\Big\Vert \mathbb{H}[\varphi]\Big\Vert_{\mathrm{H}^{-\frac{1}{2},-\frac{1}{4}}\Big(\partial \Omega\times \mathbb{R}\Big)} \lesssim \delta^{-1} \Big\Vert \varphi \Big\Vert_{\mathrm{H}^{\frac{1}{2},\frac{1}{4}}\Big(\partial \Omega\times \mathbb{R}\Big)}.
	\end{align}
	Additionally, the Steklov-Poincar\'e operator $\mathbb{A}^{\textbf{ext}} : \mathrm{H}^{\frac{1}{2},\frac{1}{4}}\Big(\partial \Omega\times \mathbb{R}\Big) \rightarrow \mathrm{H}^{-\frac{1}{2},-\frac{1}{4}}\Big(\partial \Omega\times \mathbb{R}\Big)$ for the heat equation is defined as follows:
	\begin{align}\label{aext}
		\mathbb{A}^{\textbf{ext}} := \mathbb{H}+ \Big(\frac{1}{2}I_{d} - \mathcal{K}^{*}\Big)\mathcal{S}^{-1} \Big(\frac{1}{2}I_{d} - \mathcal{K}\Big).
	\end{align}
	Accordingly, considering the bounds of the operators $\mathcal{S}^{-1}, \mathbb{H}, \mathcal{K}\ \text{and} \ \mathcal{K}^*,$ we can also conclude that the operator $\mathbb{A}^{\textbf{ext}}$ is  bounded. Next, we can see from the expression (\ref{aext}) that scaling depends on the operators $\mathcal{S}^{-1}$ and  $\mathbb{H}$ respectively, since the two other operators do not scale. 
	Therefore, we obtain the following from equation (\ref{hyperesti})
	\begin{align}\label{spoperator}
		\Big\Vert\mathbb{A}^{\textbf{ext}}\Big[\gamma^{\textbf{int}}_{0}\mathrm{U}_{\mathrm{i}}\Big]\Big\Vert_{\mathrm{H}^{-\frac{1}{2},-\frac{1}{4}}\Big(\partial \Omega\times \mathbb{R}\Big)} \lesssim \delta^{-1} \Big\Vert \gamma^{\textbf{int}}_{0}\mathrm{U}_{\mathrm{i}}\Big\Vert_{\mathrm{H}^{\frac{1}{2},\frac{1}{4}}\Big(\partial \Omega\times \mathbb{R}\Big)}.
	\end{align}
	We are now in a position to provide a priori estimate of $\gamma_{0}^{\text{int}}\mathrm{U}_{i}$. We consider (\ref{singlelayeresti1}), (\ref{volesti}), (\ref{hyperesti}), (\ref{spoperator}), Lemma \ref{3.2} and plugging these estimates into (\ref{gamma}) we obtain
	\begin{align}\label{final12}
		\Big\Vert\gamma^{\textbf{int}}_{0} \mathrm{U}_{\mathrm{i}} \Big\Vert_{\mathrm{H}^{\frac{1}{2},\frac{1}{4}}\Big(\partial \Omega\times \mathbb{R}\Big)} &\lesssim  \frac{\gamma_{\mathrm{m}}}{\gamma_{\mathrm{p}}} \alpha_{\mathrm{p}}^{-\frac{1}{2}}\delta^{-1}  \Big\Vert \gamma^{\textbf{int}}_{0}\mathrm{U}_{\mathrm{i}}\Big\Vert_{\mathrm{H}^{\frac{1}{2},\frac{1}{4}}\Big(\partial \Omega\times \mathbb{R}\Big)} +  \frac{\omega \cdot \boldsymbol{\Im}(\varepsilon_\mathrm{p})}{2\pi\gamma_{\mathrm{p}}}\alpha_{\mathrm{p}}^{\frac{1}{4}} \delta^2  \Big\Vert|\mathrm{E}|^{2}\Big\Vert_{\mathrm{L}^2\Big(\Omega\Big)},
	\end{align}
	and then
	\begin{align}
		\Big(1 - \frac{\gamma_{\mathrm{m}}}{\gamma_{\mathrm{p}}}\alpha_{\mathrm{p}}^{-\frac{1}{2}}\delta^{-1}\Big)\Big\Vert \gamma^{\textbf{int}}_{0}\mathrm{U}_{\mathrm{i}} \Big\Vert_{\mathrm{H}^{\frac{1}{2},\frac{1}{4}}\Big(\partial \Omega\times \mathbb{R}\Big)} \lesssim \frac{\omega \cdot \boldsymbol{\Im}(\varepsilon_\mathrm{p})}{2\pi\gamma_{\mathrm{p}}}\alpha_{\mathrm{p}}^{\frac{1}{4}} \delta^2 \Big\Vert|\mathrm{E}|^{2}\Big\Vert_{\mathrm{L}^2\Big(\Omega\Big)}.
	\end{align}
	Hence, due to the assumptions in (\ref{Assumption-heat-coefficients}), we have  
	\begin{equation}\label{condition-apprio-estimate}
		\frac{\gamma_{\mathrm{m}}}{\gamma_{\mathrm{p}}}\alpha_{\mathrm{p}}^{-\frac{1}{2}}\delta^{-1}<1,
	\end{equation} 
	for $\delta<<1$, which implies that
	\begin{align}\label{gamma0}
		\Big\Vert \gamma^{\textbf{int}}_{0}\mathrm{U}_{\mathrm{i}} \Big\Vert_{\mathrm{H}^{\frac{1}{2},\frac{1}{4}}\Big(\partial \Omega\times \mathbb{R}\Big)} \lesssim \frac{\omega \cdot \boldsymbol{\Im}(\varepsilon_\mathrm{p})}{2\pi\gamma_{\mathrm{p}}}\alpha_{\mathrm{p}}^{\frac{1}{4}} \delta^2 \Big\Vert|\mathrm{E}|^{2}\Big\Vert_{\mathrm{L}^2\Big(\Omega \Big)}.
	\end{align}
	Thus, it completes the proof of the first assertion of Proposition \ref{pr4}.
	\bigbreak
	\noindent
	Next, the a priori estimation of $\gamma^{\textbf{int}}_{1}\mathrm{U}_{\mathrm{i}}$ can be obtain with the help of (\ref{gamma0}). In order to do so, we consider the first integral equation in the system (\ref{IEfinal}), and obtain
	\begin{align}\label{maineq2}
		\nonumber
		\Big(\frac{1}{2}I_{d} - \mathcal{K}^{*}_{\alpha_{\mathrm{p}}}\Big)\Bigg[\gamma^{\textbf{int}}_{1}\mathrm{U}_{\mathrm{i}}\Bigg](\mathrm{x},t) = \mathbb{H}_{\alpha_{\mathrm{p}}}\Big[\gamma^{\textbf{int}}_{0}\mathrm{U}_{\mathrm{i}}\Big](\mathrm{x},t) + \frac{\omega \cdot \boldsymbol{\Im}(\varepsilon_\mathrm{p})}{2\pi\gamma_{\mathrm{p}}}\gamma^{\textbf{int}}_{1,\mathrm{x}}\mathcal{V}\Big[ |\mathrm{E}|^{2} \Big](\mathrm{x},t) \\ \text{i.e.} \quad \quad
		\gamma^{\textbf{int}}_{1} \mathrm{U}_{\mathrm{i}}= \Bigg[\frac{1}{2}I_{d} - \mathcal{K}^*_{\alpha_{\mathrm{p}}}\Bigg]^{-1}\Bigg[\mathbb{H}_{\alpha_\mathrm{p}}\Big[\gamma^{\textbf{int}}_{0}\mathrm{U}_{\mathrm{i}}\Big]\Bigg] + \frac{\omega \cdot \boldsymbol{\Im}(\varepsilon_\mathrm{p})}{2\pi\gamma_{\mathrm{p}}}\Bigg[\frac{1}{2}I_{d} - \mathcal{K}^{*}_{\alpha_{\mathrm{p}}}\Bigg]^{-1} \Bigg[ \gamma^{\textbf{int}}_{1,\mathrm{x}}\mathcal{V}\Big[ |\mathrm{E}|^{2} \Big]\Bigg]
	\end{align}
	from which we derive
	\begin{align}\label{gamma1h}
		\nonumber
		\Big\Vert  \gamma^{\textbf{int}}_{1}\mathrm{U}_{\mathrm{i}} \Big\Vert_{\mathrm{H}^{-\frac{1}{2},-\frac{1}{4}}\Big(\partial \Omega\times \mathbb{R}\Big)} &=  \Big\Vert\Big[\frac{1}{2}I_{d} - \mathcal{K}^*_{\alpha_{\mathrm{p}}}\Big]^{-1}\Big\Vert_{\mathcal{L}\Big(\mathrm{H}^{-\frac{1}{2},-\frac{1}{4}}(\partial \Omega\times \mathbb{R})\Big)} \Big\Vert \mathbb{H}_{\alpha_\mathrm{p}}\Big[\gamma^{\textbf{int}}_{0}\mathrm{U}_{\mathrm{i}}\Big] \Big\Vert_{\mathrm{H}^{-\frac{1}{2},-\frac{1}{4}}\Big(\partial \Omega\times \mathbb{R}\Big)}  \\&+ \frac{\omega \cdot \boldsymbol{\Im}(\varepsilon_\mathrm{p})}{2\pi\gamma_{\mathrm{p}}} \Big\Vert\Big[\frac{1}{2}I_{d} - \mathcal{K}^*_{\alpha_\mathrm{p}}\Big]^{-1}\Big\Vert_{\mathcal{L}\Big(\mathrm{H}^{-\frac{1}{2},-\frac{1}{4}}(\partial \Omega\times \mathbb{R)}\Big)} \Big\Vert \gamma^{\textbf{int}}_{1,\mathrm{x}}\mathcal{V}\Big[ |\mathrm{E}|^{2} \Big] \Big\Vert_{\mathrm{H}^{-\frac{1}{2},-\frac{1}{4}}\Big(\partial \Omega\times \mathbb{R}\Big)}.
	\end{align}
	Then, we deduce from Lemma \ref{3.6}, Lemma \ref{3.11}, Corollary \ref{c4} i.e. the continuity of the operator $\gamma^{\textbf{int}}_{1,\xi}\overline{\mathcal{V}}: \mathrm{L}^2(\mathrm{B}\times\mathbb{R}) \to \mathrm{H}^{-\frac{1}{2},-\frac{1}{4}}\Big(\partial  B\times \mathbb{R}\Big)$, and observing the fact $|\mathrm{E}|^2$ is independent of the time t in $(0,\mathrm{T}_0)$
	\begin{align}\label{volesti1}
		\nonumber
		\Big\Vert\gamma^{\textbf{int}}_{1,\mathrm{x}} \mathcal{V}\Big[ |\mathrm{E}|^{2}\rchi_{(0,\mathrm{T}_{0})} \Big]\Big\Vert_{\mathrm{H}^{-\frac{1}{2},-\frac{1}{4}}\Big(\partial \Omega\times \mathbb{R}\Big)}
		&\leq \nonumber
		\alpha^{\frac{1}{2}}_{\mathrm{p}} \delta^{\frac{5}{2}}\Big\Vert\gamma^{\textbf{int}}_{1,\xi} \overline{\mathcal{V}}\Big[ |\mathrm{E}|^{2}\rchi_{(0,\mathrm{T}_{0,\delta})} \Big]\Big\Vert_{\mathrm{H}^{-\frac{1}{2},-\frac{1}{4}}\Big(\partial  B\times \mathbb{R}\Big)}
		\\ &\lesssim \nonumber
		\alpha^{\frac{1}{2}}_{\mathrm{p}} \delta^{\frac{5}{2}}   \Big\Vert|\mathrm{E}|^{2}\rchi_{(0,\mathrm{T}_{0,\delta})}\Big\Vert_{\mathrm{L}^2\Big( B\times \mathbb{R}\Big)} \\&\lesssim\mathrm{T}_0\alpha^{\frac{1}{2}}_{\mathrm{p}}\delta^{\frac{3}{2}}\Big\Vert|\mathrm{E}|^{2}\Big\Vert_{\mathrm{L}^2\Big(\Omega\Big)}.
	\end{align}
	Combining (\ref{gamma1h}) and (\ref{volesti1}), we deduce that
	\begin{align}\label{final2}
		\Big\Vert \gamma^{\textbf{int}}_{1} \mathrm{U}_{\mathrm{i}} \Big\Vert_{\mathrm{H}^{-\frac{1}{2},-\frac{1}{4}}\Big(\partial \Omega\times \mathbb{R}\Big)} &\lesssim \delta^{-1}\Big\Vert \gamma^{\textbf{int}}_{0}\mathrm{U}_{\mathrm{i}}\Big\Vert_{\mathrm{H}^{\frac{1}{2},\frac{1}{4}}\Big(\partial \Omega\times \mathbb{R}\Big)} +  \frac{\omega \cdot \boldsymbol{\Im}(\varepsilon_\mathrm{p})}{2\pi\gamma_{\mathrm{p}}}\alpha^{\frac{1}{2}}_{\mathrm{p}}\delta^{\frac{3}{2}} \Big\Vert|\mathrm{E}|^{2}\Big\Vert_{\mathrm{L}^2\Big(\Omega\Big)}.
	\end{align}
	Consequently, after plugging (\ref{gamma0}) in the previous expression, we obtain
	\begin{align}\label{gamma1}
		\Big\Vert \gamma^{\textbf{int}}_{1}\mathrm{U}_{\mathrm{i}} \Big\Vert_{\mathrm{H}^{-\frac{1}{2},-\frac{1}{4}}\Big(\partial \Omega\times \mathbb{R}\Big)} \lesssim \frac{\omega \cdot \boldsymbol{\Im}(\varepsilon_\mathrm{p})}{2\pi\gamma_{\mathrm{p}}}\Big(\alpha^{\frac{1}{4}}_\mathrm{p}\delta + \alpha^{\frac{1}{2}}_\mathrm{\mathrm{p}}\delta^{\frac{3}{2}}\Big) \Big\Vert|\mathrm{E}|^{2}\Big\Vert_{\mathrm{L}^2\Big(\Omega\Big)}.
	\end{align}
	Therefore, it completes the proof of second assertion of Proposition \ref{pr4}.
	\subsection{Estimation of the Heat Potential's Dominating Term}
	The heat propagation outside and close to the nanoparticle can be estimated using the following relation:
	\begin{align}\label{1stIE}
		\mathrm{U}_{\mathrm{e}}(\xi,t)&=   -\mathcal{S}\Big[\gamma^{\textbf{ext}}_{1}\mathrm{U}_{\mathrm{e}}\Big](\xi,t) + \mathcal{D}\Big[\gamma^{\textbf{ext}}_{0}\mathrm{U}_{\mathrm{e}}\Big](\xi,t),\quad \quad \text{for} \ (\xi,t) \in \mathbb{R}^2\setminus\overline{\Omega} \times (0,\mathrm{T}),
	\end{align}
	Rewriting the above expression explicitly with the transmission condition as (\ref{eq:heat1}), we obtain
	\begin{align}
		\mathrm{U}_{\mathrm{e}}(\xi,t) \nonumber &= -\frac{\gamma_{\mathrm{p}}}{\gamma_{\mathrm{m}}}\frac{1}{\alpha_\mathrm{m}}\int_{0}^t\int_{\partial\Omega}\Phi^{\textbf{e}}(\xi,t;\mathrm{y},\tau)\gamma^{\textbf{ext}}_{1}\mathrm{U}_{\mathrm{e}}(\mathrm{y},\tau)d\sigma_{\mathrm{y}}d\tau + \frac{1}{\alpha_\mathrm{m}}\int_{0}^t\int_{\partial\Omega}\gamma^{\textbf{ext}}_{1,\mathrm{y}}\Phi^{\textbf{e}}(\xi,t;\mathrm{y},\tau)\gamma^{\textbf{ext}}_{0}\mathrm{U}_{\mathrm{e}}(\mathrm{y},\tau)d\sigma_{\mathrm{y}}d\tau \\ \nonumber &= -\frac{\gamma_{\mathrm{p}}}{\gamma_{\mathrm{m}}}\frac{1}{\alpha_\mathrm{m}} \int_{0}^t\int_{\partial\Omega}\Phi^{\textbf{e}}(\xi,t;\mathrm{z},\tau)\gamma^{\textbf{int}}_{1}\mathrm{U}_{\mathrm{i}}(\mathrm{y},\tau)d\sigma_\mathrm{y}d\tau \\ \nonumber &- \frac{\gamma_{\mathrm{p}}}{\gamma_{\mathrm{m}}}\frac{1}{\alpha_\mathrm{m}}\int_{0}^t\int_{\partial\Omega}\Big[\Phi^{\textbf{e}}(\xi,t;\mathrm{z},\tau)-\Phi^{\textbf{e}}(\xi,t;\mathrm{y},\tau)\Big]\gamma^{\textbf{int}}_{1}\mathrm{U}_{\mathrm{i}}(\mathrm{y},\tau)d\sigma_{\mathrm{y}}d\tau \\ \nonumber &+ \frac{1}{\alpha_\mathrm{m}}\int_{0}^t\int_{\partial\Omega}\gamma^{\textbf{int}}_{1,\mathrm{y}}\Phi^{\textbf{e}}(\xi,t;\mathrm{y},\tau)\gamma^{\textbf{int}}_{0}\mathrm{U}_{\mathrm{i}}(\mathrm{y},\tau)d\sigma_{\mathrm{y}}d\tau.
	\end{align}
	Now, noticing the fact that, for fixed $(\xi,t) \in \mathbb{R}^2\backslash\overline{\Omega} \times (0,\mathrm{T})$, the function $\Phi^{\textbf{e}}(\xi,t;\mathrm{y},\tau)$ is sufficiently smooth with respect to $(\mathrm{y},\tau) \in \partial\Omega\times(0,t)$, it follows from Taylor's series expansion for $\mathrm{z} \in \Omega$  and by duality pairing that
	\begin{align}
		\nonumber
		\textbf{err}^{(1)} &:= \Bigg|\int_{0}^t\int_{\partial\Omega}\Big[\Phi^{\textbf{e}}(\xi,t;\mathrm{z},\tau)-\Phi^{\textbf{e}}(\xi,t;\mathrm{y},\tau)\Big]\gamma^{\textbf{int}}_{1}\mathrm{U}_{\mathrm{i}}(\mathrm{y},\tau)d\sigma_{\mathrm{y}}d\tau \Bigg| \\ \nonumber &\lesssim \mathcal{O} \Bigg(\delta \Big\Vert\gamma^{\textbf{int}}_{1}\mathrm{U}_{\mathrm{i}} \Big\Vert_{\mathrm{H}^{-\frac{1}{2},-\frac{1}{4}}\Big(\partial \Omega\times (0,\mathrm{t})\Big)}\Big\Vert \nabla \Phi^{\textbf{e}}(\xi,t;\mathrm{z},\cdot) \Big\Vert_{\mathrm{H}^{\frac{1}{2},\frac{1}{4}}\Big(\partial \Omega\times (0,t)\Big)}\Bigg)
		\\\nonumber &\lesssim \mathcal{O} \Bigg(\delta \Big\Vert\gamma^{\textbf{int}}_{1}\mathrm{U}_{\mathrm{i}} \Big\Vert_{\mathrm{H}^{-\frac{1}{2},-\frac{1}{4}}\Big(\partial \Omega\times \mathbb{R}\Big)}\Big\Vert \nabla \Phi^{\textbf{e}}(\xi,t;\mathrm{z},\cdot) \Big\Vert_{\mathrm{H}^{\frac{1}{2},\frac{1}{4}}\Big(\partial \Omega\times (0,t)\Big)}\Bigg).
	\end{align}
	Therefore, we deduce that 
	\begin{align}
		\nonumber
		\textbf{err}^{(1)} = \mathcal{O} \Bigg(\delta \Big\Vert\gamma^{\textbf{int}}_{1}\mathrm{U}_{\mathrm{i}} \Big\Vert_{\mathrm{H}^{-\frac{1}{2},-\frac{1}{4}}\Big(\partial \Omega\times\mathbb{R}\Big)}\Big\Vert \nabla \Phi^{\textbf{e}}(\xi,t;\mathrm{z},\cdot) \Big\Vert_{\mathrm{H}^{\frac{1}{2},\frac{1}{4}}\Big( \partial \Omega\times(0,\mathrm{t})\Big)}\Bigg).
	\end{align}
	Then, applying interpolation theory and observing that $\nabla \Phi^{\textbf{e}}(\xi,t;\mathrm{z},\cdot)$ are independent of space variables, i.e. $\xi$ is outside $\overline{\Omega}$, we derive the following:
	\begin{align}\label{phi-}
		\nonumber
		\textbf{err}^{(1)} &= \mathcal{O} \Bigg(\delta \Big\Vert\gamma^{\textbf{int}}_{1}\mathrm{U}_{\mathrm{i}} \Big\Vert_{\mathrm{H}^{-\frac{1}{2},-\frac{1}{4}}\Big(\partial \Omega\times \mathbb{R}\Big)}\Big\Vert \nabla \Phi^{\textbf{e}}(\xi,t;\mathrm{z},\cdot) \Big\Vert^{\frac{1}{2}}_{\mathrm{L}^2\Big(\partial \Omega\times (0,\mathrm{t})\Big)}\Big\Vert \nabla \Phi^{\textbf{e}}(\xi,t;\mathrm{z},\cdot) \Big\Vert^{\frac{1}{2}}_{\mathrm{H}^{1,\frac{1}{2}}\Big(\partial \Omega\times (0,t)\Big)}\Bigg)
		\\ &=\mathcal{O} \Bigg(\delta^{\frac{3}{2}} \Big\Vert\gamma^{\textbf{int}}_{1}\mathrm{U}_{\mathrm{i}} \Big\Vert_{\mathrm{H}^{-\frac{1}{2},-\frac{1}{4}}\Big(\partial \Omega\times \mathbb{R}\Big)}\Big\Vert \nabla \Phi^{\textbf{e}}(\xi,t;\mathrm{z},\cdot) \Big\Vert^{\frac{1}{2}}_{\mathrm{L}^2(0,\mathrm{t})}\Big\Vert \nabla \Phi^{\textbf{e}}(\xi,t;\mathrm{z},\cdot) \Big\Vert^{\frac{1}{2}}_{\mathrm{H}^{1,\frac{1}{2}}(0,\mathrm{t})}\Bigg).
	\end{align}
	Let us estimate the integral 
	\begin{equation}
		\nonumber
		\textbf{err}^{(2)} := \displaystyle\frac{1}{\alpha_\mathrm{m}}\int_{0}^t\int_{\partial\Omega}\gamma^{\textbf{int}}_{1,\mathrm{y}}\Phi^{\textbf{e}}(\xi,t;\mathrm{y},\tau)\gamma^{\textbf{int}}_{0}\mathrm{U}_{\mathrm{i}}(\mathrm{y},\tau)d\sigma_{\mathrm{y}}d\tau.
	\end{equation}
	As $\xi \in \mathbb{R}^2\setminus \overline{\Omega}$, due to duality between the function space $\mathrm{L}^2\Big(\partial \Omega\times (0,\mathrm{t)}\Big)$, based on the standard Sobolev embedding $ \mathrm{H}^{-\frac{1}{2},-\frac{1}{4}}\Big(\partial  \mathrm{B}\times \mathbb{R}\Big) \xhookrightarrow{\mathrm{i}} \mathrm{L}^2\Big(\partial \mathrm{B}\times\mathbb{R}\Big)$, and using Lemma \ref{3.6} we obtain the following estimate:
	\begin{align}
		\nonumber
		\textbf{err}^{(2)} &:= \Bigg|\displaystyle\frac{1}{\alpha_\mathrm{m}}\int_{0}^t\int_{\partial\Omega}\gamma^{\textbf{int}}_{1,\mathrm{y}}\Phi^{\textbf{e}}(\xi,t;\mathrm{y},\tau)\gamma^{\textbf{int}}_{0}\mathrm{U}_{\mathrm{i}}(\mathrm{y},\tau)d\sigma_{\mathrm{y}}d\tau \Bigg| 
		\\ \nonumber &\lesssim
		\Big\Vert \gamma^{\textbf{int}}_{0}\mathrm{U}_{\mathrm{i}}\Big\Vert_{\mathrm{L}^2\Big(\partial \Omega \times (0,\mathrm{t})\Big)} \Big\Vert\gamma^{\textbf{int}}_{1,\mathrm{y}}\Phi^{\textbf{e}}(\xi,t;\mathrm{y},\cdot)\Big\Vert_{\mathrm{L}^2\Big(\partial \Omega \times (0,\mathrm{t})\Big)}
		\\ \nonumber &\lesssim 
		\Big\Vert \gamma^{\textbf{int}}_{0}\mathrm{U}_{\mathrm{i}}\Big\Vert_{\mathrm{L}^2\Big(\partial \Omega \times\mathbb{R}\Big)} \Big\Vert\gamma^{\textbf{int}}_{1,\mathrm{y}}\Phi^{\textbf{e}}(\xi,t;\mathrm{y},\cdot)\Big\Vert_{\mathrm{L}^2\Big(\partial \Omega \times(0,\mathrm{t})\Big)}
		\\ \nonumber &\lesssim 
		\alpha^{\frac{1}{4}}_\mathrm{p}\delta^{\frac{3}{2}}
		\Big\Vert \gamma^{\textbf{int}}_{0}\mathrm{U}_{\mathrm{i}}\Big\Vert_{\mathrm{H}^{\frac{1}{2},\frac{1}{4}}\Big(\partial \mathrm{B}\times\mathbb{R}\Big)} \alpha^{\frac{1}{2}}_{\mathrm{m}}\delta^{\frac{3}{2}}\Big\Vert\gamma^{\textbf{int}}_{1,\mathrm{y}}\Phi^{\textbf{e}}(\xi,t;\mathrm{y},\cdot)\Big\Vert_{\mathrm{L}^2\Big(\partial \Omega \times(0,\mathrm{t})\Big)}
		\\ \nonumber&\lesssim 
		\Big\Vert \gamma^{\textbf{int}}_{0}\mathrm{U}_{\mathrm{i}}\Big\Vert_{\mathrm{H}^{\frac{1}{2},\frac{1}{4}}\Big(\partial \Omega \times \mathbb{R}\Big)} \Big\Vert \gamma^{\textbf{int}}_{1,\mathrm{y}}\Phi^{\textbf{e}}(\xi,t;\mathrm{y},\cdot)\Big\Vert_{\mathrm{L}^2\Big(\partial \Omega \times(0,\mathrm{t})\Big)}.
	\end{align}
	Analogously, we deduce that
	\begin{align}\label{main1}
		\textbf{err}^{(2)} := \mathcal{O} \Bigg(\Big\Vert \gamma^{\textbf{int}}_{0}\mathrm{U}_{\mathrm{i}}\Big\Vert_{\mathrm{H}^{\frac{1}{2},\frac{1}{4}}\Big(\partial \Omega \times \mathbb{R}\Big)} \Big\Vert \gamma^{\textbf{int}}_{1,\mathrm{y}}\Phi^{\textbf{e}}(\xi,t;\mathrm{y},\cdot)\Big\Vert_{\mathrm{L}^2\Big(\partial \Omega \times (0,\mathrm{t})\Big)}\Bigg).
	\end{align}
	\noindent
	Consequently, due to (\ref{phi-}), (\ref{main1}), we derive the explicit form of $\mathrm{U}_{\mathrm{e}}(\xi,t) \ \text{for}\ (\xi,t) \in \mathbb{R}^2\setminus\Omega \times (0,\mathrm{T})$
	\begin{align}\label{mainformula1}
		\mathrm{U}_{\mathrm{e}}(\xi,\mathrm{t}) =-\frac{\gamma_{\mathrm{p}}}{\gamma_{\mathrm{m}}}\frac{1}{\alpha_\mathrm{m}} \int_{0}^t\int_{\partial\Omega}\Phi(\xi,t;\mathrm{z},\tau)\gamma^{\textbf{int}}_{1}\mathrm{U}_{\mathrm{i}}(\mathrm{y},\tau)d\sigma_\mathrm{y}d\tau + \textbf{err}^{(1)} + \textbf{err}^{(2)}.
	\end{align}
	To derive an explicit form of the first term in the approximation of $ \mathrm{U}_{\mathrm{e}}(\xi,t)$, we start with the first integral equation of the system (\ref{IEfinal})
	\begin{align}\label{maineq2}
		\Bigg(\frac{1}{2}I_{d} - \mathcal{K}^{*}_{\alpha_\mathrm{p}}\Bigg)\Bigg[\gamma^{\textbf{int}}_{1}\mathrm{U}_{\mathrm{i}}\Bigg] = \mathbb{H}_{\alpha_\mathrm{p }}\Big[\gamma^{\textbf{int}}_{0}\mathrm{U}_{\mathrm{i}}\Big] + \frac{\omega \cdot \boldsymbol{\Im}(\varepsilon_\mathrm{p})}{2\pi\gamma_{\mathrm{p}}}\gamma^{\textbf{int}}_{1,\mathrm{x}}\mathcal{V}\Big[ |\mathrm{E}|^{2} \Big].
	\end{align}
	We multiply this integral equation with $\Phi^{\textbf{e}}(\xi,t;z,\tau)$ and integrate with respect to the spatial boundary $\partial\Omega$ to obtain the following:
	\begin{align}
		\nonumber
		&\frac{1}{2}\int_{\partial\Omega}\Phi^{\textbf{e}}(\xi,t;z,\tau)\gamma^{\textbf{int}}_{1}\mathrm{U}_{\mathrm{i}}(\mathrm{y},\tau)d\sigma_\mathrm{y}- \int_{\partial\Omega} \mathcal{K}^*_\text{Lap}\Big[\gamma^{\textbf{int}}_{1}\mathrm{U}_{\mathrm{i}}\Big]\Phi^{\textbf{e}}(\xi,t;z,\tau)d\sigma_\mathrm{y} \\ \nonumber&- \int_{\partial\Omega}\Big(\mathcal{K}^{*}_{\alpha_\mathrm{p}}-\mathcal{K}^*_\text{Lap}\Big)\Big[\gamma^{\textbf{int}}_{1}\mathrm{U}_{\mathrm{i}}\Big]\Phi^{\textbf{e}}(\xi,t;z,\tau)d\sigma_\mathrm{y} = \int_{\partial\Omega}\Phi^{\textbf{e}}(\xi,t;z,\tau) \mathbb{H}_{\alpha_\mathrm{p }}\Big[\gamma^{\textbf{int}}_{0}\mathrm{U}_{\mathrm{i}}\Big](\mathrm{x},\mathrm{t}) d\sigma_\mathrm{y} \\&+ \frac{\omega \cdot \boldsymbol{\Im}(\varepsilon_\mathrm{p})}{2\pi\gamma_{\mathrm{p}}} \int_{\partial\Omega}\Phi^{\textbf{e}}(\xi,t;z,\tau) \gamma^{\textbf{int}}_{1,\mathrm{y}}\mathcal{V}\Big[ |\mathrm{E}|^{2} \Big](\mathrm{y},\tau) d\sigma_\mathrm{y}
	\end{align}
	where we denote by $\mathcal{K}^*_\text{Lap}$ the spatial adjoint of the double layer operator of the Laplacian (i.e. the Neumann-Poincar\'e operator). Then, as $\Phi^{\textbf{e}}(\xi,t;z,\tau)$ is a constant, in the the space variable 'y', then, refer to \cite{steinbuch}, $\mathcal{K}_{\textbf{Lap}}\Big[\Phi^{\textbf{e}}(\xi,t;z,\tau)\Big](y) = -\frac{1}{2}\Phi^{\textbf{e}}(\xi,t;z,\tau)$, i.e. $-\frac{1}{2}$ is an eigenvalue of the Neumman-Poincar\'e operator. This allows us to rewrite the above equation after integrating in time domain as follows:
	\begin{align}
		\nonumber
		&\int_{0}^{t}\int_{\partial\Omega}\Phi^{\textbf{e}}(\xi,t;z,\tau)\gamma^{\textbf{int}}_{1}\mathrm{U}_{\mathrm{i}}(\mathrm{y},\tau)d\sigma_\mathrm{y}d\tau - \int_{0}^{t}\int_{\partial\Omega}\gamma^{\textbf{int}}_{1}\mathrm{U}_{\mathrm{i}}(\mathrm{y},\tau)\Big(\mathcal{K}_{\alpha_\mathrm{p}}-\mathcal{K}_\text{Lap}\Big)\Big[\Phi^{\textbf{e}}(\xi,t;z,\cdot)\Big](\mathrm{y},\tau)d\sigma_\mathrm{y}d\tau \\ \nonumber &= \frac{1}{\alpha_\mathrm{p}}\int_{0}^{t}\int_{\partial\Omega}\Phi^{\textbf{e}}(\xi,t;z,\tau)  \gamma^{\textbf{int}}_{1,\mathrm{y}}\int_{0}^\tau\int_{\partial\Omega}\gamma^{\textbf{int}}_{1,\mathrm{v}}\Phi (\mathrm{y},\tau;\mathrm{v},\mathrm{s}) \gamma^{\textbf{int}}_{0}\mathrm{U}_{\mathrm{i}}(\mathrm{v},\mathrm{s})d\sigma_\mathrm{v}d\mathrm{s} d\sigma_\mathrm{y}d\tau \\&+ \frac{1}{\alpha_\mathrm{p}}\frac{\omega \cdot \boldsymbol{\Im}(\varepsilon_\mathrm{p})}{2\pi\gamma_{\mathrm{p}}} \int_{0}^{t}\int_{\partial\Omega}\Phi^{\textbf{e}}(\xi,t;z,\tau) \gamma^{\textbf{int}}_{1,\mathrm{y}}\int_{0}^{\tau}\int_{\Omega}\Phi(\mathrm{y},\tau;\mathrm{v},\mathrm{s})|\mathrm{E}|^{2}(\mathrm{v}) d\mathrm{v}d\mathrm{s} d\sigma_\mathrm{y}d\tau.
	\end{align}
	Now, we can represent the term $\Big(\mathcal{K}_{\alpha_\mathrm{p}}-\mathcal{K}_\text{Lap}\Big)\Big[\Phi^{\textbf{e}}(\xi,t;z,\cdot)\Big]$ as
	\begin{align}
		\Big(\mathcal{K}_{\alpha_\mathrm{p}}-\mathcal{K}_\text{Lap}\Big)\Big[\Phi^{\textbf{e}}(\xi,t;z,\cdot)\Big](\mathrm{y},\tau) = \int_{\partial\Omega} \frac{(\mathrm{y}-\mathrm{v})\cdot\nu_{\mathrm{v}}}{2\pi|\mathrm{y}-\mathrm{v}|^2}\Big[\varphi(\mathrm{v},\mathrm{y},\mathrm{t},\tau)-\Phi^{\textbf{e}}(\xi,t;z,\tau)\Big]d\sigma_\mathrm{v}
	\end{align}
	where we set
	\begin{equation}\label{defvarphi}
		\varphi(\mathrm{v},\mathrm{y},t, \tau) := \displaystyle\int_{0}^{\tau}\dfrac{\alpha |\mathrm{y}-\mathrm{v}| ^2}{4(\mathrm{s}-\tau)^2}\textbf{exp}\Big(-\dfrac{\alpha|\mathrm{y}-\mathrm{v}| ^2}{4(\mathrm{s}-\tau)}\Big)\Phi^{\textbf{e}}(\xi,\mathrm{t};\mathrm{z},\mathrm{s})d\mathrm{s}.
	\end{equation}
	Consequently, we derive
	\begin{align}\label{sinires1}
		\nonumber
		\textbf{Expression} &:=
		\int_{0}^{t}\int_{\partial\Omega}\Phi^{\textbf{e}}(\xi,t;z,\tau)\gamma^{\textbf{int}}_{1}\mathrm{U}_{\mathrm{i}}(\mathrm{y},\tau)d\sigma_\mathrm{y}d\tau \\ \nonumber &= \frac{1}{\alpha_\mathrm{p}}\frac{\omega \cdot \boldsymbol{\Im}(\varepsilon_\mathrm{p})}{2\pi\gamma_{\mathrm{p}}} \int_{0}^{t}\int_{\partial\Omega}\Phi^{\textbf{e}}(\xi,t;z,\tau)\gamma^{\textbf{int}}_{1,\mathrm{y}}\int_{0}^{\tau}\int_{\Omega}\Phi(\mathrm{y},\tau; \mathrm{v},\mathrm{s})|\mathrm{E}|^{2}(\mathrm{v}) d\mathrm{v}d\mathrm{s} d\sigma_\mathrm{y}d\tau \nonumber\\&+ \int_{0}^{t}\int_{\partial\Omega}\gamma^{\textbf{int}}_{1}\mathrm{U}_{\mathrm{i}}(\mathrm{y},\tau)\int_{\partial\Omega} \frac{(\mathrm{y}-\mathrm{v})\cdot\nu_{\mathrm{v}}}{2\pi|\mathrm{y}-\mathrm{v}|^2}\Big[\varphi(\mathrm{v},\mathrm{y},\mathrm{t},\tau)-\Phi^{\textbf{e}}(\xi,t;z,\tau)\Big]d\sigma_\mathrm{v} d\sigma_\mathrm{y} d\tau \nonumber \\&+ \Big\langle \Phi^{\textbf{e}}(\xi,\mathrm{t};z,\mathrm{\tau)},\mathbb{H}_{\alpha_\mathrm{p}}\Big[\gamma^{\textbf{int}}_{0}\mathrm{T}_{\mathrm{i}}\Big]\Big\rangle_{\partial\Omega\times(0,\mathrm{t})}.
	\end{align}
	In the following expression, we apply integration by parts, Leibniz's rule for integration, change of variables and observe the fact that the Newtonian heat potential satisfies the corresponding in-homogeneous equation $(\alpha\partial_\mathrm{t}-\Delta)\mathcal{V}(\mathrm{x},\mathrm{t}) = f(\mathrm{x},\mathrm{t})$ to rewrite the leading term of (\ref{sinires1}) as follows: 
	\begin{align}\label{volumeesti1}
		\nonumber
		\textbf{First Term} &:= \frac{1}{\alpha_\mathrm{p}}\frac{\omega \cdot \boldsymbol{\Im}(\varepsilon_\mathrm{p})}{2\pi\gamma_{\mathrm{p}}} \int_{0}^{t}\int_{\partial\Omega}\Phi^{\textbf{e}}(\xi,\mathrm{t};\mathrm{z},\tau) \gamma^{\textbf{int}}_{1,\mathrm{y}}\int_{0}^{\tau}\int_{\Omega}\Phi(\mathrm{y},\tau; \mathrm{v},\mathrm{s})|\mathrm{E}|^{2}(\mathrm{v}) d\mathrm{v}d\mathrm{s} d\sigma_\mathrm{y}d\tau
		\nonumber 
		\\ &=\frac{\omega \cdot \boldsymbol{\Im}(\varepsilon_\mathrm{p})}{2\pi\gamma_{\mathrm{p}}}\frac{1}{\alpha_\mathrm{p}} \int_{0}^{t}\Phi^{\textbf{e}}(\xi,t;\mathrm{z},\tau)\int_{\partial\Omega}\gamma^{\textbf{int}}_{1,\mathrm{y}}\int_{0}^{\tau}\int_{\Omega}\Phi(\mathrm{y},\tau;\mathrm{v},\mathrm{s})|\mathrm{E}|^{2}(\mathrm{v}) d\mathrm{v}d\mathrm{s} d\sigma_\mathrm{y}d\tau
		\nonumber \\ &= \int_{0}^{t}\Phi^{\textbf{e}}(\xi,t;\mathrm{z},\tau)\int_{\Omega}\Delta \frac{1}{\alpha_\mathrm{p}}\int_{0}^{\tau}\int_{\Omega}\Phi(\mathrm{y},\tau;\mathrm{v},\mathrm{s})\frac{\omega \cdot \boldsymbol{\Im}(\varepsilon_\mathrm{p})}{2\pi\gamma_{\mathrm{p}}}|\mathrm{E}|^{2}(\mathrm{v}) d\mathrm{v}d\mathrm{s} d\mathrm{y}d\tau
		\nonumber       \\ &= \int_{0}^{t}\Phi^{\textbf{e}}(\xi,t;z,\tau)\int_{\Omega}\alpha_\mathrm{p}\partial_{\tau}\frac{1}{\alpha_\mathrm{p}}\int_{0}^{\tau}\int_{\Omega}\Phi(\mathrm{y},\tau;\mathrm{v},\mathrm{s})\frac{\omega \cdot \boldsymbol{\Im}(\varepsilon_\mathrm{p})}{2\pi\gamma_{\mathrm{p}}}|\mathrm{E}|^{2}(\mathrm{v}) d\mathrm{v}d\mathrm{s} d\mathrm{y}d\tau
		\nonumber       \\ &- \frac{\omega \cdot \boldsymbol{\Im}(\varepsilon_\mathrm{p})}{2\pi\gamma_{\mathrm{p}}}\int_{0}^{t}\int_{\Omega} |\mathrm{E}|^{2}(\mathrm{y})\Phi^{\textbf{e}}(\xi,t;\mathrm{z},\tau) d\mathrm{y}d\tau 
		\nonumber \\ &= \int_{0}^{t}\int_{\Omega}\Phi^{\textbf{e}}(\xi,t;\mathrm{z},\tau)\int_{0}^{\tau}\int_{\Omega}\Phi(\mathrm{v},\mathrm{s})\frac{\omega \cdot \boldsymbol{\Im}(\varepsilon_\mathrm{p})}{2\pi\gamma_{\mathrm{p}}}\partial_{\tau}|\mathrm{E}|^{2}(\mathrm{v}-\mathrm{y}) d\mathrm{v}d\mathrm{s} d\mathrm{y}d\tau
		\nonumber       \\ &+ \frac{\omega \cdot \boldsymbol{\Im}(\varepsilon_\mathrm{p})}{2\pi\gamma_{\mathrm{p}}} \int_{0}^{t}\int_{\Omega}\Phi^{\textbf{e}}(\xi,t;\mathrm{z},\tau)\int_{\Omega}\Phi(\mathrm{v},\tau)|\mathrm{E}|^{2}(\mathrm{v}-\mathrm{y}) d\mathrm{v}d\mathrm{y}d\tau
		\nonumber       \\ &-
		-\frac{\omega \cdot \boldsymbol{\Im}(\varepsilon_\mathrm{p})}{2\pi\gamma_{\mathrm{p}}}\int_{0}^{t} \int_{\Omega}\Phi^{\textbf{e}}(\xi,t;\mathrm{z},\tau) |\mathrm{E}|^{2}(\mathrm{y}) d\mathrm{y} d\tau
		\nonumber       \\ &= -\frac{\omega \cdot \boldsymbol{\Im}(\varepsilon_\mathrm{p})}{2\pi\gamma_{\mathrm{p}}}\Bigg[\int_{0}^{t} \int_{\Omega}\Phi^{\textbf{e}}(\xi,t;\mathrm{z},\tau) |\mathrm{E}|^{2}(\mathrm{y}) d\mathrm{y} d\tau - 
		\int_{0}^{t}\int_{\Omega}\Phi^{\textbf{e}}(\xi,t;\mathrm{z},\tau)\int_{\Omega}\Phi(\mathrm{y},\tau;\mathrm{v})|\mathrm{E}|^{2}(\mathrm{v}) d\mathrm{v}d\mathrm{y}d\tau \Bigg].
	\end{align}
	Furthermore, we consider the second term of the previous expression \textbf{'First Term'} (\ref{volumeesti1}), i.e.
	\begin{align}\label{initial}
		\frac{\omega \cdot \boldsymbol{\Im}(\varepsilon_\mathrm{p})}{2\pi\gamma_{\mathrm{p}}}\int_{0}^{t}\int_{\Omega}\Phi^{\textbf{e}}(\xi,t;z,\tau)\int_{\Omega}\Phi(\mathrm{y},\tau;\mathrm{v})|\mathrm{E}|^{2}(\mathrm{v}) d\mathrm{v}d\mathrm{y}d\tau
	\end{align}
	and derive its estimate. To do that let us recall the  initial heat potential for $|\mathrm{E}|^2 \in \mathrm{L}^2\Big(\Omega \Big)$,
	$\displaystyle\mathbb{I}\Big[|\mathrm{E}|^2\Big](\mathrm{y},\tau)= \int_{\Omega}\Phi(\mathrm{y},\tau;\mathrm{v})|\mathrm{E}|^2(\mathrm{v})d\mathrm{v}.$
	\bigbreak
	\noindent
	Therefore, with the help of Lemma \ref{4.3}, due to duality pairing between the spaces $H^{1,\frac{1}{2}}\Big(\Omega\times(0,\mathrm{t})\Big)$ and $H^{-1,-\frac{1}{2}}\Big(\Omega\times(0,\mathrm{t})\Big)$, and by Sobolev embedding $\mathrm{L}^2\Big(B\times\mathbb{R}_+\Big) \xhookrightarrow {\mathrm{i}} H^{-1,-\frac{1}{2}}\Big(B\times\mathbb{R}_+\Big) $ the following estimate holds
	\begin{align}
		\nonumber
		\textbf{err}^{(3)} &:= \frac{\omega \cdot \boldsymbol{\Im}(\varepsilon_\mathrm{p})}{2\pi\gamma_{\mathrm{p}}}\int_{0}^{t}\int_{\Omega}\Phi^{\textbf{e}}(\xi,t;z,\tau)\int_{\Omega}\Phi(\mathrm{y},\tau;\mathrm{v})|\mathrm{E}|^{2}(\mathrm{v}) d\mathrm{v}d\mathrm{y}d\tau  \nonumber \\ \nonumber&\lesssim \frac{\omega \cdot \boldsymbol{\Im}(\varepsilon_\mathrm{p})}{2\pi\gamma_{\mathrm{p}}}\Big\Vert \Phi^{\textbf{e}}(\xi,t;z,\cdot)\Big\Vert_{H^{-1,-\frac{1}{2}}\Big(\Omega\times(0,\mathrm{t})\Big)} \Big\Vert \mathbb{I}_{0}\Big[|\mathrm{E}|^2\Big]\Big\Vert_{\mathrm{H}^{1,\frac{1}{2}}\Big(\Omega\times(0,\mathrm{t})\Big)}
		\\ &\lesssim \nonumber
		\frac{\omega \cdot \boldsymbol{\Im}(\varepsilon_\mathrm{p})}{2\pi\gamma_{\mathrm{p}}}\Big\Vert \rchi_{(0,\mathrm{T}_0)}\Phi^{\textbf{e}}(\xi,\mathrm{t};z,\cdot)\Big\Vert_{H^{-1,-\frac{1}{2}}\Big(\Omega\times\mathbb{R}_+\Big)} \Big\Vert \mathbb{I}_{0}\Big[|\mathrm{E}|^2\Big]\Big\Vert_{\mathrm{H}^{1,\frac{1}{2}}\Big(\Omega \times\mathbb{R_+}\Big)}
		\\ &\lesssim \nonumber
		\frac{\omega \cdot \boldsymbol{\Im}(\varepsilon_\mathrm{p})}{2\pi\gamma_{\mathrm{p}}}\alpha^{\frac{1}{2}}_\mathrm{m}\delta^2\Big\Vert \rchi_{(0,\mathrm{T}_{0,\delta})}\hat{\Phi}^{\textbf{e}}(\xi,\Tilde{t};z,\cdot)\Big\Vert_{H^{-1,-\frac{1}{2}}(B\times\mathbb{R}_+)} \alpha^{\frac{1}{2}}_\mathrm{p}\delta\Big\Vert \mathbb{I}_{0}\Big[|\hat{\mathrm{E}}|^2\Big]\Big\Vert_{\mathrm{H}^{1,\frac{1}{2}}(B\times\mathbb{R}_+)}
		\\ &\lesssim \frac{\omega \cdot \boldsymbol{\Im}(\varepsilon_\mathrm{p})}{2\pi\gamma_{\mathrm{p}}}\alpha^{\frac{1}{2}}_\mathrm{m}\delta^2\Big\Vert \rchi_{(0,\mathrm{T}_{0,\delta})}\hat{\Phi}^{\textbf{e}}(\xi,\Tilde{t};z,\cdot)\Big\Vert_{\mathrm{L}^2\Big(B\times\mathbb{R}_+\Big)} \alpha^{\frac{1}{2}}_\mathrm{p}\delta\Big\Vert |\hat{\mathrm{E}}|^2 \Big\Vert_{\mathrm{L}^2\Big (B\Big)}.
	\end{align}
	Observing that $|\hat{\mathrm{E}}|^2$ is independent of the time $t$ in $(0,\mathrm{T}_{0,\delta})$ and $\hat{\Phi}^{\textbf{e}}$ is independent on the space variable, we can rewrite the above expression as
	\begin{align}
		\nonumber
		\textbf{err}^{(3)} &:= \frac{\omega \cdot \boldsymbol{\Im}(\varepsilon_\mathrm{p})}{2\pi\gamma_{\mathrm{p}}}\int_{0}^{t}\int_{\Omega}\Phi^{\textbf{e}}(\xi,t;z,\tau)\int_{\Omega}\Phi(\mathrm{y},\tau;\mathrm{v})|\mathrm{E}|^{2}(\mathrm{v}) d\mathrm{v}d\mathrm{y}d\tau  \\ \nonumber &\lesssim \frac{\omega \cdot \boldsymbol{\Im}(\varepsilon_\mathrm{p})}{2\pi\gamma_{\mathrm{p}}}\alpha^{\frac{1}{2}}_\mathrm{m}\delta^2\Big\Vert \hat{\Phi}^{\textbf{e}}(\xi,\Tilde{t};z,\cdot)\Big\Vert_{\mathrm{L}^2\Big(0,\mathrm{T}_{0,\delta}\Big)} \alpha^{\frac{1}{2}}_\mathrm{p}\delta\Big\Vert |\hat{\mathrm{E}}|^2 \Big\Vert_{\mathrm{L}^2\Big(B\Big)}.
	\end{align}
	Consequently, when we scale back to $\Omega$, the following expression can be derived
	\begin{align}\label{initialesti}
		\nonumber
		\textbf{err}^{(3)} &:= \frac{\omega \cdot \boldsymbol{\Im}(\varepsilon_\mathrm{p})}{2\pi\gamma_{\mathrm{p}}}\int_{0}^{t}\int_{\Omega}\Phi^{\textbf{e}}(\xi,t;z,\tau)\int_{\Omega}\Phi(\mathrm{y},\tau;\mathrm{v})|\mathrm{E}|^{2}(\mathrm{v}) d\mathrm{v}d\mathrm{y}d\tau \\ \nonumber
		&\lesssim \nonumber \frac{\omega \cdot \boldsymbol{\Im}(\varepsilon_\mathrm{p})}{2\pi\gamma_{\mathrm{p}}}\delta\Big\Vert \Phi^{\textbf{e}}(\xi,t;z,\cdot)\Big\Vert_{\mathrm{L}^2\Big(0,T_0\Big)} \alpha^{\frac{1}{2}}_\mathrm{p}\Big\Vert |\mathrm{E}|^2 \Big\Vert_{\mathrm{L}^2\Big(\Omega\Big)}
		\\ &= \mathcal{O}\Bigg(\frac{\omega \cdot \boldsymbol{\Im}(\varepsilon_\mathrm{p})}{2\pi\gamma_{\mathrm{p}}}\alpha_\mathrm{p}^{\frac{1}{2}}\delta\Big\Vert |\mathrm{E}|^2 \Big\Vert_{\mathrm{L}^2\Big(\Omega\Big)}\sqrt{\int_{0}^{T_{0}}\Big|\Phi^{\textbf{e}}(\xi,t;z,\tau)\Big|^2}d\tau \Bigg).
	\end{align}
	Thus, we deduce from (\ref{volumeesti1}) and (\ref{initialesti}) that
	\begin{align}\label{E1}
		\nonumber
		\textbf{First Term} &:= \frac{1}{\alpha_\mathrm{p}}\frac{\omega \cdot \boldsymbol{\Im}(\varepsilon_\mathrm{p})}{2\pi\gamma_{\mathrm{p}}}\int_{0}^{t}\int_{\partial\Omega}\Phi^{\textbf{e}}(\xi,t;\mathrm{z},\tau)\gamma^{\textbf{int}}_{1,\mathrm{y}}\int_{0}^{\tau}\int_{\Omega}\Phi(\mathrm{y},\tau;\mathrm{v},\mathrm{s})|\mathrm{E}|^{2}(\mathrm{v}) d\mathrm{v}d\mathrm{s} d\sigma_{\mathrm{y}}d\tau \\  &= -\frac{\omega \cdot \boldsymbol{\Im}(\varepsilon_\mathrm{p})}{2\pi\gamma_{\mathrm{p}}}\int_{0}^{t}\int_{\Omega} \Phi^{\textbf{e}}(\xi,t;\mathrm{z},\tau) |\mathrm{E}|^{2}(\mathrm{y}) d\mathrm{y}d\tau + \textbf{err}^{(3)}.
	\end{align}
	Next, we proceed to estimate the second term of (\ref{sinires1}) defined as:
	\begin{align}
		\textbf{Second Term} := \frac{1}{2\pi}\int_{0}^{t}\int_{\partial\Omega}\gamma^{\textbf{int}}_{1}\mathrm{U}_{\mathrm{i}}(\mathrm{y},\tau)\int_{\partial\Omega} \frac{(\mathrm{y}-\mathrm{v})\cdot\nu_{\mathrm{v}}}{2\pi|\mathrm{y}-\mathrm{v}|^2}\Big[\varphi(\mathrm{v},\mathrm{y},\mathrm{t},\tau)-\Phi^{\textbf{e}}(\xi,t;z,\tau)\Big]d\sigma_\mathrm{v} d\sigma_\mathrm{y} d\tau.
	\end{align}
	To do that, we define
	\begin{align}
		\mathbb{E}(\mathrm{y},t) := \int_{\partial\Omega} \frac{(\mathrm{y}-\mathrm{v})\cdot\nu_{\mathrm{v}}}{2\pi|\mathrm{y}-\mathrm{v}|^2}\Big[\varphi(\mathrm{v},\mathrm{y},\mathrm{t},\tau)-\Phi^{\textbf{e}}(\xi,t;z,\tau)\Big]d\sigma_\mathrm{v}.
	\end{align}
	Then, based on the duality pairing between the spaces $\mathrm{H}^{\frac{1}{2},\frac{1}{4}}$ and $\mathrm{H}^{-\frac{1}{2},-\frac{1}{4}}$, we obtain
	\begin{align}\label{secondterm}
		\nonumber
		\textbf{Second Term} &:=  \Bigg|\int_{0}^{t}\int_{\partial\Omega}\gamma^{\textbf{int}}_{1}\mathrm{U}_{\mathrm{i}}(\mathrm{y},\tau)\int_{\partial\Omega} \frac{(\mathrm{y}-\mathrm{v})\cdot\nu_{\mathrm{v}}}{2\pi|\mathrm{y}-\mathrm{v}|^2}\Big[\varphi(\mathrm{v},\mathrm{y},\mathrm{t},\tau)-\Phi^{\textbf{e}}(\xi,t;z,\tau)\Big]d\sigma_\mathrm{v} d\sigma_\mathrm{y} d\tau \Bigg| \\ \nonumber &\lesssim \Big\Vert \gamma^{\textbf{int}}_{1}\mathrm{U}_{\mathrm{i}}\Big\Vert_{\mathrm{H}^{-\frac{1}{2},-\frac{1}{4}}\Big(\partial \Omega\times (0,\mathrm{t})\Big)}\Big\Vert \mathbb{E}(\mathrm{y},\mathrm{t})\Big\Vert_{\mathrm{H}^{\frac{1}{2},\frac{1}{4}}\Big(\partial \Omega\times (0,\mathrm{t})\Big)} \\ &\lesssim \Big\Vert \gamma^{\textbf{int}}_{1}\mathrm{U}_{\mathrm{i}}\Big\Vert_{\mathrm{H}^{-\frac{1}{2},-\frac{1}{4}}\Big(\partial \Omega\times \mathbb{R}\Big)}\Big\Vert \mathbb{E}(\mathrm{y},\mathrm{t})\Big\Vert_{\mathrm{H}^{\frac{1}{2},\frac{1}{4}}\Big(\partial \Omega\times (0,\mathrm{t})\Big)}.
	\end{align}
	We state the following lemma. As its justification is rather lengthy but elementary, we postpone it to the Appendix. 
	\begin{lemma}\label{4.1}
		We set $\varphi(\mathrm{v},\mathrm{y},\mathrm{t},\tau)$ as follows:
		\begin{equation}\label{defvarphi}
			\varphi(\mathrm{v},\mathrm{y},\mathrm{t}, \tau) := \displaystyle\int_{0}^{\tau}\dfrac{\alpha |\mathrm{y}-\mathrm{v}| ^2}{4(\mathrm{s}-\tau)^2}\textbf{exp}\Big(-\dfrac{\alpha|\mathrm{y}-\mathrm{v}| ^2}{4(\mathrm{s}-\tau)}\Big)\Phi^{\textbf{e}}(\xi,\mathrm{t};\mathrm{z},\mathrm{s})d\mathrm{s}.
		\end{equation}
		Then we have
		\begin{equation}
			\varphi(\mathrm{v},\mathrm{y},\mathrm{t}, \tau)-\Phi^{\textbf{e}}(\xi,\mathrm{t};\mathrm{z},\tau) = \mathcal{O}\Bigg(\sqrt{\alpha}|\mathrm{y}-\mathrm{v}| \ \Vert \partial_\mathrm{s}\Phi^{\textbf{e}}(\xi,\mathrm{t};\mathrm{z},\cdot)\Vert_{\mathrm{H}^{-\frac{1}{4}}(0,t)}\Bigg),
		\end{equation}
		for $\mathrm{x}, \mathrm{y}$ such that $|\mathrm{v}-\mathrm{y}| <<1$ and $\mathrm{t} \in (0,\mathrm{T}]$ uniformly with respect to $\Omega$. 
	\end{lemma}
	\noindent
	Recall that if $\Omega$ is a bounded open subset of $\mathbb{R}^n$ of class $\mathcal{C}^{1,\alpha}$, with $\alpha\in (0,1]$, then there exist a constant $c_{\Omega,\alpha}>0$ such that
	$
	|(\mathrm{y}-\mathrm{v})\cdot\nu_{\mathrm{v}}| \le c_{\Omega,\alpha} |\mathrm{y} -\mathrm{v}|^{1+\alpha}.$
	Using Lemma \ref{4.1}, the duality pairing on spatial domain, and the Sobolev embedding $ \mathrm{L}^2\Big(\partial \mathrm{B}\times\mathbb{R}\Big) \xhookrightarrow{\mathrm{i}} \mathrm{H}^{-\frac{1}{2},-\frac{1}{4}}\Big(\partial \mathrm{B}\times \mathbb{R}\Big)$, we obtain
	\begin{align}\label{epp1}
		\nonumber
		\mathbb{V}_1 := \Big\Vert \mathbb{E}(\mathrm{y},\mathrm{t})\Big\Vert_{\mathrm{L}^2\Big(\partial \Omega\times (0,\mathrm{t})\Big)}^2 &= \displaystyle \int_{0}^{t}\int_{\partial\Omega}\Bigg|\int_{\partial\Omega} \frac{(\mathrm{y}-\mathrm{v})\cdot\nu_{\mathrm{v}}}{2\pi|\mathrm{y}-\mathrm{v}|^2}\Big[\varphi(\mathrm{v},\mathrm{y},\mathrm{t},\tau)-\Phi^{\textbf{e}}(\xi,t;z,\tau)\Big]d\sigma_\mathrm{v}\Bigg|^2 d\sigma_\mathrm{y}d\mathrm{t}
		\\ \nonumber &= \mathcal{O}\Bigg(\int_{0}^{t}\int_{\partial\Omega}\int_{\partial\Omega} \frac{|(\mathrm{y}-\mathrm{v})\cdot\nu_{\mathrm{v}}|^2}{|\mathrm{y}-\mathrm{v}|^2} \alpha_\mathrm{p}|\mathrm{y}-\mathrm{v}|^2\Big\Vert \partial_\mathrm{s}\Phi^{\textbf{e}}(\xi,t;z,\cdot)\Big\Vert_{\mathrm{H}^{-\frac{1}{4}}(0,t)} d\sigma_\mathrm{v} d\sigma_\mathrm{y}dt\Bigg)
		\\ \nonumber &=  \mathcal{O}\Bigg(\alpha_\mathrm{p}\delta^{\frac{9}{2}}\int_{0}^{t}\int_{\partial\Omega}\Big\Vert \partial_\mathrm{s}\Phi^{\textbf{e}}(\xi,t;z,\cdot)\Big\Vert_{\mathrm{H}^{-\frac{1}{2}, -\frac{1}{4}}\Big(\partial \Omega\times(0,\mathrm{t})\Big)}d\sigma_\mathrm{y}dt\Bigg)
		\\ \nonumber &=  \mathcal{O}\Bigg(\alpha_\mathrm{p}\delta^{\frac{9}{2}}\int_{0}^{t}\int_{\partial\Omega}\Big\Vert \rchi_{(0,\mathrm{T_0})} \partial_\mathrm{s}\Phi^{\textbf{e}}(\xi,t;z,\cdot)\Big\Vert_{\mathrm{H}^{-\frac{1}{2}, -\frac{1}{4}}\Big(\partial \Omega\times \mathbb{R}\Big)}d\sigma_\mathrm{y}dt\Bigg)
		\\ \nonumber &= \mathcal{O}\Bigg(\alpha_\mathrm{p}\delta^{\frac{9}{2}}\int_{0}^{t}\int_{\partial\Omega}\alpha^{-\frac{1}{2}}_{\mathrm{m}}\delta^{-\frac{1}{2}}\Big\Vert \rchi_{(0,\mathrm{T_{0,\delta}})}\partial_{\Tilde{\mathrm{s}}}\hat{\Phi}^{\textbf{e}}(\xi,t;z,\cdot)\Big\Vert_{\mathrm{H}^{-\frac{1}{2},-\frac{1}{4}}\Big(\partial B\times \mathbb{R}\Big)}d\sigma_\mathrm{y}dt\Bigg)
		\\ \nonumber &= \mathcal{O}\Bigg(\alpha_\mathrm{p}\delta^{\frac{9}{2}}\int_{0}^{t}\int_{\partial\Omega}\alpha^{-\frac{1}{2}}_{\mathrm{m}}\delta^{-\frac{1}{2}}\Big\Vert \rchi_{(0,\mathrm{T}_{0,\delta})}\partial_{\Tilde{\mathrm{s}}}\hat{\Phi}^{\textbf{e}}(\xi,t;z,\cdot)\Big\Vert_{\mathrm{L}^2\Big(\partial B\times \mathbb{R}\Big)}d\sigma_\mathrm{y}dt\Bigg)
		\\ \nonumber &= \mathcal{O}\Bigg(\alpha_\mathrm{p}\delta^{\frac{9}{2}}\int_{0}^{t}\int_{\partial\Omega}\Big\Vert \partial_\mathrm{s}\Phi^{\textbf{e}}(\xi,t;z,\cdot)\Big\Vert_{\mathrm{L}^2\Big(\partial \Omega\times (0,\mathrm{T}_{0})\Big)}d\sigma_\mathrm{y}dt\Bigg)
		\\&= \nonumber \mathcal{O}\Bigg(\alpha_\mathrm{p}\delta^{\frac{11}{2}} \int_{0}^{t}\Big\Vert \partial_\mathrm{s}\Phi^{\textbf{e}}(\xi,t;z,\cdot)\Big\Vert_{\mathrm{L}^2\Big(\partial \Omega\times (0,\mathrm{T}_0)\Big)}dt\Bigg) \\
		&= \mathcal{O}\Bigg(\alpha_\mathrm{p}\delta^6 \int_{0}^{t}\Big\Vert \partial_\mathrm{s}\Phi^{\textbf{e}}(\xi,t;z,\cdot)\Big\Vert_{\mathrm{L}^2(0,T_0)}dt\Bigg).
	\end{align}
	The above expression needs to be estimated with respect to the $H^{1,\frac{1}{2}}\Big(\partial \Omega\times \mathbb{R}\Big)$-norm.
	As $\mathcal{K}_{\textbf{Lap}}\Big[\Phi^{\textbf{e}}(\xi,t;z,\tau)\Big](y) = -\frac{1}{2}\Phi^{\textbf{e}}(\xi,t;z,\tau)$, then
	\begin{align}\label{klap}
		\Big\Vert \mathcal{K}_{\text{Lap}}\Big[\Phi^{\textbf{e}}(\xi,t;z,\tau)\Big]\Big\Vert_{\mathrm{H}^{1,\frac{1}{2}}\Big(\partial \Omega\times \mathbb{R}\Big)}^2 = \frac{1}{4} \Big\Vert\Phi^{\textbf{e}}(\xi,t;z,\tau)\Big\Vert_{\mathrm{H}^{1,\frac{1}{2}}\Big(\partial \Omega\times \mathbb{R}\Big)}^2.
	\end{align}
	Also, from the continuity of the operator $\mathcal{K}: H^{1,\frac{1}{2}}\Big(\partial \Omega\times \mathbb{R}\Big) \rightarrow H^{1,\frac{1}{2}}\Big(\partial \Omega\times \mathbb{R}\Big)$, see \cite{hofman}, we obtain
	\begin{align}\label{kheat}
		\Big\Vert \mathcal{K}\Big[\Phi^{\textbf{e}}(\xi,t;z,\tau)\Big]\Big\Vert_{\mathrm{H}^{1,\frac{1}{2}}\Big(\partial \Omega\times \mathbb{R}\Big)}^2 \lesssim \Big\Vert \Phi^{\textbf{e}}(\xi,t;z,\tau)\Big\Vert_{\mathrm{H}^{1,\frac{1}{2}}\Big(\partial \Omega\times \mathbb{R}\Big)}^2.
	\end{align}
	We set $\mathbb{E} := \Big(\mathcal{K}_{\alpha_\mathrm{p}}- \mathcal{K}_{\text{Lap}}\Big)\Big[\Phi^{\textbf{e}}(\xi,t;z,\tau)\Big].$
	Hence, with the help of the previous two expression (\ref{klap}) and (\ref{kheat}), considering the independence of $\Phi^{\textbf{e}}(\xi,t;z,\tau)$ on space variable, and observing the fact that we are only concerned with the time interval $(0,\mathrm{T_0})$, we derive that
	\begin{align}\label{ep2}
		\nonumber
		\mathbb{V}_2 &:= \Big\Vert \mathbb{E}(\mathrm{y},t)\Big\Vert_{\mathrm{H}^{1,\frac{1}{2}}\Big(\partial \Omega\times (0,\mathrm{t)}\Big)}^2 \\ \nonumber &\lesssim \Big\Vert \mathbb{E}(\mathrm{y},t)\Big\Vert_{\mathrm{H}^{1,\frac{1}{2}}\Big(\partial \Omega\times \mathbb{R}\Big)}^2 \lesssim\Big\Vert \Phi^{\textbf{e}}(\xi,t;z,\cdot)\Big\Vert_{\mathrm{H}^{1,\frac{1}{2}}\Big(\partial \Omega\times \mathbb{R}\Big)}^2
		\\ &= \mathcal{O}\Bigg(\delta \int_{0}^{\mathrm{T_0}} \Big[|\nabla_{\textbf{tan}}\Phi^{\textbf{e}}(\xi,t;z,\tau)|^2 + \Big|\Phi^{\textbf{e}}(\xi,t;z,\tau)\Big|^2 + |\partial_{t}^{\frac{1}{2}}\Phi^{\textbf{e}}(\xi,t;z,\tau)|^2\Big]d\tau\Bigg).
	\end{align}
	Thus, by interpolation, we obtain, based on (\ref{epp1}) and (\ref{ep2})
	\begin{align}\label{s2}
		\nonumber
		\Big\Vert \mathbb{E}(\mathrm{y},t)\Big\Vert_{\mathrm{H}^{\frac{1}{2},\frac{1}{4}}\Big(\partial \Omega\times \mathbb{R}\Big)} &\lesssim \Big\Vert \mathbb{E}(\mathrm{y},t)\Big\Vert_{\mathrm{H}^{1,\frac{1}{2}}\Big(\partial \Omega\times (0,\mathrm{t)}\Big)}^{\frac{1}{2}}\cdot \Big\Vert \mathbb{E}(\mathrm{y},t)\Big\Vert_{\mathrm{L}^2\Big(\partial \Omega\times (0,\mathrm{t)}\Big)}^{\frac{1}{2}}
		\\ &= \mathcal{O}\Bigg(\alpha_\mathrm{p}^{\frac{1}{4}} \delta^{\frac{7}{4}} \sqrt[4]{\mathbb{V}_1 \cdot \mathbb{V}_2}\Bigg).
	\end{align}
	Consequently, we provide the following estimate with the help of previous estimate (\ref{s2}) and inserting it into (\ref{secondterm})
	\begin{align*}
		\nonumber
		\textbf{Second Term} &:= \Bigg|\int_{0}^{t}\int_{\partial\Omega}\gamma^{\textbf{int}}_{1}\mathrm{U}_{\mathrm{i}}(\mathrm{y},\tau)\int_{\partial\Omega} \frac{(\mathrm{y}-\mathrm{v})\cdot\nu_{\mathrm{v}}}{2\pi|\mathrm{y}-\mathrm{v}|^2}\Big[\varphi(\mathrm{v},\mathrm{y},\mathrm{t},\tau)-\Phi^{\textbf{e}}(\xi,t;z,\tau)\Big]d\sigma_\mathrm{v} d\sigma_\mathrm{y} d\tau \Bigg| 
		\\ &= \mathcal{O}\Bigg(\alpha_\mathrm{p}^{\frac{1}{4}} \delta^{\frac{7}{4}} \Big\Vert \gamma^{\textbf{int}}_{1}\mathrm{U}_{\mathrm{i}}\Big\Vert_{\mathrm{H}^{-\frac{1}{2},-\frac{1}{4}}\Big(\partial \Omega\times \mathbb{R}\Big)} \sqrt[4]{\mathbb{V}_1 \cdot \mathbb{V}_2}\Bigg). 
	\end{align*}
	Additionally, we redefine the term above with the following notation
	\begin{align}\label{secondtermfinal}
		\textbf{err}^{(4)} := \mathcal{O}\Bigg(\alpha_\mathrm{p}^{\frac{1}{4}} \delta^{\frac{7}{4}} \Big\Vert \gamma^{\textbf{int}}_{1}\mathrm{U}_{\mathrm{i}}\Big\Vert_{\mathrm{H}^{-\frac{1}{2},-\frac{1}{4}}\Big(\partial \Omega\times \mathbb{R}\Big)} \sqrt[4]{\mathbb{V}_1 \cdot \mathbb{V}_2}\Bigg).
	\end{align}
	As for the third term of (\ref{sinires1}), in the same way as elliptic equations, the bi-linear form of the hyper-singular heat operator can be represented as a weakly singular operator. Moreover, we state the following result, see \cite[Theorem 6.1]{costabel} for details.
	\begin{lemma}\label{5.5}
		Let $\partial_{T}$ and $\partial_{\mathrm{t}}$ denote the derivative with respect to the arc length on $\Gamma$ and the time derivative, respectively. let $\underline{\mathrm{n}}$ denote the exterior normal vector. Then
		\begin{align}
			\Big\langle \Phi^{\textbf{e}}(\xi,\mathrm{t};z,\mathrm{\tau)},\mathbb{H}\Big[\gamma^{\textbf{int}}_{0}\mathrm{U}_{\mathrm{i}}\Big]\Big\rangle_{(\partial\Omega)_{\mathrm{T}}} \nonumber&= \Big\langle \partial_{\mathrm{T}}\Phi^{\textbf{e}}(\xi,\mathrm{t};z,\mathrm{\tau)},\mathcal{S}\Big[\partial_{\mathrm{T}}\gamma^{\textbf{int}}_{0}\mathrm{U}_{\mathrm{i}}\Big]\Big\rangle_{(\partial\Omega)_{\mathrm{T}}} - \Big\langle \Phi^{\textbf{e}}(\xi,\mathrm{t};z,\mathrm{\tau)}\underline{\mathrm{n}},\partial_{t}\mathcal{S}\Big[\underline{\mathrm{n}}\gamma^{\textbf{int}}_{0}\mathrm{U}_{\mathrm{i}}\Big]\Big\rangle_{(\partial\Omega)_{\mathrm{T}}},
		\end{align}
		where $ \Phi^{\textbf{e}}(\xi,\mathrm{t};z,\mathrm{\tau)}$ and $\gamma^{\textbf{int}}_{0}\mathrm{U}_{\mathrm{i}} \in \mathrm{H}^{\frac{1}{2},\frac{1}{4}}\Big(\partial \Omega \times \mathbb{R}\Big).$
	\end{lemma}
	\noindent
	It was noted in \cite{costabel} that we should use the following identity together with the distributional time derivative on $\mathbb{R}\times \Gamma$, $
	\Big\langle \Phi^{\textbf{e}}(\xi,\mathrm{t};z,\mathrm{\tau)}\underline{\mathrm{n}},\partial_{t}\mathcal{S}\Big[\underline{\mathrm{n}}\gamma^{\textbf{int}}_{0}\mathrm{U}_{\mathrm{i}}\Big]\Big\rangle_{(\partial\Omega)_{\mathrm{T}}} = \Big\langle \partial_{t}\Phi^{\textbf{e}}(\xi,\mathrm{t};z,\mathrm{\tau)}\underline{\mathrm{n}},\mathcal{S}\Big[\underline{\mathrm{n}}\gamma^{\textbf{int}}_{0}\mathrm{U}_{\mathrm{i}}\Big]\Big\rangle_{(\partial\Omega)_{\mathrm{T}}}$. There is an alternative representation of the above formula given in \cite[Theorem 4.1]{graz1}. \\
	\bigbreak
	\noindent
	If we set $\psi \displaystyle\equiv \Phi(\xi,t;z,\tau)$, then from Lemma \ref{5.5} we deduce
	\begin{align}\label{bilinear}
		\nonumber
		\textbf{Third Term} &:=\Big\langle \Phi^{\textbf{e}}(\xi,\mathrm{t};z,\mathrm{\tau)},\mathbb{H}_{\alpha_\mathrm{p}}\Big[\gamma^{\textbf{int}}_{0}\mathrm{U}_{\mathrm{i}}\Big]\Big\rangle_{(\partial\Omega)_{\mathrm{T}}} \\ &= -\Big\langle \mathcal{S}_{\alpha_\mathrm{p}}[\gamma^{\textbf{int}}_{0}\mathrm{T}_{\mathrm{i}}\underline{\mathrm{n}}], \frac{\partial}{\partial \mathrm{t}}\Phi^{\textbf{e}}(\xi,t;z,\tau)
		\underline{\mathrm{n}}\Big\rangle_{(\partial\Omega)_{\mathrm{T}}}.
	\end{align}
	\noindent
	In the next step, we use Lemma \ref{3.6}, Lemma \ref{4.7}, duality of the function space $\mathrm{L}^2\Big(\Omega\Big)$, Sobolev embeddings $\mathrm{L}^2\Big(\partial B\times\mathbb{R}\Big) \xhookrightarrow {\mathrm{i}} H^{-1,-\frac{1}{2}}\Big(\partial B\times\mathbb{R}\Big)$, $\mathrm{H}^{\frac{1}{2},\frac{1}{4}}\Big(\partial B\times\mathbb{R}\Big) \xhookrightarrow{\mathrm{i}} \mathrm{L}^2\Big(\partial B\times\mathbb{R}\Big)$ and continuity of the single layer operator $\overline{\mathcal{S}}: H^{-1,-\frac{1}{2}} {\Big(\partial B\times\mathbb{R}\Big)} \to \mathrm{L}^2 {\Big(\partial B\times\mathbb{R}\Big)}$,  to obtain the following estimate:
	\begin{align}
		\nonumber
		\textbf{Third Term} &:= \Big\langle \mathcal{S}_{\alpha_\mathrm{p}}[\gamma^{\textbf{int}}_{0}\mathrm{U}_{\mathrm{i}}\underline{\mathrm{n}}], \frac{\partial}{\partial \mathrm{t}}\Phi^{\textbf{e}}(\xi,t;z,\tau)
		\underline{\mathrm{n}}\Big\rangle_{\partial\Omega\times(0,\mathrm{t})}
		\\ \nonumber &\lesssim 
		\Big\Vert \frac{\partial}{\partial \mathrm{t}}\Phi^{\textbf{e}}(\xi,t;z,\tau)
		\underline{\mathrm{n}}\Big\Vert_{\mathrm{L}^2\Big(\partial\Omega\times(0,\mathrm{t})\Big)} \Big\Vert \mathcal{S}_{\alpha_\mathrm{p}}\Big[\gamma^{\textbf{int}}_{0}\mathrm{U}_{\mathrm{i}}\underline{\mathrm{n}}\Big]\Big\Vert_{\mathrm{L}^2\Big(\partial\Omega\times(0,\mathrm{t})\Big)}
		\nonumber 
		\\ &\lesssim  \nonumber
		\Big\Vert\rchi_{(0,\mathrm{T}_0)} \frac{\partial}{\partial \mathrm{t}}\Phi^{\textbf{e}}(\xi,t;z,\tau)
		\underline{\mathrm{n}}\Big\Vert_{\mathrm{L}^2\Big(\partial \Omega\times\mathbb{R}\Big)}\Big\Vert \mathcal{S}_{\alpha_\mathrm{p}}\Big[\underline{\mathrm{n}}\gamma^{\textbf{int}}_{0}\mathrm{U}_{\mathrm{i}}\Big]\Big\Vert_{\mathrm{L}^2\Big(\partial \Omega\times \mathbb{R}\Big)}
		\nonumber 
		\\ &\lesssim  \nonumber
		\alpha_{\mathrm{m}}^{-\frac{1}{2}}\delta^{-\frac{1}{2}}\Big\Vert\rchi_{(0,\mathrm{T}_{0,\delta})}  \frac{\partial}{\partial \Tilde{t}}\Hat{\Phi}^{\textbf{e}}(\xi,\Tilde{t};z,\Tilde{\tau})
		\underline{\mathrm{n}}\Big\Vert_{\mathrm{L}^2\Big(\partial B\times \mathbb{R}\Big)} \alpha_{\mathrm{p}}^{\frac{1}{2}}\delta^{\frac{5}{2}}\Big\Vert \overline{\mathcal{S}}_{\alpha_\mathrm{p}}\Big[\underline{\mathrm{n}}\gamma^{\textbf{int}}_{0}\hat{\mathrm{T}}_{\mathrm{i}}\Big]\Big\Vert_{\mathrm{L}^2\Big(\partial B\times \mathbb{R}\Big)}
		\\ &\lesssim \nonumber
		\alpha_{\mathrm{m}}^{-\frac{1}{2}}\delta^{-\frac{1}{2}}\Big\Vert \rchi_{(0,\mathrm{T}_{0,\delta})}\frac{\partial}{\partial \Tilde{t}}\Hat{\Phi}^{\textbf{e}}(\xi,\Tilde{t};z,\Tilde{\tau})
		\underline{\mathrm{n}}\Big\Vert_{\mathrm{L}^2\Big(\partial B\times\mathbb{R}\Big)} \alpha_{\mathrm{p}}^{\frac{1}{2}}\delta^{\frac{5}{2}}\Big\Vert \underline{\mathrm{n}}\gamma^{\textbf{int}}_{0}\hat{\mathrm{U}}_{\mathrm{i}}\Big\Vert_{H^{-1,-\frac{1}{2}}\Big(\partial B\times\mathbb{R}\Big)}
		\\ &\lesssim \nonumber
		\alpha_{\mathrm{m}}^{-\frac{1}{2}}\delta^{-\frac{1}{2}}\Big\Vert \rchi_{(0,\mathrm{T}_{0,\delta})}\frac{\partial}{\partial \Tilde{t}}\Hat{\Phi}^{\textbf{e}}(\xi,\Tilde{t};z,\Tilde{\tau})
		\underline{\mathrm{n}}\Big\Vert_{\mathrm{L}^2\Big(\partial B\times\mathbb{R}\Big)} \alpha_{\mathrm{p}}^{\frac{1}{2}}\delta^{\frac{5}{2}}\Big\Vert \underline{\mathrm{n}}\gamma^{\textbf{int}}_{0}\hat{\mathrm{U}}_{\mathrm{i}}\Big\Vert_{\mathrm{L}^2\Big(\partial B\times\mathbb{R}\Big)}
		\\ &\lesssim \nonumber
		\alpha_{\mathrm{m}}^{-\frac{1}{2}}\delta^{-\frac{1}{2}}\Big\Vert \rchi_{(0,\mathrm{T}_{0,\delta})}\frac{\partial}{\partial \Tilde{t}}\Hat{\Phi}^{\textbf{e}}(\xi,\Tilde{t};z,\Tilde{\tau})
		\underline{\mathrm{n}}\Big\Vert_{\mathrm{L}^2\Big(\partial B\times\mathbb{R}\Big)} \alpha_{\mathrm{p}}^{\frac{1}{2}}\delta^{\frac{5}{2}}\Big\Vert \underline{\mathrm{n}}\gamma^{\textbf{int}}_{0}\hat{\mathrm{U}}_{\mathrm{i}}\Big\Vert_{\mathrm{H}^{\frac{1}{2},\frac{1}{4}}\Big(\partial B\times\mathbb{R}\Big)}
		\\ &\lesssim \nonumber
		\Big\Vert \frac{\partial}{\partial t}\Phi^{\textbf{e}}(\xi,t;z,\tau)
		\underline{\mathrm{n}}\Big\Vert_{\mathrm{L}^2\Big(\partial \Omega\times(0,\mathrm{T}_0)\Big)} \delta\Big\Vert \underline{\mathrm{n}}\gamma^{\textbf{int}}_{0}\mathrm{U}_{\mathrm{i}}\Big\Vert_{\mathrm{H}^{\frac{1}{2},\frac{1}{4}}\Big(\partial \Omega\times \mathbb{R}\Big)}
		\\ &= \mathcal{O}\Bigg(\delta^{\frac{3}{2}}\Big\Vert \gamma^{\textbf{int}}_{0}\mathrm{U}_{\mathrm{i}}\Big\Vert_{\mathrm{H}^{\frac{1}{2},\frac{1}{4}}\Big(\partial \Omega\times \mathbb{R}\Big)}\sqrt{\int_{0}^{\mathrm{T}_0}\Big|\partial_{t}\Phi^{\textbf{e}}(\xi,t;z,\tau)\Big|^2 d\tau}\Bigg).
	\end{align}
	To further simplify the above expression, we redefine it with the following notation
	\begin{align}\label{E2}
		\textbf{err}^{(5)} := \mathcal{O}\Bigg(\delta^{\frac{3}{2}}\Big\Vert \gamma^{\textbf{int}}_{0}\mathrm{U}_{\mathrm{i}}\Big\Vert_{\mathrm{H}^{\frac{1}{2},\frac{1}{4}}\Big(\partial \Omega\times \mathbb{R}\Big)}\sqrt{\int_{0}^{\mathrm{T}_0}\Big|\partial_{t}\Phi^{\textbf{e}}(\xi,t;z,\tau)\Big|^2 d\tau}\Bigg).
	\end{align}
	Now, considering (\ref{E1}), (\ref{secondtermfinal}) and (\ref{E2}) and inserting these into (\ref{sinires1}), we obtain
	\begin{align}\label{beforefinal}
		\nonumber
		\int_{0}^{t}\int_{\partial\Omega}\Phi^{\textbf{e}}(\xi,t;z,\tau)\gamma^{\textbf{int}}_{1}\mathrm{U}_{\mathrm{i}}(\mathrm{y},\tau)d\sigma_\mathrm{y}d\tau &= -\frac{\omega \cdot \boldsymbol{\Im}(\varepsilon_\mathrm{p})}{2\pi\gamma_{\mathrm{p}}}\int_{0}^{t}\int_{\Omega}\Phi^{\textbf{e}}(\xi,t;z,\tau)|\mathrm{E}|^{2}(\mathrm{y}) d\mathrm{y} d\tau + \textbf{err}^{(3)} + \textbf{err}^{(4)}\\  &+ \textbf{err}^{(5)}.
	\end{align}
	We recall the previously derived expression
	\begin{align}\label{mainformula1}
		\mathrm{U}_{\mathrm{e}}(\xi,t) &=\frac{\gamma_{\mathrm{p}}}{\gamma_{\mathrm{m}}}\frac{1}{\alpha_\mathrm{m}} \Bigg[-\int_{0}^t\int_{\partial\Omega}\Phi^{\textbf{e}}(\xi,t;z,\tau)\gamma^{\textbf{int}}_{1}\mathrm{U}_{\mathrm{i}}(\mathrm{y},\tau)d\sigma_\mathrm{y}d\tau  + \textbf{err}^{(1)} + \textbf{err}^{(2)}\Bigg].
	\end{align}
	Therefore, using Proposition \ref{pr4}, considering the regime $\frac{1}{\gamma_\mathrm{p}} \sim \delta^2$, $\alpha_{\mathrm{m}} \sim 1$, inserting (\ref{beforefinal}) in (\ref{mainformula1}) we deduce the desired expression of the heat potential around the inserted plasmonic nanoparticle for $(\xi,t) \in (\mathbb{R}^2\setminus \overline{\Omega})_T $ and $\mathrm{z} \in \Omega$:
	\begin{align}\label{formula}
		\mathrm{U}_{\mathrm{e}}(\xi,t) &= \frac{\gamma_{\mathrm{p}}}{\gamma_{\mathrm{m}}}\frac{1}{\alpha_\mathrm{m}}\Bigg[\frac{\omega \cdot \boldsymbol{\Im}(\varepsilon_\mathrm{p})}{2\pi\gamma_{\mathrm{p}}}\int_{0}^{t}\int_{\Omega} \Phi^{\textbf{e}}(\xi,t;z,\tau)|\mathrm{E}|^{2}(\mathrm{y}) d\mathrm{y}d\tau + \textbf{err}^{(1)} + \textbf{err}^{(2)} + \textbf{err}^{(3)} + \textbf{err}^{(4)} + \textbf{err}^{(5)}\Bigg],
	\end{align}
	where
	\begin{align}
		\nonumber
		\textbf{err}^{(1)} = \mathcal{O}\Bigg(\delta^{7}\sqrt{\mathcal{K}^{(\mathrm{T_0})}_{\mathrm{r}}}\frac{1}{|\xi-z|^{\frac{7}{2}-2\mathrm{r}}}\Big\Vert|\mathrm{E}|^{2}\Big\Vert_{\mathrm{L}^2(\Omega)}\Bigg),\
		\textbf{err}^{(2)} = \mathcal{O}\Bigg(\delta^{7}\sqrt{\mathcal{K}^{(\mathrm{T_0})}_{\mathrm{r}}}\frac{1}{|\xi-\cdot|^{2-2\mathrm{r}}}\Big\Vert|\mathrm{E}|^{2}\Big\Vert_{\mathrm{L}^2(\Omega)}\Bigg)
	\end{align}
	\begin{equation*}
		\textbf{err}^{(3)} = \mathcal{O}\Bigg(\delta^{4}\Big\Vert |\mathrm{E}|^2 \Big\Vert_{\mathrm{L}^2\Big(\Omega\Big)}\sqrt{\int_{0}^{\mathrm{T_0}}\Big|\Phi^{\textbf{e}}(\xi,t;z,\tau)\Big|^2}d\tau \Bigg),\ 
		\textbf{err}^{(4)} = \mathcal{O}\Bigg(\delta^{\frac{23}{4}}\Big\Vert|\mathrm{E}|^{2}\Big\Vert_{\mathrm{L}^2(\Omega)}\sqrt[4]{\mathbb{V}_1\cdot\mathbb{V}_2}\Bigg)
	\end{equation*}
	and
	\begin{equation*}
		\textbf{err}^{(5)} = \mathcal{O}\Bigg(\delta^{6}\Big\Vert|\mathrm{E}|^{2}\Big\Vert_{\mathrm{L}^2(\Omega)}\sqrt{\int_{0}^{\mathrm{T_0}}\Big|\partial_{t}\Phi^{\textbf{e}}(\xi,t;z,\tau)\Big|^2 d\tau}\Bigg).     
	\end{equation*}
	
	\noindent
	Next, we derive a simplified form of the above error terms. In order to do so, we compute the classical singularity estimate for the fundamental solution $\Phi(\xi,t;z,\tau)$ and we also refer to \cite[Chapter 1]{friedman}, \cite[Chapter~9]{kress}. The key is to use the inequality $\mathrm{s}^\mathrm{r}e^{-\mathrm{s}}\le \mathrm{r}^\mathrm{r}e^{-\mathrm{m}}$ with $0<\mathrm{s},\mathrm{r}<\infty$ and $\mathrm{s}:= \alpha|\mathrm{x}-\mathrm{y}|^2/4(t-\tau)$. Thus we deduce the following singular estimates
	\begin{align}
		\Big|\Phi(\mathrm{x},t;\mathrm{y},\tau)\Big| \lesssim \frac{\alpha^\mathrm{r}}{(t-\tau)^r}\frac{1}{|\mathrm{x}-\mathrm{y}|^{2-2r}}, \ \ \ \mathrm{r}<1,
	\end{align}
	\begin{align}
		\Big|\partial_{\mathrm{x}_{i}}\Phi(\mathrm{x},t;\mathrm{y},\tau)\Big| \lesssim \frac{\alpha^\mathrm{r}}{(t-\tau)^\mathrm{r}}\frac{1}{|\mathrm{x}-\mathrm{y}|^{3-2\mathrm{r}}}, \ \ \ \mathrm{r}<2,\ i=1,2,
	\end{align}
	\begin{align}\label{n2}
		\Big|\partial_\nu\Phi(\mathrm{x},t;\mathrm{y},\tau)\Big| \lesssim \frac{\alpha^\mathrm{r}}{(t-\tau)^r}\frac{1}{|\mathrm{x}-\mathrm{y}|^{2-2r}}, \ \ \ \mathrm{r}<2, \mbox{ for } x, y \in \partial \Omega,
	\end{align}
	\begin{align}
		\Big|\partial_{t}\Phi(\mathrm{x},t;\mathrm{y},\tau)\Big| \lesssim \frac{\alpha^{1-r}}{(t-\tau)^\mathrm{r}}\frac{1}{|\mathrm{x}-\mathrm{y}|^{4-2\mathrm{r}}}, \ \ \ r<2,
	\end{align}
	and
	\begin{align}
		\Big|\partial_{t}^{\frac{1}{2}}\Phi(\mathrm{x},t;\mathrm{y},\tau)\Big| \lesssim \frac{\alpha^r}{(t-\tau)^\mathrm{r}}\frac{1}{|\mathrm{x}-\mathrm{y}|^{3-2\mathrm{r}}}, \ \ \ r<2.
	\end{align}
	From the above estimates we deduce that $\displaystyle
	\int_{0}^{\mathrm{T_0}}|\Phi^{\textbf{e}}(\xi,t;z,\tau)|^2d\tau \lesssim \frac{1}{|\xi-z|^{4-4\mathrm{r}}} \int_{0}^{\mathrm{T_0}}\frac{1}{(t-\tau)^{2\mathrm{r}}}d\tau$.  
	We denote by $\mathcal{K}^{(\mathrm{T_0})}_{\mathrm{r}}:= \sup_{t\in(0,T)}\displaystyle \int_{0}^{\mathrm{T_0}}\frac{1}{(t-\tau)^{2\mathrm{r}}}dt \ \ < +\infty \ \ \text{for} \ \mathrm{r}<\frac{1}{2}$, then we have
	$ \displaystyle
	\sqrt{\int_{0}^{\mathrm{T_0}}|\Phi^{\textbf{e}}(\xi,t;z,\tau)|^2d\tau} = \mathcal{O}\Bigg(\sqrt{\mathcal{K}^{(\mathrm{T_0})}_{\mathrm{r}}}\frac{1}{|\xi-z|^{2-2\mathrm{r}}} \Bigg)$
	and therefore, 
	\begin{align}\label{ff1}
		\textbf{err}^{(3)} := \mathcal{O}\Bigg(\varepsilon^4\sqrt{\mathcal{K}^{(\mathrm{T_0})}_{\mathrm{r}}}\frac{1}{|\xi-z|^{2-2\mathrm{r}}}\Big\Vert|\mathrm{E}|^{2}\Big\Vert_{\mathrm{L}^2(\Omega)} \Bigg).
	\end{align}
	We also note that $\nabla_{\text{tan}}u := \nabla u -(\partial_{\nu}u)\nu$. Analogously, we find
	\begin{align}\label{int1}
		\sqrt[4]{\int_{0}^\mathrm{T_0} \Big[|\nabla_{\textbf{tan}}\Phi^{\textbf{e}}(\xi,t;z,\tau)|^2 + \Phi^{\textbf{e}}(\xi,t;z,\tau)^2 + \Big(\partial_{t}^{\frac{1}{2}}\Phi^{\textbf{e}}(\xi,t;z,\tau)\Big)^2\Big] d\tau} = \mathcal{O}\Bigg(\sqrt[4]{\mathcal{K}^{(\mathrm{T_0})}_{\mathrm{r}}}\frac{1}{|\xi-z|^{\frac{3}{2}-\mathrm{r}}}\Bigg).
	\end{align}
	Based on the previous estimates, we also obtain the estimate
	\begin{align}\label{int2}
		\sqrt[4]{\int_{0}^{\mathrm{T_0}}\int_{0}^{\mathrm{T_0}} |\partial_{t}\Phi^{\textbf{e}}(\xi,t;z,\tau)|^2d\tau dt} = \mathcal{O}\Bigg(\sqrt[4]{\mathcal{S}^{(\mathrm{T_0})}_{\mathrm{r}}}\frac{1}{|\xi-z|^{2-2\mathrm{r}}}\Bigg),
	\end{align}
	where we denote by $
	\mathcal{S}^{(\mathrm{T_0})}_{\mathrm{r}} \equiv \sup_{t\in(0,T)}\displaystyle \int_{0}^{\mathrm{T_0}}\int_{0}^{\mathrm{T_0}}\frac{1}{(t-\tau)^{2\mathrm{r}}}d\tau \ \ < +\infty \ \ \text{for} \ \mathrm{r}<1.$
	Thus we obtain from (\ref{int1}) and (\ref{int2}) $
	\sqrt[4]{\mathbb{V}_1 \cdot \mathbb{V}_2} = \mathcal{O}\Bigg(\sqrt[4]{\mathcal{K}^{(\mathrm{T_0})}_{\mathrm{r}}\mathcal{S}^{(\mathrm{T_0})}_{\mathrm{r}}}\frac{1}{|\xi-z|^{\frac{7}{2}-2\mathrm{r}}}\Bigg),$
	and then
	\begin{align}\label{ff2}
		\textbf{err}^{(4)} = \mathcal{O}\Bigg(\varepsilon^{\frac{23}{4}}\sqrt[4]{\mathcal{K}^{(\mathrm{T_0})}_{\mathrm{r}}\mathcal{S}^{(\mathrm{T_0})}_{\mathrm{r}}}\frac{1}{|\xi-z|^{\frac{7}{2}-2\mathrm{r}}}\Big\Vert|\mathrm{E}|^{2}\Big\Vert_{\mathrm{L}^2(\Omega)}\Bigg).
	\end{align}
	Next, we further derive the estimate $ \displaystyle
	\sqrt{\int_{0}^{\mathrm{T_0}}|\partial_{t}\Phi^{\textbf{e}}(\xi,t;z,\tau)|^2d\tau} = \mathcal{O}\Bigg(\sqrt{\mathcal{K}^{(\mathrm{T_0})}_{\mathrm{r}}}\frac{1}{|\xi-z|^{4-2\mathrm{r}}} \Bigg).$
	Therefore, the following can be obtained
	\begin{align}\label{ff3}
		\textbf{err}^{(5)} = \mathcal{O}\Bigg(\varepsilon^6\sqrt{\mathcal{K}^{(\mathrm{T_0})}_{\mathrm{r}}}\frac{1}{|\xi-z|^{4-2\mathrm{r}}} \Big\Vert|\mathrm{E}|^{2}\Big\Vert_{\mathrm{L}^2(\Omega)}\Bigg).
	\end{align}
	Furthermore, with the similar analysis as above we obtain
	\begin{align}\label{5.53}
		\textbf{err}^{(1)} = \mathcal{O}\Bigg(\delta^{7}\sqrt{\mathcal{K}^{(\mathrm{T_0})}_{\mathrm{r}}}\frac{1}{|\xi-z|^{\frac{7}{2}-2\mathrm{r}}}\Big\Vert|\mathrm{E}|^{2}\Big\Vert_{\mathrm{L}^2(\Omega)}\Bigg)
	\end{align}
	and base on the estimate (\ref{n2})
	\begin{align}\label{5.54}
		\textbf{err}^{(2)} = \mathcal{O}\Bigg(\delta^{7}\sqrt{\mathcal{K}^{(\mathrm{T_0})}_{\mathrm{r}}}\frac{1}{|\xi-\cdot|^{2-2\mathrm{r}}}\Big\Vert|\mathrm{E}|^{2}\Big\Vert_{\mathrm{L}^2(\Omega)}\Bigg).
	\end{align}
	Hence, considering $\xi \in \mathbb{R}^2\setminus\overline{\Omega}$, $\mathrm{y} \in \partial\Omega$, $z \in \Omega$, inserting (\ref{ff1}), (\ref{ff2}), (\ref{ff3}), (\ref{5.53}) and (\ref{5.54}) in (\ref{formula}), we obtain 
	\begin{align}
		\nonumber
		\mathrm{U}_{\mathrm{e}}(\xi,t) &= \frac{\gamma_{\mathrm{p}}}{\gamma_{\mathrm{m}}}\frac{1}{\alpha_\mathrm{m}}\Bigg[\frac{\omega \cdot \boldsymbol{\Im}(\varepsilon_\mathrm{p})}{2\pi\gamma_{\mathrm{p}}}\int_{0}^{t}\int_{\Omega} \Phi^{\textbf{e}}(\xi,t;z,\tau) |\mathrm{E}|^{2}(\mathrm{y}) d\mathrm{y}d\tau + \mathcal{O}\Big(\delta^{4}\sqrt{\mathcal{K}^{(\mathrm{T_0})}_{\mathrm{r}}}\frac{1}{|\xi-z|^{2-2\mathrm{r}}}\Big\Vert |\mathrm{E}|^2 \Big\Vert_{\mathrm{L}^2(\Omega)} \Big) 
		\\ \nonumber &+ \mathcal{O}\Bigg(\delta^{6}\sqrt{\mathcal{K}^{(\mathrm{T_0})}_{\mathrm{r}}}\frac{1}{|\xi-z|^{4-2\mathrm{r}}}\Big\Vert|\mathrm{E}|^{2}\Big\Vert_{\mathrm{L}^2(\Omega)}\Bigg) + \mathcal{O}\Bigg(\delta^{\frac{23}{4}}\sqrt[4]{\mathcal{K}^{(\mathrm{T_0})}_{\mathrm{r}}\mathcal{S}^{(\mathrm{T_0})}_{\mathrm{r}}}\frac{1}{|\xi-z|^{\frac{7}{2}-2\mathrm{r}}}\Big\Vert|\mathrm{E}|^{2}\Big\Vert_{\mathrm{L}^2(\Omega)}\Big) \\
		&+ \mathcal{O}\Bigg(\delta^{7}\sqrt{\mathcal{K}^{(\mathrm{T_0})}_{\mathrm{r}}}\frac{1}{|\xi-z|^{\frac{7}{2}-2\mathrm{r}}}\Big\Vert|\mathrm{E}|^{2}\Big\Vert_{\mathrm{L}^2(\Omega)}\Bigg) 
		+\mathcal{O}\Bigg(\delta^{7}\sqrt{\mathcal{K}^{(\mathrm{T_0})}_{\mathrm{r}}}\frac{1}{|\xi-\cdot|^{2-2\mathrm{r}}}\Big\Vert|\mathrm{E}|^{2}\Big\Vert_{\mathrm{L}^2(\Omega)}\Bigg)\Bigg].
	\end{align}
	\noindent
	The a priori estimates of the electric fields $\mathrm{E}$ are given in the next proposition:
	\begin{proposition}\label{prop2} 
		The electric field corresponding to either the plasmonic or dielectric nanoparticles enjoy the following a priori estimates: 
		\begin{align}\label{el1}
			\Big\Vert|\mathrm{E}|^{2}\Big\Vert_{\mathrm{L}^2\Big(\Omega\Big)} = \mathcal{O}\Big(\delta^{-1}\Big),
		\end{align}
		and 
		\begin{align}\label{el2}
			\Big\Vert |\mathrm{E}|^2 \Big\Vert_{\mathrm{L}^2\Big(\Omega\Big)} = \mathcal{O}\Big(\delta|\log\delta|^{\frac{3\mathrm{h}}{2}}\Big),
		\end{align}
		respectively.
	\end{proposition}
	\noindent
	\textbf{Proof.} See Section \ref{apriori} and Section \ref{dielectric} respectively. 
	\bigbreak
	\noindent
	Let us take $\xi \in \mathbb{R}^3 \setminus{\overline {\Omega}}$ such that $dist(\xi, \Omega)\sim \delta^p$ and then $\vert \xi- z\vert\sim \delta^p+\delta$.\\
	\noindent
	Using (\ref{el1}) in Proposition \ref{prop2}, with the  choice of $\mathrm{p}$ and $\mathrm{r}<\frac{1}{2}$, such that $\frac{1+2\mathrm{p}(1-\mathrm{r})}{2} < \mathrm{h}<1 $ we deduce that
	\begin{align}
		\mathrm{U}_{\mathrm{e}}(\xi,t) &= \frac{\gamma_{\mathrm{p}}}{\gamma_{\mathrm{m}}}\frac{1}{\alpha_\mathrm{m}}\Bigg[\frac{\omega \cdot \boldsymbol{\Im}(\varepsilon_\mathrm{p})}{2\pi\gamma_{\mathrm{p}}}\int_{0}^{t}\int_{\Omega} \Phi^{\textbf{e}}(\xi,t;\mathrm{z},\tau)|\mathrm{E}|^{2}(\mathrm{y}) d\mathrm{y}d\tau + \mathcal{O}\Bigg(\frac{\omega \cdot \boldsymbol{\Im}(\varepsilon_\mathrm{p})}{2\pi}\delta^{3-2\mathrm{p}(1-\mathrm{r})}\sqrt{\mathcal{K}^{(\mathrm{T_0})}_{\mathrm{r}}} \Bigg)\Bigg].
	\end{align}
	This completes the first part of the proof of Theorem \ref{th2}.
	\bigbreak
	\noindent
	Next, under the estimate (\ref{el2}) of Proposition \ref{prop2} and with the condition $2\mathrm{p}(1-\mathrm{r})<1$ we deduce
	\begin{align}
		\mathrm{U}_{\mathrm{e}}(\xi,t) &= \frac{\gamma_{\mathrm{p}}}{\gamma_{\mathrm{m}}}\frac{1}{\alpha_\mathrm{m}}\Bigg[\frac{\omega \cdot \boldsymbol{\Im}(\varepsilon_\mathrm{p})}{2\pi\gamma_{\mathrm{p}}}\int_{0}^{t}\int_{\Omega} \Phi^{\textbf{e}}(\xi,t;\mathrm{z},\tau)|\mathrm{E}|^{2}(\mathrm{y}) d\mathrm{y}d\tau + \mathcal{O}\Bigg(\frac{\omega \cdot \boldsymbol{\Im}(\varepsilon_\mathrm{p})}{2\pi}\delta^{5-2\mathrm{p}(1-\mathrm{r})}|\log\delta|^{\frac{3\mathrm{h}}{2}}\sqrt{\mathcal{K}^{(\mathrm{T_0})}_{\mathrm{r}}}\Bigg)\Bigg].
	\end{align}
	This completes the proof of the second part of the Theorem \ref{th2}.
	
	\section{\Large\textbf{Proof of Theorem \ref{th1}}}\label{sec5}
	\noindent
	In this section, we show the asymptotic analysis of the solution to (\ref{eq:helmholtz}) as $\delta\to 0$ when a plasmonic nanoparticle occupy a bounded domain $\Omega = \mathrm{z} + \delta B.$ We begin to transform the equation (\ref{eq:helmholtz}) into an integral equation. For this, let us recall the fundamental solution of the Helmholtz equation in two dimensions
	\begin{equation*}
		\mathbb{G}^{(\mathrm{k})}(\mathrm{x},\mathrm{y}) = \dfrac{i}{4}\mathcal{H}^{(1)}_{0}(\mathrm{k}|\mathrm{x}-\mathrm{y}|) \quad \text{with} \quad \mathrm{x}\neq \mathrm{y},
	\end{equation*}
	where $\mathcal{H}^{(1)}_{0}$ is the Hankel function of the first kind of order zero.
	The unique solution of the equation (\ref{eq:helmholtz}) leads to the Lippmann–Schwinger equation
	\begin{equation}\label{eq:lm}
		\mathrm{H}(\mathrm{x}) - \alpha \int_\Omega \nabla\mathbb{G}^{(\mathrm{k})}(\mathrm{x},\mathrm{y})\cdot\nabla \mathrm{H}(\mathrm{y}) d\mathrm{y} = \mathrm{H}^{\textbf{in}}(\mathrm{x}),
	\end{equation}
	where $\alpha = \dfrac{1}{\varepsilon_{c}(\omega)} - \dfrac{1}{\varepsilon_{m}}.$ We then take gradient of the equation (\ref{eq:lm}) to get
	\begin{equation}{\label{eq:nabla}}
		\nabla \mathrm{H}(\mathrm{x}) -  \alpha\Big(\nabla \int_\Omega\nabla\mathbb{G}^{(\mathrm{k})}(\mathrm{x},\mathrm{y})\cdot\nabla \mathrm{H}(\mathrm{y}) d\mathrm{y}\Big)= \nabla\mathrm{H}^{\textbf{in}}(\mathrm{x})
	\end{equation}
	Then, we rewrite the above equation as 
	\begin{align}\label{eq:laplace0}
		\nabla \mathrm{H}(\mathrm{x})  -\alpha \Big(\nabla \int_\Omega\nabla\mathbb{G}^{(0)}(\mathrm{x},\mathrm{y})\cdot\nabla \mathrm{H}(\mathrm{y}) d\mathrm{y}\Big) = \nabla\mathrm{H}^{\textbf{in}}(\mathrm{x}) +  \alpha \Big(\nabla\int_\Omega\nabla\Big(\mathbb{G}^{(\mathrm{k})}-\mathbb{G}^{(0)}\Big)(\mathrm{x},\mathrm{y})\cdot\nabla \mathrm{H}(\mathrm{y}) d\mathrm{y}\Big),
	\end{align}
	where $\mathbb{G}^{(0)}(\mathrm{x},\mathrm{y})$ is the fundamental solution of the Laplacian.
	\bigbreak
	\noindent
	We study the system of integral equations (\ref{eq:laplace0}) in the Hilbert space of vector-valued function $\Big(\mathrm{L}^{2}(\Omega)\Big)^2.$ For the sake of simplicity, we use $\mathbb{L}^2(\Omega) =\Big(\mathrm{L}^{2}(\Omega)\Big)^2$. This space can be decomposed into three sub-spaces as a direct sum as following, see \cite{raveski},
	\begin{equation*}
		\mathbb{L}^{2} = \textbf{H}_{0}(\textbf{div},0) \oplus\textbf{H}_{0}(\textbf{curl},0)\oplus \nabla H_{arm},
	\end{equation*}
	where we define these three sub-spaces as follows
	\begin{equation}{\label{eq:subspaces}}
		\begin{cases}
			\textbf{H}_{0}(\textbf{div},0) = \left\{ u\in \mathbb{L}^{2}(\Omega): \nabla \cdot u = 0 \ \text{and} \ u\cdot \nu = 0 \right\},\\
			\textbf{H}_{0}(\textbf{curl},0) = \left\{ u\in \mathbb{L}^{2}(\Omega): \nabla \times u = 0 \ \text{and} \ u\times \nu = 0 \right\}, \\
			\nabla H_{arm} = \left\{ u \in \mathbb{L}^{2}(\Omega): \exists \  \varphi \ \text{s.t.} \ u = \nabla \varphi \ \text{and} \ \Delta \varphi = 0 \right\}. 
		\end{cases}
	\end{equation} 
	Consider $\big(\mathrm{e}^{(1)}_{\mathrm{n}}\big)_{\mathrm{n} \in \mathbb{N}}$ and $\big(\mathrm{e}^{(2)}_{\mathrm{n}}\big)_{\mathrm{n} \in \mathbb{N}}$ to be any orthonormal basis of the sub-spaces $\textbf{H}_{0}(\textbf{div},0)$ and $\textbf{H}_{0}(\textbf{curl},0)$ respectively. But for the sub-space $\nabla H_{arm}$, we consider the complete orthonormal basis $\big(\lambda^{(3)}_{\mathrm{n}},\mathrm{e}^{(3)}_{\mathrm{n}}\big)_{n \in \mathbb{N}}$ derived as the eigenfunctions of the magnetization operator $\mathbb{M}: \nabla H_{arm}\rightarrow \nabla H_{arm},$ \ defined by
	\begin{equation}\label{magnetic operator}
		\mathbb{M}\Big[\nabla \mathrm{H}\Big](\mathrm{x}) = \nabla \int_{\Omega}  \nabla\mathbb{G}^{(0)}(\mathrm{x},\mathrm{y})\cdot\nabla \mathrm{H}(\mathrm{y})d\mathrm{y}.
	\end{equation}
	\noindent
	In addition, we have the expansion
	\begin{equation}\label{expansion21}
		\partial_{\mathrm{ij}}\Big(\mathbb{G}^{(\mathrm{k})}-\mathbb{G}^{(0)}\Big)(\mathrm{x},\mathrm{y}) = 2a \delta_{\mathrm{ij}}\log|x-y| + 2\mathrm{a} \dfrac{(\mathrm{x}-\mathrm{y})_{\mathrm{i}}(\mathrm{x}-\mathrm{y})_{\mathrm{j}}}{|\mathrm{x}-\mathrm{y}|^{2}} + 2\mathrm{b} \delta_{\mathrm{ij}} + \mathcal{O}(|\mathrm{x}-\mathrm{y}|),\quad \mathrm{x}\sim \mathrm{y}.
	\end{equation}
	from which we deduce smoothing properties of the operator $\mathbb{T}$ defined by 
	\begin{equation}\label{T}
		\mathbb{T}\Big[\nabla \mathrm{H}\Big](\mathrm{x}) = \int_{\Omega} \nabla \nabla\cdot \Big(\mathbb{G}^{(\mathrm{k})}-\mathbb{G}^{(0)}\Big)(\mathrm{x},\mathrm{y})\nabla \mathrm{H}(\mathrm{y})d\mathrm{y}.
	\end{equation}
	In particular, we have the following lemma.
	\begin{lemma}\label{5.1}
		The operator $\mathbb{T}$ is a bounded operator from $\mathbb{L}^{2}(\Omega) \rightarrow \mathbb{H}^{2}(\Omega).$
	\end{lemma}
	\noindent
	We will then take projection of the equation (\ref{eq:laplace0}) into each of the above described sub-spaces. 
	
	\subsection{Projection on $\textbf{H}_{0}(\text{div},0)$}
	We begin by considering the orthonormal basis $\big(\mathrm{e}^{(1)}_{\mathrm{n}}\big)_{\mathrm{n} \in \mathbb{N}}$ and take the projection of the equation (\ref{eq:laplace0}) into the subspace $\textbf{H}_{0}(\textbf{div},0)$ to get
	\begin{align}\label{eq:div}
		\Big\langle \nabla \mathrm{H};\mathrm{e}^{(1)}_{\mathrm{n}}\Big\rangle_{\mathbb{L}^{2}(\Omega)} - \alpha \Big\langle \nabla \int_\Omega\nabla\mathbb{G}^{(0)}(\mathrm{x},\mathrm{y})\cdot\nabla \mathrm{H}(\mathrm{y}) d\mathrm{y};\mathrm{e}^{(1)}_{\mathrm{n}}\Big\rangle_{\mathbb{L}^{2}(\Omega)} &= \Big\langle \nabla \mathrm{H}^{\textbf{in}};\mathrm{e}^{(1)}_{\mathrm{n}}\Big\rangle_{\mathbb{L}^{2}(\Omega)}+  \alpha \Big\langle \mathbb{T}\Big[\nabla\mathrm{H}\big];\mathrm{e}^{(1)}_{\mathrm{n}}\Big\rangle_{\mathbb{L}^{2}(\Omega)} \end{align}
	Next, we evaluate the second term of the left hand of the equation (\ref{eq:div}) with the usual induced inner product of $\mathbb{L}^{2}$. Then, from the definition of the sub-space $\textbf{H}_{0}(\textbf{div},0)$ and integration by parts, we obtain
	\begin{align}\label{pr1}
		\Big\langle \nabla \int_\Omega\nabla\mathbb{G}^{(0)}(\mathrm{x},\mathrm{y})\cdot\nabla \mathrm{H}(\mathrm{y}) d\mathrm{y};\mathrm{e}^{(1)}_{\mathrm{n}}\Big\rangle_{\mathbb{L}^{2}\Big(\Omega \Big)} \nonumber&=\int_{\Omega} \mathrm{e}^{(1)}_{\mathrm{n}}\cdot  \nabla\int_\Omega\nabla \mathbb{G}^{(0)}(\mathrm{x},\mathrm{y})\cdot\nabla \mathrm{H}(\mathrm{y}) d\mathrm{y}\\ \nonumber
		&=-\int_{\Omega} \nabla\cdot \mathrm{e}^{(1)}_{\mathrm{n}}\  \int_\Omega\nabla \mathbb{G}^{(0)}(\mathrm{x},\mathrm{y})\cdot\nabla \mathrm{H}(\mathrm{y}) d\mathrm{y} \\ &+ \int_{\partial \Omega} \mathrm{e}^{(1)}_{\mathrm{n}}\cdot\nu \  \int_\Omega\nabla\mathbb{G}^{(0)}(\mathrm{x},\mathrm{y})\cdot\nabla \mathrm{H}(\mathrm{y}) d\mathrm{y} = 0.
	\end{align}
	Using the same arguments, we have
	\begin{align}\label{eq:incident}
		\Big\langle \nabla \mathrm{H}^{\textbf{in}};\mathrm{e}^{(1)}_{\mathrm{n}}\Big\rangle_{\mathbb{L}^{2}\Big(\Omega \Big)} = 0.
	\end{align}
	Hence combining equations (\ref{pr1}), (\ref{eq:incident}) and plugging these into the equation (\ref{eq:div}) we obtain
	\begin{align}\label{eq:1}
		\Big\langle \nabla\mathrm{H};\mathrm{e}^{(1)}_{\mathrm{n}}\Big\rangle_{\mathbb{L}^{2}(\Omega)} = \Big\langle \mathbb{T}\Big[\nabla\mathrm{H}\Big];\ \mathrm{e}^{(1)}_\mathrm{n} \Big\rangle_{\mathbb{L}^2(\Omega)}.  
	\end{align}
	
	\subsection{Projection on $\textbf{H}_{0}(\text{curl},0)$}
	Considering the set of orthonormal basis $\big(\mathrm{e}^{(2)}_{\mathrm{n}}\big)_{\mathrm{n} \in \mathbb{N}}$ and taking the projection of the equation (\ref{eq:laplace0}) into the subspace $\textbf{H}_{0}(\textbf{curl},0)$ we obtain
	\begin{align}\label{eq:curl}
		\Big\langle \nabla \mathrm{H};\mathrm{e}^{(2)}_{\mathrm{n}}\Big\rangle_{\mathbb{L}^{2}(\Omega)}-\alpha \Big\langle \nabla \int_\Omega\nabla\mathbb{G}^{(0)}(\mathrm{x},\mathrm{y})\cdot\nabla \mathrm{H}(\mathrm{y}) d\mathrm{y};\mathrm{e}^{(2)}_{\mathrm{n}}\Big\rangle_{\mathbb{L}^{2}(\Omega)} &= \Big\langle \nabla \mathrm{H}^{\textbf{in}};\mathrm{e}^{(2)}_{\mathrm{n}}\Big\rangle_{\mathbb{L}^{2}(\Omega)} + \alpha \Big\langle \mathbb{T}\Big[\nabla\mathrm{H}\Big];\mathrm{e}^{(2)}_{\mathrm{n}}\Big\rangle_{\mathbb{L}^{2}(\Omega)}
	\end{align}
	Next, we evaluate the second term of the left hand of the equation (\ref{eq:curl}) with the usual induced inner product of $\mathbb{L}^{2}$. As before, we consider the second sub-space $\textbf{H}_{0}(\textbf{curl},0)$. We use the vector identity $\nabla\nabla\cdot\textbf{f} = \nabla\times\nabla\times\textbf{f}+\Delta\textbf{f}$ and the vector integration formula to deduce that
	\begin{align}\label{eq:proj2}
		\nonumber
		\Big\langle \nabla\int_\Omega\nabla\mathbb{G}^{(0)}(\mathrm{x},\mathrm{y})\cdot\nabla \mathrm{H}(\mathrm{y}) d\mathrm{y}\ ;\ \mathrm{e}^{(2)}_{\mathrm{n}}\Big\rangle_{\textbf{L}^{2}(\Omega)} 
		&= \int_{\Omega} \mathrm{e}^{(2)}_{\mathrm{n}} \nabla \int_\Omega\nabla\mathbb{G}^{(0)}(\mathrm{x},\mathrm{y})\cdot\nabla \mathrm{H}(\mathrm{y}) d\mathrm{y}\\
		\nonumber
		&= \int_{\Omega} \mathrm{e}^{(2)}_{\mathrm{n}}\cdot \ \nabla \times\Big(\nabla\times \int_\Omega\mathbb{G}^{(0)}(\mathrm{x},\mathrm{y})\nabla \mathrm{H}(\mathrm{y}) d\mathrm{y}\Big) \\ \nonumber &+ \int_{\Omega} \mathrm{e}^{(2)}_{\mathrm{n}} \  \Delta\int_\Omega\mathbb{G}^{(0)}(\mathrm{x},\mathrm{y})\nabla \mathrm{H}(\mathrm{y}) d\mathrm{y}\\
		\nonumber
		&= \int_{\Omega} \nabla \times \mathrm{e}^{(2)}_{\mathrm{n}} \cdot \nabla\times \int_\Omega\mathbb{G}^{(0)}(\mathrm{x},\mathrm{y})\nabla \mathrm{H}(\mathrm{y}) d\mathrm{y} \\ \nonumber &+ \int_{\partial\Omega} \nu \times \mathrm{e}^{(2)}_{\mathrm{n}}\cdot  \nabla\times \int_\Omega\mathbb{G}^{(0)}(\mathrm{x},\mathrm{y})\nabla \mathrm{H}(\mathrm{y}) d\mathrm{y} \\ \nonumber &+ \int_{\Omega} \mathrm{e}^{(2)}_{\mathrm{n}} \Delta \int_\Omega\mathbb{G}^{(0)}(\mathrm{x},\mathrm{y})\nabla \mathrm{H}(\mathrm{y}) d\mathrm{y}
		\\ \nonumber &= \int_{\Omega} \mathrm{e}^{(2)}_{\mathrm{n}} \ \Delta \int_\Omega\mathbb{G}^{(0)}(\mathrm{x},\mathrm{y})\nabla \mathrm{H}(\mathrm{y}) d\mathrm{y} \\ &= \int_{\Omega} \mathrm{e}^{(2)}_{\mathrm{n}}\cdot \nabla \mathrm{H}(\mathrm{y}) d\mathrm{y}  = \Big\langle \nabla \mathrm{H};\mathrm{e}^{(2)}_{\mathrm{n}}\Big\rangle_{\mathbb{L}^{2}(\Omega)}.
	\end{align}
	Therefore, from the equation (\ref{eq:proj2}) and plugging it in (\ref{eq:curl}) we get
	\begin{align}\label{eq:2}
		(1-\alpha)\Big\langle \nabla \mathrm{H};\mathrm{e}^{(2)}_{\mathrm{n}}\Big\rangle_{\mathbb{L}^{2}(\Omega)} \ = \ \Big\langle \nabla \mathrm{H}^{\textbf{in}};\mathrm{e}^{(2)}_{\mathrm{n}}\Big\rangle_{\mathbb{L}^{2}\Big(\Omega \Big)} + \Big\langle \mathbb{T}\Big[ \nabla \mathrm{H}\Big];\ \mathrm{e}^{(2)}_\mathrm{n} \Big\rangle_{\mathbb{L}^{2}(\Omega)}.
	\end{align}
	\subsection{Projection on $\nabla \mathrm{H}_{\textbf{arm}}$}
	We know that the magnetization operator $\mathbb{M}$ as define by (\ref{magnetic operator}) is a compact, positive and self-adjoint operator on $\mathbb{L}^{2}$. Therefore, the eigenfunctions of $\mathbb{M}:\nabla \mathrm{H}_{\textbf{arm}} \to \nabla \mathrm{H}_{\textbf{arm}}$ form a complete orthonormal system in $\nabla \mathrm{H}_{\textbf{arm}}$. Let us consider $(\lambda^{(3)}_{\mathrm{n}},\mathrm{e}^{(3)}_{\mathrm{n}})_{\mathrm{n}\in \mathbb{N}}$ be the orhonormal eigensystem of the magnetic operator in $\nabla \mathrm{H}_{\textbf{arm}}$. More details about the properties of magnetic operator can be found in \cite{friedmanI}, \cite{friedmanII}, \cite{friedmanIII}. Using this orthonormal eigensystem in  (\ref{eq:laplace0}), we obtain
	\begin{align}\label{eq:laplace3}
		\Big\langle \nabla \mathrm{H};\mathrm{e}^{(3)}_{\mathrm{n}}\Big\rangle_{\mathbb{L}^{2}(\Omega)} - \alpha \Big\langle \nabla \int_\Omega\nabla\mathbb{G}^{(0)}(\mathrm{x},\mathrm{y})\cdot\nabla \mathrm{H}(\mathrm{y}) d\mathrm{y};\mathrm{e}^{(3)}_{\mathrm{n}}\Big\rangle_{\mathbb{L}^{2}(\Omega)} &=\Big\langle \nabla \mathrm{H}^{\textbf{in}};\mathrm{e}^{(3)}_{\mathrm{n}}\Big\rangle_{\mathbb{L}^{2}(\Omega)} + \alpha \Big\langle \mathbb{T}\Big[\nabla\mathrm{H}\Big]; \mathrm{e}^{(3)}_{\mathrm{n}}\Big\rangle_{\mathbb{L}^{2}(\Omega)} 
	\end{align}
	With the self-adjointness of the magnetic operator, we evaluate the second term of the equation (\ref{eq:laplace3})
	\begin{align}\label{eq:proj3}
		\Bigg\langle \nabla \int_\Omega\nabla\mathbb{G}^{(0)}(\mathrm{x},\mathrm{y})\cdot\nabla \mathrm{H}(\mathrm{y}) d\mathrm{y}\ ;\ \mathrm{e}^{(3)}_{\mathrm{n}}\Bigg\rangle_{\mathbb{L}^{2}\Big(\Omega \Big)} \nonumber
		&= \Big\langle \nabla \mathrm{H};\ \nabla\int_\Omega\nabla\mathbb{G}^{(0)}(\mathrm{x},\mathrm{y})\cdot \mathrm{e}^{(3)}_{\mathrm{n}}d\mathrm{y}\Big\rangle_{\textbf{L}^{2}\Big(\Omega \Big)}\\ \nonumber
		&= \Big\langle \nabla \mathrm{H};\ \lambda^{(3)}_{\mathrm{n}}\mathrm{e}^{(3)}_{\mathrm{n}}\Big\rangle_{\mathbb{L}^{2}\Big(\Omega \Big)}\\
		&= \lambda^{(3)}_{\mathrm{n}}\Big\langle \nabla \mathrm{H};\mathrm{e}^{(3)}_{\mathrm{n}}\Big\rangle_{\mathbb{L}^{2}\Big(\Omega \Big)}.
	\end{align}
	Consequently, we obtain from the equations (\ref{eq:proj3}) and plugging it into the equation (\ref{eq:laplace3})
	\begin{align}\label{eq:3}
		(1-\alpha\lambda^{(3)}_{\mathrm{n}})\Big\langle \nabla \mathrm{H};\mathrm{e}^{(3)}_{\mathrm{n}}\Big\rangle_{\textbf{L}^{2}(\Omega)} \ = \ \Big\langle \nabla \mathrm{H}^{\textbf{in}};\mathrm{e}^{(3)}_{\mathrm{n}}\Big\rangle_{\mathbb{L}^{2}(\Omega )} + \Big\langle \mathbb{T}\Big[ \nabla \mathrm{H}\Big];\ \mathrm{e}^{(3)}_\mathrm{n} \Big\rangle_{\mathbb{L}^{2}(\Omega)}.
	\end{align}
	\bigbreak
	\noindent
	Due to this, in order to estimate the error term associated with each equation arising from projection of the equation (\ref{eq:laplace0}) onto the described sub-spaces we need an a priori estimate which we provide in the following proposition
	\begin{proposition}\label{prop}
		We have 
		\begin{align}
			\Big\Vert \nabla \mathrm{H} \Big\Vert_{\mathbb{L}^{2}(\Omega)} = \mathcal{O}\Big(\delta^{1-\mathrm{h}}\Big) \quad \text{for} \quad 0<\mathrm{h}<1.
		\end{align}
	\end{proposition}
	\noindent
	\textbf{Proof.} See section \ref{apriori}.
	\bigbreak
	\noindent
	Based on these orthonormal bases we can arrive at a basis for the whole space $\mathbb{L}^{2}(\Omega)$. Consider the complete orthonormal basis $\Big\{\big(\mathrm{e}^{(1)}_{\mathrm{n}}\big), \big(\mathrm{e}^{(2)}_{\mathrm{n}}\big), \big(\mathrm{e}^{(3)}_{\mathrm{n}}\big) \Big\}_{\mathrm{n} \in \mathbb{N}}.$ In order to estimate $\nabla \mathrm{H}$, we use the Parseval's  identity
	\begin{align}\label{parseval}
		\Big\Vert \nabla \mathrm{H} \Big\Vert_{\mathbb{L}^2(\Omega)}^2 = \sum_{\substack{\mathrm{n}\in \mathbb{N}\\
				\mathrm{i}=1,2,3}}\Big|\Big\langle \nabla \mathrm{H};\mathrm{e}^{(\mathrm{i})}_{\mathrm{n}}\Big\rangle_{\mathbb{L}^{2}(\Omega)}\Big|^2.
	\end{align}
	\noindent
	As a next step, we examine the equation (\ref{eq:3}) and rewrite it in the scaled domain $\mathrm{B}$
	\begin{align}
		\Big(1-\alpha \lambda^{(3)}_\mathrm{n}\Big)\Big\langle\widehat{\nabla \mathrm{H}};\mathrm{e}^{(3)}_\mathrm{n}\Big\rangle_{\mathbb{L}^2(\mathrm{B})} = \big\langle\widehat{\nabla\mathrm{H}}^{\textbf{in}};\mathrm{e}^{(3)}_\mathrm{n}\Big\rangle_{\mathbb{L}^2(\mathrm{B})} + \alpha \delta^2\Big\langle \mathbb{T}_\delta\Big[\widehat{\nabla\mathrm{H}}\Big];\mathrm{e}^{(3)}_\mathrm{n}\Big\rangle_{\mathbb{L}^2(\mathrm{B})}
	\end{align}
	and derive the estimate
	\begin{align}\label{f1}
		\nonumber
		\Big|\Big\langle \widehat{\nabla \mathrm{H}};\mathrm{e}^{(3)}_{\mathrm{n}}\Big\rangle_{\mathbb{L}^{2}(\mathrm{B})}\Big| = \dfrac{1}{ |1 - \alpha\lambda^{(3)}_{\mathrm{n}}|} \Bigg[\Big|\Big\langle \widehat{\nabla \mathrm{H}}^{\textbf{in}};\mathrm{e}^{(3)}_{\mathrm{n}}\Big\rangle_{\mathbb{L}^{2}(\mathrm{B})}\Big| + |\alpha|\delta^2\Big|\Big\langle \mathbb{T}_\delta\Big[\widehat{ \nabla \mathrm{H}}\Big];\ \mathrm{e}^{(3)}_\mathrm{n} \Big\rangle_{\mathbb{L}^{2}(\mathrm{B})}\Big|\Bigg].\\
		\Longrightarrow \Big|\Big\langle \widehat{\nabla \mathrm{H}};\mathrm{e}^{(3)}_{\mathrm{n}}\Big\rangle_{\mathbb{L}^{2}(\mathrm{B})}\Big|^2 \lesssim \dfrac{1}{ |1 - \alpha\lambda^{(3)}_{\mathrm{n}}|^2}\Bigg[ \Big|\Big\langle \widehat{\nabla \mathrm{H}}^{\textbf{in}};\mathrm{e}^{(3)}_{\mathrm{n}}\Big\rangle_{\mathbb{L}^{2}(\mathrm{B})}\Big|^2 + \delta^4\Big|\Big\langle \mathbb{T}_\delta\Big[ \widehat{\nabla \mathrm{H}}\Big];\ \mathrm{e}^{(3)}_\mathrm{n} \Big\rangle_{\mathbb{L}^{2}(\mathrm{B})}\Big|^2\Bigg].
	\end{align}
	In a similar way, we obtain from (\ref{eq:2}) and (\ref{eq:1})
	\begin{align}\label{f2}
		\nonumber
		\Big|\Big\langle \widehat{\nabla \mathrm{H}};\mathrm{e}^{(2)}_{\mathrm{n}}\Big\rangle_{\mathbb{L}^{2}(\mathrm{B})}\Big| = \dfrac{1}{ |1 - \alpha|^2}\Bigg[ \Big|\Big\langle \widehat{\nabla \mathrm{H}}^{\textbf{in}};\mathrm{e}^{(2)}_{\mathrm{n}}\Big\rangle_{\mathbb{L}^{2}(\mathrm{B})}\Big| + |\alpha|\delta^2\Big|\Big\langle \mathbb{T}_\delta\Big[ \widehat{\nabla \mathrm{H}}\Big];\ \mathrm{e}^{(2)}_\mathrm{n} \Big\rangle_{\mathbb{L}^{2}(\mathrm{B})}\Big|\Bigg].\\
		\Longrightarrow \Big|\Big\langle \widehat{\nabla \mathrm{H}};\mathrm{e}^{(2)}_{\mathrm{n}}\Big\rangle_{\mathbb{L}^{2}(\mathrm{B})}\Big|^2 \lesssim \dfrac{1}{ |1 - \alpha|^2}\Bigg[ \Big|\Big\langle \widehat{\nabla \mathrm{H}}^{\textbf{in}};\mathrm{e}^{(2)}_{\mathrm{n}}\Big\rangle_{\mathbb{L}^{2}(\mathrm{B})}\Big|^2 + \delta^4\Big|\Big\langle \mathbb{T}\Big[ \widehat{\nabla \mathrm{H}}\Big];\ \mathrm{e}^{(2)}_\mathrm{n} \Big\rangle_{\mathbb{L}^{2}(\mathrm{B})}\Big|^2\Bigg],   
	\end{align}
	and 
	\begin{align}\label{f3}
		\Big|\Big\langle \widehat{\nabla \mathrm{H}};\mathrm{e}^{(1)}_{\mathrm{n}}\Big\rangle_{\mathbb{L}^{2}(\mathrm{B})}\Big|^2 \lesssim \delta^4\Big|\Big\langle \mathbb{T}_\delta\Big[ \widehat{\nabla \mathrm{H}}\Big];\ \mathrm{e}^{(1)}_\mathrm{n} \Big\rangle_{\mathbb{L}^{2}(\mathrm{B})}\Big|^2
	\end{align}
	respectively.
	\bigbreak
	\noindent
	Furthermore, we rewrite the expression (\ref{parseval}) as follows
	\begin{align}\label{par1}
		\displaystyle \int_\mathrm{B} |\widehat{\nabla\mathrm{H}}|^2(\mathrm{\eta})d\mathrm{\eta} = \Big|\Big\langle\nabla\mathrm{H};\mathrm{e}^{(\mathrm{i})}_\mathrm{n_0}\big\rangle_{\mathbb{L}^2(\mathrm{B})}\Big|^2 + \sum_{\substack{\mathrm{n}\ne\mathrm{n}_0\\ \mathrm{i}=1,2,3}} \Big|\Big\langle\widehat{\nabla\mathrm{H}};\mathrm{e}^{(\mathrm{i})}_\mathrm{n}\big\rangle_{\mathbb{L}^2(\mathrm{B})}\Big|^2.
	\end{align}
	Let us recall that from the Lorentz model we can deduce the following, see Section \ref{plasmonics} for details,
	\begin{align}\label{condition}
		\big|1-\alpha\lambda^{(3)}_{\mathrm{n}}\Big| = 
		\begin{cases}
			\delta^\mathrm{h} & \mathrm{n} = \mathrm{n}_0 \\
			1 & \mathrm{n} \ne \mathrm{n}_0.
		\end{cases}   
	\end{align}
	Therefore, considering the estimates (\ref{f1}), (\ref{f2}), and (\ref{f3}) and plugging these results into the previous series expansion (\ref{par1}) we obtain
	\begin{align}\label{prefinal}
		\nonumber
		\int_\mathrm{B}|\widehat{\nabla \mathrm{H}}|^2(\mathrm{\eta})d\mathrm{\eta} &= \dfrac{1}{ |1 - \alpha\lambda^{(3)}_{\mathrm{n}_{0}}|^2}\Big|\Big\langle \widehat{\nabla \mathrm{H}}^{\textbf{in}};e^{(3)}_{\mathrm{n}_{0}}\Big\rangle_{\mathbb{L}^{2}(\mathrm{B})}\Big|^2 + \dfrac{1}{ |1 - \alpha|^2}\Big|\Big\langle\widehat{ \nabla \mathrm{H}}^{\textbf{in}};e^{(2)}_{\mathrm{n}_{0}}\Big\rangle_{\mathbb{L}^{2}(\mathrm{B})}\Big|^2 + \sum_{\substack{\mathrm{n}\ne\mathrm{n}_0}} \frac{\Big|\Big\langle\widehat{\nabla\mathrm{H}}^{\textbf{in}};\mathrm{e}^{(\mathrm{3})}_\mathrm{n}\big\rangle_{\mathbb{L}^2(\mathrm{B})}\Big|^2}{|1 - \alpha\lambda^{(3)}_\mathrm{n}|^2} \\&+ \sum_{\substack{\mathrm{n}\ne\mathrm{n}_0}} \frac{\Big|\Big\langle\widehat{\nabla\mathrm{H}};\mathrm{e}^{(\mathrm{2})}_\mathrm{n}\big\rangle_{\mathbb{L}^2(\mathrm{B})}\Big|^2}{|1 - \alpha|^2} + \mathcal{O}\Big(\dfrac{1}{ |1 - \alpha\lambda^{(3)}_{\mathrm{n}_{0}}|^2}\delta^4\Big\Vert\mathbb{T}_\delta\Big[\widehat{\nabla\mathrm{H}}\Big]\Big\Vert^2_{\mathbb{L}^2(\mathrm{B})}\Big).
	\end{align}
	As a next step, we deal with the second term of the previous expression (\ref{prefinal}) i.e. $\Big\langle \nabla \mathrm{H}^{\textbf{in}};e^{(2)}_{\mathrm{n}_{0}}\Big\rangle_{\mathbb{L}^{2}(\Omega)}.$ We have the following lemma
	\begin{lemma}\label{l51}
		The eigen-function $\mathrm{e}^{(2)}_{\mathrm{n}}(\cdot)$ satisfies the following mean vanishing integral property
		\begin{align}
			\int_\Omega \mathrm{e}^{(2)}_{\mathrm{n}}(\mathrm{x)}d\mathrm{x} = 0.
		\end{align}
	\end{lemma} 
	\noindent
	Since, $\mathrm{H}^{\textbf{in}}$ is smooth, we can use Taylor's series expansion, and hence we observe that the dominating term is related to $\displaystyle\int_\Omega \mathrm{e}^{(2)}_{\mathrm{n}}(\mathrm{x)}d\mathrm{x}$, which is zero by the Lemma \ref{l51}. Thus we deduce that
	\begin{align}\label{vanish}
		\Bigg|\Big\langle \widehat{\nabla \mathrm{H}}^{\textbf{in}};\mathrm{e}^{(2)}_{\mathrm{n}_{0}}\Big\rangle_{\mathbb{L}^{2}(\mathrm{B})}\Bigg|^2 = \mathcal{O}\Big(\delta^2\Big).
	\end{align}
	Taylor's expansion provide
	\begin{align}\label{taylor}
		\Big\langle \widehat{\nabla\mathrm{H}}^{\textbf{in}};\mathrm{e}^{(3)}_{\mathrm{n}_{0}}\Big\rangle_{\mathbb{L}^2(\mathrm{B})} =  \nabla\mathrm{H}^{\textbf{in}}(\mathrm{z})\int_\mathrm{B}\Tilde{\mathrm{e}}^{(3)}_{\mathrm{n}_{0}}(\mathrm{x})d\mathrm{x} + \mathcal{O}(\delta).
	\end{align}
	Also, we have 
	\begin{align} \label{taylor1}
		\sum_{\substack{\mathrm{n}\ne\mathrm{n}_0}} \frac{\Big|\Big\langle\widehat{\nabla\mathrm{H}}^{\textbf{in}};\mathrm{e}^{(\mathrm{3})}_\mathrm{n}\Big\rangle_{\mathbb{L}^2(\mathrm{B})}\Big|^2}{|1 - \alpha\lambda^{(3)}_\mathrm{n}|^2} \sim 1.
	\end{align}
	Combining three above estimates (\ref{vanish}), (\ref{taylor}), (\ref{taylor1}) and plugging this into (\ref{prefinal}) we obtain
	\begin{align}
		\Big\Vert \widehat{\nabla\mathrm{H}}\Big\Vert^2_{\mathbb{L}^2(\mathrm{B})} = \dfrac{1}{ |1 - \alpha\lambda^{(3)}_{\mathrm{n}_{0}}|^2}\Big|\Big\langle \widehat{\nabla \mathrm{H}}^{\textbf{in}};e^{(3)}_{\mathrm{n}_{0}}\Big\rangle_{\mathbb{L}^{2}(\mathrm{B})}\Big|^2 + \mathcal{O}\Big(1\Big).
	\end{align}
	Consequently, from the fact $|\nabla \mathrm{H}|^2 = |\mathrm{E}|^2$ and when we scale back to $\Omega$, we derive that 
	\begin{align}
		\int_{\Omega}|\mathrm{E}|^2(\mathrm{y})d\mathrm{y} = \dfrac{1}{ |1 - \alpha\lambda^{(3)}_{\mathrm{n}_{0}}|^2}\Big|\Big\langle \widehat{\nabla \mathrm{H}}^{\textbf{in}};e^{(3)}_{\mathrm{n}_{0}}\Big\rangle_{\mathbb{L}^{2}(\mathrm{B})}\Big|^2 \delta^2 + \mathcal{O}\Big(\delta^{2}\Big).
	\end{align}
	Moreover, after plugging the estimate (\ref{taylor}) in the above expression, we obtain
	\begin{align}
		\int_{\Omega}|\mathrm{E}|^2(\mathrm{y})d\mathrm{y} = \dfrac{1}{ |1 - \alpha\lambda^{(3)}_{\mathrm{n}_{0}}|^2}\Bigg[|\nabla\mathrm{H}^{\textbf{in}}|^2(\mathrm{z})\Big(\int_\mathrm{B}\Tilde{\mathrm{e}}^{(3)}_{\mathrm{n}_{0}}(\mathrm{x})d\mathrm{x}\Big)^2\delta^2 + \mathcal{O}\Big(\delta^{3}\Big)\Bigg] \quad \text{with} \quad \mathrm{h}<1.
	\end{align}
	This completes the proof of Theorem \ref{th1}.
	
	\subsection{\textbf{A Priori Estimate }}\label{apriori} This section is divided into two subsections providing the proof for Proposition \ref{prop} and Proposition \ref{prop2} respectively.
	\subsubsection{\textbf{Proof of Proposition} \ref{prop}}
	\noindent
	We need to identify the exact dominating term in the equations (\ref{eq:1}), (\ref{eq:2}), and (\ref{eq:3}). We start with the equation 
	\begin{align}\label{eqofinterest}
		\nabla \mathrm{H}(\mathrm{x}) - \alpha \ \nabla\int_\Omega\nabla\mathbb{G}^{(0)}(\mathrm{x},\mathrm{y})\cdot\nabla \mathrm{H}(\mathrm{y}) d\mathrm{y} = \nabla\mathrm{H}^{\textbf{in}}(\mathrm{x}) + \alpha \mathbb{T}\Big[\nabla\mathrm{H}\Big].  
	\end{align}
	It can be showed that from the definition (\ref{expansion21}) that
	\begin{align}
		\mathbb{T}\Big[\nabla\mathrm{H}\Big] \approx \delta^2  \mathbb{T}_\delta\Big[\widehat{\nabla\mathrm{H}}\Big]
	\end{align}
	therefore, after rewriting the equation (\ref{eqofinterest}) in the scaled domain $\mathrm{B}$, we obtain
	\begin{align}\label{eqofinterest1}
		\widehat{\nabla \mathrm{H}} - \alpha \ \nabla\int_\mathrm{B}\nabla\mathbb{G}^{(0)}(\mathrm{\xi},\mathrm{\eta})\cdot\widehat{\nabla \mathrm{H}}(\mathrm{\eta}) d\mathrm{\eta} = \widehat{\nabla\mathrm{H}}^{\textbf{in}} + \alpha \delta^2 \mathbb{T}_\delta\Big[\widehat{\nabla\mathrm{H}}\Big].    
	\end{align}
	Now, we estimate $\nabla \mathrm{H}$ on the sub-space $\nabla \mathrm{H}_{\textbf{arm}}$ .
	Taking the induced $\textbf{L}^{2}$-norm for the equation (\ref{eqofinterest}), we get
	\begin{align}\label{eqofinterest12}
		\Big\langle\widehat{\nabla \mathrm{H}};\mathrm{e}^{(3)}_\mathrm{n}\Big\rangle_{\mathbb{L}^2(\mathrm{B})} - \alpha \Big\langle \nabla\int_\mathrm{B}\nabla\mathbb{G}^{(0)}(\mathrm{\xi},\mathrm{\eta})\cdot\widehat{\nabla \mathrm{H}}(\mathrm{\eta}) d\mathrm{\eta};\mathrm{e}^{(3)}_\mathrm{n}\Big\rangle_{\mathbb{L}^2(\mathrm{B})} = \big\langle\widehat{\nabla\mathrm{H}}^{\textbf{in}};\mathrm{e}^{(3)}_\mathrm{n}\Big\rangle_{\mathbb{L}^2(\mathrm{B})} + \alpha \delta^2\Big\langle \mathbb{T}_\delta\Big[\widehat{\nabla\mathrm{H}}\Big];\mathrm{e}^{(3)}_\mathrm{n}\Big\rangle_{\mathbb{L}^2(\mathrm{B})}.
	\end{align}
	Then due to the self-adjointness of the magnetization operator $\mathbb{M}$ and the orthonormal system $\Big(\lambda^{(3)}_\mathrm{n}, \mathrm{e}^{(3)}_\mathrm{n}\Big)$, we obtain
	\begin{align}
		\Big(1-\alpha \lambda^{(3)}_\mathrm{n}\Big)\Big\langle\widehat{\nabla \mathrm{H}};\mathrm{e}^{(3)}_\mathrm{n}\Big\rangle_{\mathbb{L}^2(\mathrm{B})} = \big\langle\widehat{\nabla\mathrm{H}}^{\textbf{in}};\mathrm{e}^{(3)}_\mathrm{n}\Big\rangle_{\mathbb{L}^2(\mathrm{B})} + \alpha \delta^2\Big\langle \mathbb{T}_\delta\Big[\widehat{\nabla\mathrm{H}}\Big];\mathrm{e}^{(3)}_\mathrm{n}\Big\rangle_{\mathbb{L}^2(\mathrm{B})}
	\end{align}
	and then
	\begin{align}\label{c1}
		\sum_{\mathrm{n}\in \mathbb{N}}\Big|\Big\langle\widehat{\nabla \mathrm{H}};\mathrm{e}^{(3)}_\mathrm{n}\Big\rangle_{\mathbb{L}^2(\mathrm{B})}\Big|^2 = \delta^{-2\mathrm{h}}\sum_{\mathrm{n}\in \mathbb{N}}\Big|\Big\langle\widehat{\nabla\mathrm{H}}^{\textbf{in}};\mathrm{e}^{(3)}_\mathrm{n}\Big\rangle_{\mathbb{L}^2(\mathrm{B})}\Big|^2 + |\alpha|^2 \delta^{4-2\mathrm{h}}\sum_{\mathrm{n}\in \mathbb{N}}\Big|\Big\langle \mathbb{T}_\delta\Big[\widehat{\nabla\mathrm{H}}\Big];\mathrm{e}^{(3)}_\mathrm{n}\Big\rangle_{\mathbb{L}^2(\mathrm{B})}\Big|^2.
	\end{align}
	In a similar way as above, if we consider the orhonormal systems $\Big( \mathrm{e}^{(1)}_\mathrm{n}\Big)$ and $\Big( \mathrm{e}^{(2)}_\mathrm{n}\Big)$, we deduce from (\ref{eqofinterest1})
	\begin{align}\label{c2}
		\sum_{\mathrm{n}\in \mathbb{N}}\Big|\Big\langle\widehat{\nabla \mathrm{H}};\mathrm{e}^{(1)}_\mathrm{n}\Big\rangle_{\mathbb{L}^2(\mathrm{B})}\Big|^2 = |\alpha|^2 \delta^{4}\sum_{\mathrm{n}\in \mathbb{N}}\Big|\Big\langle \mathbb{T}_\delta\Big[\widehat{\nabla\mathrm{H}}\Big];\mathrm{e}^{(3)}_\mathrm{n}\Big\rangle_{\mathbb{L}^2(\mathrm{B})}\Big|^2
	\end{align}
	and 
	\begin{align}\label{c3}
		\sum_{\mathrm{n}\in \mathbb{N}}\Big|\Big\langle\widehat{\nabla \mathrm{H}};\mathrm{e}^{(2)}_\mathrm{n}\Big\rangle_{\mathbb{L}^2(\mathrm{B})}\Big|^2 = \frac{1}{|1-\alpha|^2}\sum_{\mathrm{n}\in \mathbb{N}}\Big|\Big\langle\widehat{\nabla\mathrm{H}}^{\textbf{in}};\mathrm{e}^{(2)}_\mathrm{n}\Big\rangle_{\mathbb{L}^2(\mathrm{B})}\Big|^2 + \frac{|\alpha|^2 \delta^{4}}{|1-\alpha|^2}\sum_{\mathrm{n}\in \mathbb{N}}\Big|\Big\langle \mathbb{T}_\delta\Big[\widehat{\nabla\mathrm{H}}\Big];\mathrm{e}^{(2)}_\mathrm{n}\Big\rangle_{\mathbb{L}^2(\mathrm{B})}\Big|^2,
	\end{align}
	respectively.
	\bigbreak
	\noindent
	Therefore, combining the expressions (\ref{c1}), (\ref{c2}), (\ref{c3}) and due to the boundedness of the operator $\mathbb{T}: \mathbb{L}^2(\mathrm{B}) \to \mathbb{H}^2(\mathrm{B})$, we obtain
	\begin{align}
		\Big\Vert \widehat{\nabla\mathrm{H}}\Big\Vert^2_{\mathbb{L}^2(\mathrm{B})} \lesssim \delta^{-2\mathrm{h}}\Big\Vert \widehat{\nabla\mathrm{H}}^{\textbf{in}}\Big\Vert^2_{\mathbb{L}^2(\mathrm{B})} + \delta^{4-2\mathrm{h}} \Big\Vert \widehat{\nabla\mathrm{H}}\Big\Vert^2_{\mathbb{L}^2(\mathrm{B})}
	\end{align}
	hence
	\begin{align}\label{apr}
		\Big\Vert \widehat{\nabla\mathrm{H}}\Big\Vert^2_{\mathbb{L}^2(\mathrm{B})} \lesssim \delta^{-2\mathrm{h}}\Big\Vert \widehat{\nabla\mathrm{H}}^{\textbf{in}}\Big\Vert^2_{\mathbb{L}^2(\mathrm{B})}.
	\end{align}
	Moreover, as the incident wave is smooth enough, so we deduce 
	\begin{equation}\label{p9}
		\Big\Vert \nabla \mathrm{H}^{\text{in}}\Big\Vert^{2}_{\mathbb{L}^{2}\Big(\Omega \Big)} = \mathcal{O}(\delta^{2}).
	\end{equation}
	Combining (\ref{apr}) and (\ref{p9}), we get
	\begin{align}\label{a priori}
		\Big\Vert\nabla \mathrm{H}\Big\Vert_{\mathbb{L}^{2}\Big(\Omega \Big)}^{2} = \mathcal{O}\Big(\delta^{2-2h}\Big)  \ \text{with the condition} \ h<1.
	\end{align}
	This completes the proof of Proposition \ref{prop}.
	\bigbreak
	\noindent
	Next, we want to make a short remark, which will be essential to describe the proof of Theorem \ref{th2}. As $|\nabla \mathrm{H}|^2 = |\mathrm{E}|^2$, we obtain
	\begin{align}
		\Big\Vert E \Big\Vert_{{L}^2\Big(\Omega \Big)} = \mathcal{O}\Big(\delta^{1-h}\Big) \quad \text{with the condition} \quad \mathrm{h}<1.
	\end{align}
	
	\subsubsection{\textbf{Proof of Proposition \ref{prop2}}} Let us state the following lemma before proceeding to give the proof of Proposition \ref{prop2}.
	\begin{lemma} \label{6.2}
		Suppose $0<\delta\le 1$ and $\Omega = \delta \mathrm{B} + x^*.$ Then for $\nabla\varphi \in \mathbb{L}^2(\Omega)$, we have
		\begin{align}\label{lemma 3.1}
			\Vert \nabla\varphi \Vert_{\mathbb{L}^2(\Omega)} = \Vert \nabla \hat{\varphi} \Vert_{\mathbb{L}^2(\mathrm{B})} \quad \text{and} \quad \Vert \nabla\varphi \Vert_{\mathbb{L}^4(\Omega)} = \delta^{-\frac{1}{2}}\Vert \nabla \hat{\varphi} \Vert_{\mathbb{L}^4(\mathrm{B})}.
		\end{align} 
	\end{lemma}
	\noindent
	\textbf{Proof.} Let us suppose $\mathrm{x} = \delta \xi + \mathrm{z}$. Then we obtain for $\nabla\varphi \in \mathbb{L}^2(\Omega)$
	\begin{align}
		\nonumber
		\Vert \nabla\varphi \Vert^2_{\mathbb{L}^2(\Omega)} &=  \int_{\Omega}|\nabla\varphi(\mathrm{x})|^2 d\mathrm{x}
		\\ \nonumber &= \delta^2 \int_{\mathrm{B}}\delta^{-2}|\nabla_{\xi}\varphi(\delta \xi + z)|^2 d\xi
		\\ &=  \Vert \nabla \hat{\varphi} \Vert^2_{\mathbb{L}^2(\mathrm{B})}
	\end{align}
	and for $\nabla\varphi \in \mathbb{L}^4(\Omega)$
	\begin{align}
		\nonumber
		\Vert \nabla\varphi \Vert^4_{\mathbb{L}^4(\Omega)} &=  \int_{\Omega}|\nabla\varphi(\mathrm{x})|^4 d\mathrm{x}
		\\ \nonumber &= \delta^2 \int_{\mathrm{B}}\delta^{-4}|\nabla_{\xi}\varphi(\delta \xi + z)|^4 d\xi
		\\ &=  \delta^{-2}\Vert \nabla \hat{\varphi} \Vert^4_{\mathbb{L}^4(\mathrm{B})},
	\end{align}
	which leads to (\ref{lemma 3.1}). 
	\bigbreak
	\noindent
	Next, we start with the Lippmann-Schwinger equation in order to estimate $\Big\Vert |\mathrm{E}|^2\Big\Vert_{\mathbb{L}^2(\Omega)}$:
	
	\begin{align}\label{eq:laplace}
		\mathrm{H}(\mathrm{x})  -\alpha \int_\Omega\nabla\mathbb{G}^{(0)}(\mathrm{x},\mathrm{y})\cdot\nabla \mathrm{H}(\mathrm{y}) d\mathrm{y} = \mathrm{H}^{\textbf{in}}(\mathrm{x}) +  \alpha \int_\Omega\nabla\Big(\mathbb{G}^{(\mathrm{k})}-\mathbb{G}^{(0)}\Big)(\mathrm{x},\mathrm{y})\cdot\nabla \mathrm{H}(\mathrm{y}) d\mathrm{y}.
	\end{align}
	From the asymptotic expression of the Hankel function as $|\mathrm{x}-\mathrm{y}|\to 0$, we obtain the expansions:
	\begin{align}
		\Big(\mathbb{G}^{(\mathrm{k})} - \mathbb{G}^{(0)}\Big)(\mathrm{x},\mathrm{y}) = a |\mathrm{x}-\mathrm{y}|^2\log{|\mathrm{x}-\mathrm{y}|} + b |\mathrm{x}-\mathrm{y}|^2 + \mathcal{O}(|\mathrm{x}-\mathrm{y}|^3), 
	\end{align}
	and
	\begin{align}
		\nabla\Big(\mathbb{G}^{(\mathrm{k})}- \mathbb{G}^{(0)}\Big)(\mathrm{x},\mathrm{y}) = \mathrm{C} |\mathrm{x}-\mathrm{y}|\log{|\mathrm{x}-\mathrm{y}|} + \mathrm{D} |\mathrm{x}-\mathrm{y}| + \mathcal{O}(|\mathrm{x}-\mathrm{y}|^2),
	\end{align}
	as $|\mathrm{x}-\mathrm{y}|\to 0$ with some constants where $\mathrm{a}$ and $\mathrm{b}$, $\mathrm{C}$ and $\mathrm{D}$. Then, we write the above expression in the domain $\mathrm{B}$ i.e. after scaling we obtain
	\begin{align}\label{eq:laplace12}
		\nonumber
		\Tilde{\mathrm{H}}(\xi)  -\alpha \int_\mathrm{B}\nabla_{\eta}\mathbb{G}^{(0)}(\xi,\eta)\cdot\nabla \Tilde{\mathrm{H}}(\eta) d\eta &= \Tilde{\mathrm{H}}^{in}(\xi) +  \alpha\delta^2\Big[ \int_\mathrm{B}|\xi-\eta|\log{|\xi-\eta|}\nabla \Tilde{\mathrm{H}}(\eta) d\eta \\ &+ \log{\delta} \int_\mathrm{B}|\xi-\eta|\nabla \Tilde{\mathrm{H}}(\eta) d\eta
		+ \delta \int_\mathrm{B}|\xi-\eta|\nabla \Tilde{\mathrm{H}}(\eta) d\eta\Big].
	\end{align}
	The Newtonian operator $\mathcal{V}: \mathbb{L}^2(\mathrm{B}) \to \mathbb{H}^2(\mathrm{B})$ is continuous and then we have the continuity of $\nabla\cdot \mathcal{V}: \mathbb{L}^2(\mathrm{B}) \to \mathbb{H}^1(\mathrm{B})$. 
	Consequently, since we have $\Vert\nabla \mathrm{H}\Vert_{\mathbb{L}^{2}(\Omega)}= \Vert\nabla \Tilde{\mathrm{H}}\Vert_{\mathbb{L}^{2}(\mathrm{B})} = \mathcal{O}\Big(\delta^{1-h}\Big)$ and taking $\mathbb{H}^1(\mathrm{B})$-norm of the equation (\ref{eq:laplace12}), we derive
	\begin{align}\label{order1}
		\nonumber
		\Big\Vert \Tilde{\mathrm{H}} \Big\Vert_{\mathbb{H}^1(\mathrm{B})} &\lesssim \Big\Vert \Tilde{\mathrm{H}}^{\textbf{in}} \Big\Vert_{\mathbb{H}^1(\mathrm{B})} + |\alpha| \Big\Vert \nabla \Tilde{\mathrm{H}} \Big\Vert_{\mathbb{L}^2(\mathrm{B})} + |\alpha|\delta^2\Big[ \Big\Vert \nabla \Tilde{\mathrm{H}} \Big\Vert_{\mathbb{L}^2(\mathrm{B})} + |\log{\delta}|\Big\Vert \nabla \Tilde{\mathrm{H}} \Big\Vert_{\mathbb{L}^2(\mathrm{B})} +  \delta\Big\Vert \nabla \Tilde{\mathrm{H}} \Big\Vert_{\mathbb{L}^2(\mathrm{B})}\Big]
		\\ \nonumber &= \mathcal{O}\Big(1\Big) + \mathcal{O}(\delta^{1-\mathrm{h}}) + \mathcal{O}\Big(\delta^{3-\mathrm{h}}\Big) + \mathcal{O}\Big(|\log{\delta}|\delta^{3-\mathrm{h}}\Big) + \mathcal{O}\Big(\delta^{4-\mathrm{h}}\Big)
		\\ &= \mathcal{O}\Big(1\Big).
	\end{align}
	We state the scaled equation (\ref{eq:laplace12})
	\begin{align}\label{eq:laplaceB}
		\nonumber
		\Tilde{\mathrm{H}}(\xi)  -\alpha \int_\mathrm{B}\nabla_{\eta}\mathbb{G}^{(0)}(\xi,\eta)\cdot\nabla \Tilde{\mathrm{H}}(\eta) d\eta &= \Tilde{\mathrm{H}}^{in}(\xi) +  \alpha \delta^2\Big[ \int_\mathrm{B}|\xi-\eta|\log{|\xi-\eta|}\nabla \Tilde{\mathrm{H}}(\eta) d\eta \\ &+ \log{\delta} \int_\mathrm{B}|\xi-\eta|\nabla \Tilde{\mathrm{H}}(\eta) d\eta
		+ \delta \int_\mathrm{B}|\xi-\eta|\nabla \Tilde{\mathrm{H}}(\eta) d\eta\Big]
	\end{align}
	and rewrite it as follows:
	\begin{align}\label{eq:laplaceB}
		\nonumber
		\Tilde{\mathrm{H}}(\xi) + \alpha \int_{\partial \mathrm{B}}\mathbb{G}^{(0)}(\xi,\eta)\partial_{\nu_\eta}\Tilde{\mathrm{H}}(\eta) d\eta + \alpha \delta^2 \delta_{c}\mu_{m}\omega^2 \int_{\mathrm{B}}\mathbb{G}^{(0)}(\xi,\eta)\Tilde{\mathrm{H}}(\eta) d\eta  &= \Tilde{\mathrm{H}}^{in}(\xi) +  \alpha\delta^2\Big[ \delta \int_\mathrm{B}|\xi-\eta|\nabla \Tilde{\mathrm{H}}(\eta) \\ \nonumber &+ \log{\delta} \int_\mathrm{B}|\xi-\eta|\nabla \Tilde{\mathrm{H}}(\eta) d\eta \\ 
		&+ \int_\mathrm{B}|\xi-\eta|\log{|\xi-\eta|}\nabla \Tilde{\mathrm{H}}(\eta) d\eta   d\eta\Big].
	\end{align}
	Taking the normal derivative, as $\xi \to \partial \mathrm{B}$ from inside $\mathrm{B}$, from the jump relations of the derivative of the single layer potential we obtain
	\begin{align}\label{eq:laplaceB}
		\nonumber
		\Big(I+\frac{\alpha}{2}\Big)\partial_{\nu}\Tilde{\mathrm{H}} &= - \alpha \mathbb{K}^{*}_{0}\Big[\partial_{\nu}\Tilde{\mathrm{H}}\Big]  -\alpha \delta^2 \varepsilon_{\mathrm{p}}\mu_{m}\omega^2 \partial_{\nu}\int_\mathrm{B}\mathbb{G}^{(0)}(\xi,\eta)\Tilde{\mathrm{H}}(\eta) d\eta + \partial_{\nu}\Tilde{\mathrm{H}}^{in} +  \alpha \delta^2\Big[\delta \partial_{\nu}\int_\mathrm{B}|\xi-\eta|\nabla \Tilde{\mathrm{H}}(\eta) d\eta \\ &+ \partial_{\nu}\int_\mathrm{B}|\xi-\eta|\ln{|\xi-\eta|}\nabla \Tilde{\mathrm{H}}(\eta) d\eta +  \log{\delta} \partial_{\nu}\int_\mathrm{B}|\xi-\eta|\nabla \Tilde{\mathrm{H}}(\eta) d\eta\Big]
	\end{align}
	then
	\begin{align}\label{eq:laplaceB2}
		\nonumber
		\partial_{\nu}\Tilde{\mathrm{H}} &= - \frac{\alpha}{1+\frac{\alpha}{2}} \mathbb{K}^{*}_{0}\Big[\partial_{\nu}\Tilde{\mathrm{H}}\Big]+ \frac{1}{1+\frac{\alpha}{2}}\partial_{\nu}\Tilde{\mathrm{H}}^{in} +\frac{\alpha \delta^2}{1+\frac{\alpha}{2}}\Bigg[ -\varepsilon_{\mathrm{p}}\mu_{\mathrm{m}}\omega^2 \partial_{\nu}\int_{\mathrm{B}}\mathbb{G}^{(0)}(\xi,\eta)\Tilde{\mathrm{H}}(\eta) d\eta  \\ &+  \partial_{\nu}\int_\mathrm{B}|\xi-\eta|\ln{|\xi-\eta|}\nabla \Tilde{\mathrm{H}}(\eta) d\eta +  \log{\delta} \partial_{\nu}\int_\mathrm{B}|\xi-\eta|\nabla \Tilde{\mathrm{H}}(\eta) d\eta
		+ \delta \partial_{\nu}\int_\mathrm{B}|\xi-\eta|\nabla \Tilde{\mathrm{H}}(\eta) d\eta\Bigg].
	\end{align}
	Next, we state the following lemma, for details we refer to \cite[Theorem 3.6, p. 43]{CK}.
	\begin{lemma}
		If $\mathrm{B}$ is of class $\mathcal{C}^2$, then we have the following regularity
		\begin{align}
			\Bigg\Vert\mathbb{K}^{*}_{0}\Big[\partial_{\nu}\Tilde{\mathrm{H}}\Big]\Big\Vert_{\mathbb{H}^{\frac{1}{2}}(\partial \mathrm{B})} \lesssim \Big\Vert\partial_{\nu}\Tilde{\mathrm{H}}\Big\Vert_{\mathbb{H}^{-\frac{1}{2}}(\partial \mathrm{B})}.
		\end{align}
	\end{lemma} 
	\noindent
	Also, for $\Tilde{\mathrm{H}} \in \mathbb{L}^2(\mathrm{B})$, we have $\displaystyle\int_{B}\mathbb{G}^{(0)}(\xi,\eta)\Tilde{\mathrm{H}}(\eta) d\eta \in \mathbb{H}^2(\mathrm{B})$ and then $\displaystyle\partial_{\nu}\int_{\mathrm{B}}\mathbb{G}^{(0)}(\xi,\eta)\Tilde{\mathrm{H}}(\eta) d\eta \in \mathbb{H}^{\frac{1}{2}}(\mathrm{B})$. Thus, after taking the $\mathbb{H}^{\frac{1}{2}}(\mathrm{B})$-norm on both side of the equation (\ref{eq:laplaceB2}), observing the fact $\Big|\frac{\alpha}{1+\frac{\alpha}{2}}\Big| \sim 1$, and since the integrals on the right hand side of the expression (\ref{eq:laplaceB2}) are also bounded, we deduce the estimate
	\begin{align}\label{order2}
		\nonumber
		\Big\Vert \partial_{\nu}\Tilde{\mathrm{H}} \Big\Vert_{\mathbb{H}^{\frac{1}{2}}(\partial \mathrm{B})} &= \mathcal{O}\Bigg( \Big\Vert \partial_{\nu}\Tilde{\mathrm{H}} \Big\Vert_{\mathbb{H}^{-\frac{1}{2}}(\partial \mathrm{B})} +  \Big\Vert \partial_{\nu}\Tilde{\mathrm{H}}^{\textbf{in}} \Big\Vert_{\mathbb{H}^{\frac{1}{2}}(\partial \mathrm{B})} + \delta^2\Big\Vert \Tilde{\mathrm{H}} \Big\Vert_{\mathbb{L}^2(\mathrm{B})} + \delta^2\Big\Vert \nabla \Tilde{\mathrm{H}} \Big\Vert_{\mathbb{L}^2(\mathrm{B})} + \delta^2|\log{\delta}|\Big\Vert \nabla \Tilde{\mathrm{H}} \Big\Vert_{\mathbb{L}^2(\mathrm{B})} \\ &+  \delta^3\Big\Vert \nabla \Tilde{\mathrm{H}} \Big\Vert_{\mathbb{L}^2(\mathrm{B})}\Bigg).
	\end{align}
	Consequently,
	\begin{align}\label{order2}
		\nonumber
		\Big\Vert \partial_{\nu}\Tilde{\mathrm{H}} \Big\Vert_{\mathbb{H}^{\frac{1}{2}}(\partial \mathrm{B})} &= \mathcal{O}\Bigg( \Big\Vert \Tilde{\mathrm{H}} \Big\Vert_{\mathbb{H}^1(\mathrm{B})} +  \Big\Vert \partial_{\nu}\Tilde{\mathrm{H}}^{\textbf{in}} \Big\Vert_{\mathbb{H}^{\frac{1}{2}}(\partial \mathrm{B})} + \delta^2\Big\Vert \Tilde{\mathrm{H}} \Big\Vert_{\mathbb{L}^2(\mathrm{B})} + \delta^2\Big\Vert \nabla \Tilde{\mathrm{H}} \Big\Vert_{\mathbb{L}^2(\mathrm{B})} + \delta^2|\log{\delta}|\Big\Vert \nabla \Tilde{\mathrm{H}} \Big\Vert_{\mathbb{L}^2(\mathrm{B})} \\ &+  \delta^3\Big\Vert \nabla \Tilde{\mathrm{H}} \Big\Vert_{\mathbb{L}^2(\mathrm{B})}\Bigg).
	\end{align}
	Hence, based on the estimate (\ref{order1}), $\Big\Vert\nabla \mathrm{H}\Big\Vert_{\mathbb{L}^{2}\Big(\Omega \Big)}= \Big\Vert\nabla \Tilde{\mathrm{H}}\Big\Vert_{\mathbb{L}^{2}\Big(\mathrm{B} \Big)} = \mathcal{O}\Big(\delta^{1-h}\Big)$, we obtain
	\begin{align}
		\nonumber
		\Big\Vert \partial_{\nu}\Tilde{\mathrm{H}} \Big\Vert_{\mathbb{H}^{\frac{1}{2}}(\partial \mathrm{B})} &= \mathcal{O}\Big(1\Big) + \mathcal{O}\Big(\delta^2\Big) + \mathcal{O}\Big(|\log{\delta}|\delta^{3-h}\Big) + \mathcal{O}\Big(\delta^{3-h}\Big) + \mathcal{O}\Big(\delta^{4-h}\Big)
		\\ &= \mathcal{O}\Big(1\Big).
	\end{align}
	Moreover, from the Sobolev embedding inequality and considering the well-posed problem satisfied by $\Tilde{\mathrm{H}}$, we deduce
	\begin{align}
		\Big\Vert \nabla \Tilde{\mathrm{H}}\Big\Vert_{\mathbb{L}^4(\mathrm{B})} \lesssim \Big\Vert \Tilde{\mathrm{H}} \Big\Vert_{H^{\frac{3}{2}}(\mathrm{B})} &\lesssim \Big\Vert \partial_{\nu}\Tilde{\mathrm{H}} \Big\Vert_{H^{\frac{1}{2}}(\partial \mathrm{B})} + \delta^2 \Big\Vert \Tilde{\mathrm{H}} \Big\Vert_{\mathbb{L}^2(\partial \mathrm{B})} =\mathcal{O}\Big(1\Big).
	\end{align}
	Hence, by Lemma \ref{6.2}, when we scale back to $\Omega$, we obtain
	\begin{align}
		\nonumber
		\Big\Vert \nabla \mathrm{H}\Big\Vert_{\mathbb{L}^4(\Omega)} &= \delta^{-\frac{1}{2}}\Big\Vert \nabla \Tilde{\mathrm{H}}\Big\Vert_{\mathbb{L}^4(\mathrm{B})}
		\\ & = \mathcal{O}\Big(\delta^{-\frac{1}{2}}\Big).
	\end{align}
	Consequently, from the above estimate and due to the fact $|\mathrm{E}|^2 = |\nabla \mathrm{H}|^2$ we obtain
	\begin{align}
		\Big\Vert |\mathrm{E}|^2\Big\Vert_{\mathrm{L}^2(\Omega)} = \mathcal{O}\Big(\delta^{-1}\Big).
	\end{align}
	This completes the proof of Proposition \ref{prop2}.
	
	\section{Proof of Theorem \ref{1.2}}\label{dielectric}
	We recall the related model
	\begin{align}\label{eq1}
		\Delta \mathrm{E} + \omega^2\mu_\mathrm{m}\varepsilon \mathrm{E} = 0 
	\end{align}
	where the electric permittivity $\varepsilon = {\varepsilon_{p}}\rchi_{\Omega} +\varepsilon_\mathrm{m}\rchi_{\mathbb{R}^{2}\setminus\overline{\Omega}}$, $\varepsilon_\mathrm{m} = \varepsilon_\infty \varepsilon'_\mathrm{m}$ is the permittivity of the host medium and $\mu_\mathrm{m} = \mu_{\infty}\mu'_\mathrm{m}$ is the permeability of the host medium respectively, both are assumed to be constant, independent of the frequency $\omega$ of the incident wave. Here, $\varepsilon_\infty$ and $\mu_\infty$ as the permittivity and permeability of the vacuum. For the permittivity, $\varepsilon_\mathrm{p}(\omega)$ of the Dielectric nanoparticle, we will use the so-called Lorenz model for the permittivity, which can be described as follows
	\begin{align}
		\varepsilon_\mathrm{p}(\omega) = \varepsilon_\infty\Bigg[1+\dfrac{\omega_\mathrm{p}^2}{\omega_0^2-\omega^2-i\gamma\omega} \Bigg]
	\end{align}
	where where $\omega^2_\mathrm{p}$ is the electric plasma frequency, $\omega^2_0$ is the undamped resonance frequency and $\gamma$ is the electric damping parameter.
	\hfill \break
	Furthermore, let us recall that the volumetric Logarithmic Operator $\displaystyle \int_\Omega -\frac{1}{2\pi}\log|\mathrm{x}-\mathrm{y}|\mathrm{E}(\mathrm{y})d\mathrm{y}$ has eigenvalues $\lambda_\mathrm{n}^{(\boldsymbol{\ell})}$ and the corresponding eigenfunctions $\mathrm{e}_\mathrm{n}^{(\boldsymbol{\ell})}$. The  following properties are needed, namely, for a particle $\Omega$ of radius $\delta <<1$, $\displaystyle\int_\Omega\mathrm{e}_\mathrm{n}^{(\boldsymbol{\ell})}(\mathrm{x})d\mathrm{x} \neq 0$ and $\lambda_\mathrm{n}^{(\boldsymbol{\ell})}\sim \delta^2|\log\delta|.$ These properties are are justified for $n=1$, when $\Omega$ is the disc of radius $\delta$, see \cite{Ahcene}.
	\bigbreak
	\noindent
	The unique solution of the problem (\ref{eq1}) satisfies also the following Lippmann-Schwinger equation
	\begin{align}\label{LS}
		\mathrm{E}(\mathrm{x}) - \omega^2\mu_\mathrm{m}\Big(\varepsilon_\mathrm{p}(\omega,\gamma)-\varepsilon_\mathrm{m}\Big)\int_\Omega \mathbb{G}^{(\mathrm{k})}(\mathrm{x},\mathrm{y})\mathrm{E}(\mathrm{y})d\mathrm{y} = \mathrm{E}^{\textbf{in}}(\mathrm{x})
	\end{align}
	where we recall the fundamental solution, $\mathbb{G}^{(\mathrm{k})}$
	\begin{equation*}
		\mathbb{G}^{(k)}(\mathrm{x},\mathrm{y}) = \dfrac{i}{4}\mathcal{H}^{(1)}_{0}(k|\mathrm{x}-\mathrm{y}|), \quad \textbf{for} \quad \mathrm{x}\neq \mathrm{y}.
	\end{equation*}
	Moreover, we denote by $\bm{\tau_\mathrm{p}} := \varepsilon_\mathrm{p}(\omega,\gamma)-\varepsilon_\mathrm{m}$ and we also use the fact that $|1-\omega^2\mu_\mathrm{m}\varepsilon_\mathrm{p}\lambda_\mathrm{n_0}^{(\boldsymbol{\ell})}| \sim |\log\delta|^{-\mathrm{h}}$, see Section \ref{dielectric}. The  equation (\ref{LS}) can be rewritten as follows:
	\begin{align}\label{LS}
		\mathrm{E}(\mathrm{x}) - \omega^2\mu_\mathrm{m}\bm{\tau_\mathrm{p}}\int_\Omega \mathbb{G}^{(\mathrm{0})}(\mathrm{x},\mathrm{y})\mathrm{E}(\mathrm{y})d\mathrm{y} = \mathrm{E}^{\textbf{in}}(\mathrm{x}) + \omega^2\mu_\mathrm{m}\bm{\tau_\mathrm{p}}\int_\Omega \Big(\mathbb{G}^{(\mathrm{k})}-\mathbb{G}^{(\mathrm{0})}\Big)(\mathrm{x},\mathrm{y})\mathrm{E}(\mathrm{y})d\mathrm{y}.
	\end{align}
	Using the asymptotic expression $|\mathrm{x}-\mathrm{y}|\to 0$
	\begin{align}
		\mathbb{G}^{(k)}(\mathrm{x},\mathrm{y}) - \mathbb{G}^{(0)}(\mathrm{x},\mathrm{y}) = \bm{a} |\mathrm{x}-\mathrm{y}|^2\ln{|\mathrm{x}-\mathrm{y}|} + \bm{b} |\mathrm{x}-\mathrm{y}|^2 + \mathcal{O}(|\mathrm{x}-\mathrm{y}|^3) 
	\end{align}
	in the equation (\ref{LS}), rewritten in scaled domain $\mathrm{B}$ with $|\mathrm{B}| = \mathcal{O}(1)$ and $\Omega = \delta\mathrm{B} + \mathrm{z}$, we obtain
	\begin{align}\label{scale}
		\Tilde{\mathrm{E}}(\mathrm{\xi}) - \omega^2\mu_\mathrm{m}\bm{\tau_\mathrm{p}}\delta^2\int_\mathrm{B} \mathbb{G}^{(\mathrm{0})}(\mathrm{\xi},\mathrm{\eta})\Tilde{\mathrm{E}}(\mathrm{\eta})d\mathrm{\eta}\nonumber &= \Tilde{\mathrm{E}}^{\textbf{in}}(\mathrm{\xi}) - \omega^2\mu_\mathrm{m}\bm{\tau_\mathrm{p}}\delta^2\Bigg[\log\delta\int_\mathrm{B}\Tilde{\mathrm{E}}(\mathrm{\eta})d\mathrm{\eta} 
		+ \delta^2\log\delta\int_\mathrm{B}|\xi-\eta|^2\Tilde{\mathrm{E}}(\mathrm{\eta})d\mathrm{\eta}
		\\ \nonumber &+ \delta^2\int_\mathrm{B}|\xi-\eta|^2\ln|\xi-\eta|\Tilde{\mathrm{E}}(\mathrm{\eta})d\mathrm{\eta} + \delta^2\int_\mathrm{B}|\xi-\eta|^2\Tilde{\mathrm{E}}(\mathrm{\eta})d\mathrm{\eta} \\ &+ \delta^3 \int_\mathrm{B}|\xi-\eta|^3\Tilde{\mathrm{E}}(\mathrm{\eta})d\mathrm{\eta}\Bigg].    
	\end{align}
	Moreover, we use the local $\mathrm{H}^2$-regularity of the solution of (\ref{eq1}) in $\Omega$. The following inequality is a consequence of interpolation inequalities based on Gagliardo-Nirenberg, we refer to \cite[p. 313]{heim}.
	\begin{lemma}
		Let $\mathrm{B} \subset \mathbb{R}^\mathrm{N}$ be a regular bounded open set. Let $ 1\le \mathrm{q}\le \mathrm{p} < \infty$. Then the following inequality follows
		\begin{align}
			\Big\Vert u\Big\Vert_{\mathrm{L}^\mathrm{p}\Big(\mathrm{B}\Big)} \lesssim \Big\Vert u\Big\Vert_{\mathrm{L}^\mathrm{q}\Big(\mathrm{B}\Big)}^{1-\mathrm{r}}\Big\Vert u\Big\Vert_{W^{1,\mathrm{N}}\Big(\mathrm{B}\Big)}^\mathrm{r} \quad \forall u \in W^{1,\mathrm{N}}\Big(\mathrm{B}\Big), \quad \textbf{where} \quad \mathrm{r} = 1-\frac{\mathrm{q}}{\mathrm{p}}
		\end{align}
		In particular, for $\mathrm{N}=2, \mathrm{p}=4, \mathrm{q}=2$ and $\mathrm{r}=\frac{1}{2}$ we obtain
		\begin{align}
			\Big\Vert u\Big\Vert_{L^\mathrm{4}\Big(\mathrm{B}\Big)} \lesssim \Big\Vert u\Big\Vert_{L^\mathrm{2}\Big(\mathrm{B}\Big)}^{\frac{1}{2}}\Big\Vert u\Big\Vert_{\mathrm{H}^1\Big(\mathrm{B}\Big)}^\frac{1}{2}.
		\end{align}
	\end{lemma}
	Next, we use the continuity property of the integral $\displaystyle \int_\mathrm{B} \mathbb{G}^{(\mathrm{0})}(\mathrm{\xi},\mathrm{\eta})\Tilde{\mathrm{E}}(\mathrm{\eta})d\mathrm{\eta}$ from $\mathrm{L}^2\Big(\mathrm{B}\Big)\to \mathrm{H}^2\Big(\mathrm{B}\Big)$ and the smoothness of the remaining integral of the equation (\ref{scale}) to obtain
	\begin{align}\label{es1}
		\nonumber
		\Big\Vert \Tilde{\mathrm{E}} \Big\Vert_{\mathrm{H}^2\Big(B\Big)} &\lesssim \Big\Vert \Tilde{\mathrm{E}}^{\textbf{in}} \Big\Vert_{\mathrm{H}^2\Big(B\Big)} + \omega^2\mu_\mathrm{m}\bm{\tau_\mathrm{p}}\delta^2 \Bigg[ \Big\Vert \Tilde{\mathrm{E}} \Big\Vert_{\mathrm{L}^2\Big(\mathrm{B}\Big)} + |\ln\delta| \Big\Vert \Tilde{\mathrm{E}} \Big\Vert_{\mathrm{L}^2\Big(\mathrm{B}\Big)} + \delta^2 |\ln{\delta}|\Big\Vert  \Tilde{\mathrm{E}} \Big\Vert_{\mathrm{L}^2\Big(\mathrm{B}\Big)} +  \delta^2\Big\Vert \Tilde{\mathrm{E}} \Big\Vert_{\mathrm{L}^2\Big(\mathrm{B}\Big)}
		\\  &+\delta^3\Big\Vert \Tilde{\mathrm{E}} \Big\Vert_{\mathrm{L}^2\Big(\mathrm{B}\Big)}
		\Bigg].
	\end{align}
	We have the following estimate as derived in \cite{Ahcene},
	\begin{align}\label{p1}
		\Big\Vert \mathrm{E } \Big\Vert_{\mathrm{L}^2\Big(\Omega\Big)} = \mathcal{O}\Big(\delta|\log{\delta}|^\mathrm{h}\Big).
	\end{align}
	Therefore with the help of (\ref{p1}), we derive that 
	\begin{align}
		\Big\Vert \Tilde{\mathrm{E}} \Big\Vert_{\mathrm{H}^2\Big(B\Big)} = \mathcal{O}\Big(1\Big).
	\end{align}
	Hence by interpolation, we obtain 
	\begin{align}
		\Big\Vert |\mathrm{E}|^2 \Big\Vert_{\mathrm{L}^2\Big(\Omega\Big)} = \mathcal{O}\Big(\delta|\log\delta|^{\frac{3\mathrm{h}}{2}}\Big).
	\end{align}
	Next, we state the following result, obtained in \cite{Ahcene}, to derive the dominant term of $\displaystyle\int_\Omega|\mathrm{E}|^2(\mathrm{y})d\mathrm{y}$.
	\begin{lemma}
		We have the following approximation
		\begin{align}
			\displaystyle\int_\Omega|\mathrm{E}|^2(\mathrm{y})d\mathrm{y} = \dfrac{|\mathrm{E}^{\textbf{in}}|^2(\mathrm{z)}\Big(\displaystyle\int_\Omega \mathrm{e}_\mathrm{n_0}^{(\boldsymbol{\ell})}(\mathrm{y})d\mathrm{y}\Big)^2}{|1-\omega^2\mu_\mathrm{m}\varepsilon_\mathrm{p}\lambda_\mathrm{n_0}^{(\boldsymbol{\ell})}|^2} + \mathcal{O}\Big(\delta^2|\log\delta|^{3\mathrm{h}-1}\Big).
		\end{align}
	\end{lemma}
	\section{\textbf{Appendix}}\label{sec6}
	\subsection{\textbf{Derivation of the System of Boundary Integral Equations:}}
	As a first step, we rewrite the representation (\ref{rep formula}) as follows
	\begin{align}\label{1stIE1}
		\mathrm{U}_{\mathrm{i}}(\mathrm{x},t) &=   \mathcal{S}_{\alpha_\mathrm{p}}\Big[\gamma^{\textbf{int}}_{1}\mathrm{U}_{\mathrm{i}}\Big](\mathrm{x},t) - \mathcal{D}_{\alpha_\mathrm{p}}\Big[\gamma^{\textbf{int}}_{0}\mathrm{U}_{\mathrm{i}}\Big](\mathrm{x},t) + \frac{\omega \cdot \boldsymbol{\Im}(\varepsilon_\mathrm{p})}{2\pi\gamma_{\mathrm{p}}}\mathcal{V}\Big[ |\mathrm{E}|^{2} \Big](\mathrm{x},t) , \ \text{for} \ \mathrm{x} \in \Omega \times (0,\mathrm{T}).
	\end{align}
	By using the jump relation of the heat potential we can obtain the Neumann trace of the (\ref{1stIE1})
	\begin{align}
		\nonumber
		\gamma^{\textbf{int}}_{1}\mathrm{U}_{\mathrm{i}}(\mathrm{x},t) &=   \gamma^{\textbf{int}}_{1,\mathrm{x}}\mathcal{S}_{\alpha_\mathrm{p}}\Big[\gamma^{\textbf{int}}_{1}\mathrm{U}_{\mathrm{i}}\Big](\mathrm{x},t)- \gamma^{\textbf{int}}_{1,\mathrm{x}}\mathcal{D}_{\alpha_\mathrm{p}}\Big[\gamma^{\textbf{int}}_{0}\mathrm{U}_{\mathrm{i}}\Big](\mathrm{x},t) + \frac{\omega \cdot \boldsymbol{\Im}(\varepsilon_\mathrm{p})}{2\pi\gamma_{\mathrm{p}}}\gamma^{\textbf{int}}_{1,\mathrm{x}}\mathcal{V}\Big[ |\mathrm{E}|^{2} \Big](\mathrm{x},t) \\ \nonumber &=  \Big(\frac{1}{2}+\mathcal{K^*_{\alpha_\mathrm{p}}}\Big)\Big[\gamma^{\textbf{int}}_{1}\mathrm{U}_{\mathrm{i}}\Big](\mathrm{x},t) - \gamma^{\textbf{int}}_{1,\mathrm{x}}\mathcal{D}_{\alpha_\mathrm{p}}\Big[\gamma^{\textbf{int}}_{0}\mathrm{U}_{\mathrm{i}}\Big](\mathrm{x},t) + \frac{\omega \cdot \boldsymbol{\Im}(\varepsilon_\mathrm{p})}{2\pi\gamma_{\mathrm{p}}}\gamma^{\textbf{int}}_{1,\mathrm{x}}\mathcal{V}\Big[ |\mathrm{E}|^{2} \Big](\mathrm{x},t),
	\end{align}
	where $\mathcal{K^*_{\alpha_\mathrm{p}}}$ is the adjoint of double layer heat operator. Next, we denote $\mathcal{H}_{\alpha_\mathrm{p}}$ as the hypersingular boundary integral operator with density $\psi$
	\begin{align}\label{hypersingular}
		\mathbb{H}_{\alpha_\mathrm{p}}\Big[\psi\Big](\mathrm{x},\mathrm{t}) := -\gamma^{\textbf{int}}_{1,\mathrm{x}}\mathcal{D}_{\alpha_\mathrm{p}}\Big[\psi\Big](\mathrm{x},t) =-\gamma^{\textbf{int}}_{1,\mathrm{x}} \frac{1}{\alpha}\displaystyle\int_0^{t}\int_{\partial\Omega} \gamma^{\textbf{int}}_{1,\mathrm{y}}\Phi(\mathrm{x},t;\mathrm{y},\tau) \psi(y,\tau)d\sigma_\mathrm{y}d\tau.
	\end{align}
	Hence, we obtain
	\begin{align}
		\Bigg(\frac{1}{2}I_{d} - \mathcal{K}^{*}_{\alpha_\mathrm{p}}\Bigg)\Bigg[\gamma^{\textbf{int}}_{1}\mathrm{U}_{\mathrm{i}}\Bigg](\mathrm{x},t) = \mathbb{H}_{\alpha_\mathrm{p}}\Big[\gamma^{\textbf{int}}_{0}\mathrm{U}_{\mathrm{i}}\Big](\mathrm{x},t) + \frac{\omega \cdot \boldsymbol{\Im}(\varepsilon_\mathrm{p})}{2\pi\gamma_{\mathrm{p}}}\gamma^{\textbf{int}}_{1,\mathrm{x}}\mathcal{V}\Big[ |\mathrm{E}|^{2} \Big](\mathrm{x},t).
	\end{align}
	In a similar way, taking the dirichlet trace of the equation (\ref{1stIE1}) and using the jump relation of heat potential, we get
	\begin{align}
		\nonumber
		\gamma^{\textbf{int}}_{0,\mathrm{x}}\mathrm{U}_{\mathrm{i}}(\mathrm{x},t) &=   \gamma^{\textbf{int}}_{0,\mathrm{x}}\mathcal{S}_{\alpha_\mathrm{p}}\Big[\gamma^{\textbf{int}}_{1}\mathrm{U}_{\mathrm{i}}\Big](\mathrm{x},t) - \gamma^{\textbf{int}}_{0,\mathrm{x}}\mathcal{D}_{\alpha_\mathrm{p}}\Big[\gamma^{\textbf{int}}_{0,\mathrm{x}}\mathrm{U}_{\mathrm{i}}\Big](\mathrm{x},t)+ \frac{\omega \cdot \boldsymbol{\Im}(\varepsilon_\mathrm{p})}{2\pi\gamma_{\mathrm{p}}}\gamma^{\textbf{int}}_{0,\mathrm{x}}\mathcal{V}\Big[ |\mathrm{E}|^{2} \Big](\mathrm{x},t) \\ \nonumber &= \mathcal{S}_{\alpha_\mathrm{p}}\Big[\gamma^{\textbf{int}}_{1}\mathrm{U}_{\mathrm{i}}\Big](\mathrm{x},t)- \Big(-\frac{1}{2}+\mathcal{K}_{\alpha_\mathrm{p}}\Big)\Big[\gamma^{\textbf{int}}_{0}\mathrm{U}_{\mathrm{i}}\Big](\mathrm{x},t) + \frac{\omega \cdot \boldsymbol{\Im}(\varepsilon_\mathrm{p})}{2\pi\gamma_{\mathrm{p}}}\gamma^{\textbf{int}}_{0,\mathrm{x}}\mathcal{V}\Big[ |\mathrm{E}|^{2} \Big](\mathrm{x},t),
	\end{align}
	where $\mathcal{K}$ is the double layer heat operator. Hence, we get
	\begin{align}
		\Bigg(\frac{1}{2}I_{d} + \mathcal{K}_{\alpha_\mathrm{p}}\Bigg)\Bigg[\gamma^{\textbf{int}}_{0}\mathrm{U}_{\mathrm{i}}\Bigg](\mathrm{x},t) = \mathcal{S}_{\alpha_\mathrm{p}}\Big[\gamma^{\textbf{int}}_{1}\mathrm{U}_{\mathrm{i}}\Big](\mathrm{x},t)+ \frac{\omega \cdot \boldsymbol{\Im}(\varepsilon_\mathrm{p})}{2\pi\gamma_{\mathrm{p}}}\gamma^{\textbf{int}}_{0}\mathcal{V}\Big[ |\mathrm{E}|^{2} \Big](\mathrm{x},t).
	\end{align} 
	\noindent
	Thus, using the internal Dirichlet problem along with the initial condition, we arrive at the integral representations
	\begin{align}\label{2ndIE}
		\nonumber
		\Big(\frac{1}{2}I_{d} + \mathcal{K}_{\alpha_\mathrm{p}}\Big)\Big[\gamma^{\textbf{int}}_{0}\mathrm{U}_{\mathrm{i}}\Big](\mathrm{x},t) = \mathcal{S}_{\alpha_\mathrm{p}}\Big[\gamma^{\textbf{int}}_{1}\mathrm{U}_{\mathrm{i}}\Big](\mathrm{x},t) + \frac{\omega \cdot \boldsymbol{\Im}(\varepsilon_\mathrm{p})}{2\pi\gamma_{\mathrm{p}}}\gamma^{\textbf{int}}_{0}\mathcal{V}\Big[ |\mathrm{E}|^{2} \Big](\mathrm{x},t),
		\\
		\Big(\frac{1}{2}I_{d} - \mathcal{K}_{\alpha_\mathrm{p}}^{*}\Big)\Big[\gamma^{\textbf{int}}_{1}\mathrm{U}_{\mathrm{i}}\Big](\mathrm{x},t) =  \mathbb{H}_{\alpha_\mathrm{p}}\Big[\gamma^{\textbf{int}}_{0}\mathrm{U}_{\mathrm{i}}\Big](\mathrm{x},t) + \frac{\omega \cdot \boldsymbol{\Im}(\varepsilon_\mathrm{p})}{2\pi\gamma_{\mathrm{p}}}\gamma^{\textbf{int}}_{1,\mathrm{x}}\mathcal{V}\Big[ |\mathrm{E}|^{2} \Big](\mathrm{x},t).
	\end{align}
	In a similar manner, for exterior Dirichlet Problem we arrive to the integral representations 
	\begin{align}\label{3ndIE}
		\nonumber
		\Big(\frac{1}{2}I_{d} - \mathcal{K}_{\alpha_\mathrm{m}}\Big)\Big[\gamma^{\textbf{ext}}_{0}\mathrm{U}_{\mathrm{e}}\Big](\mathrm{x},t) = -\mathcal{S}_{\alpha_\mathrm{m}}\Big[\gamma^{\textbf{ext}}_{1}\mathrm{U}_{\mathrm{e}}\Big](\mathrm{x},t),\\
		\Big(\frac{1}{2}I_{d} + \mathcal{K}_{\alpha_\mathrm{m}}^{*}\Big)\Big[\gamma^{\textbf{ext}}_{1}\mathrm{U}_{\mathrm{e}}\Big](\mathrm{x},t) = -\mathcal{H}_{\alpha_\mathrm{m}}\Big[\gamma^{\textbf{ext}}_{0}\mathrm{U}_{\mathrm{e}}\Big](\mathrm{x},t),
	\end{align}
	where $\alpha_\mathrm{m} = \frac{\rho_{m}C_{m}}{\gamma_{m}}$. From the first boundary-integral equation of (\ref{3ndIE}), we can have
	\begin{equation}\label{a11}
		\gamma^{\textbf{ext}}_{1}\mathrm{U}_{\mathrm{e}} = -\mathcal{S}_{\alpha_\mathrm{m}}^{-1} \Big(\frac{1}{2}I_{d} - \mathcal{K}_{\alpha_\mathrm{m}}\Big)\Big[\gamma^{\textbf{ext}}_{0}\mathrm{U}_{\mathrm{e}}\Big].
	\end{equation}
	Also, from the second boundary-integral equation of (\ref{3ndIE}), we have
	\begin{equation}\label{a22}
		\gamma^{\textbf{ext}}_{1}\mathrm{U}_{\mathrm{e}} = -\mathbb{H}_{\alpha_\mathrm{m}}\Big[\gamma^{\textbf{ext}}_{0}\mathrm{U}_{\mathrm{e}}\Big](\mathrm{x},t) + \Big(\frac{1}{2}I_{d} - \mathcal{K}_{\alpha_\mathrm{m}}^{*}\Big)\Big[\gamma^{\textbf{ext}}_{1}\mathrm{U}_{\mathrm{e}}\Big](\mathrm{x},t).
	\end{equation}
	Now replacing (\ref{a11}) in (\ref{a22}), we get
	\begin{equation}\label{a3}
		\gamma^{\textbf{ext}}_{1}\mathrm{U}_{\mathrm{e}} = -\mathbb{H}_{\alpha_\mathrm{m}}\Big[\mathrm{U}_{\mathrm{e}}\Big](\mathrm{x},t) - \Big(\frac{1}{2}I_{d} - \mathcal{K}_{\alpha_\mathrm{m}}^{*}\Big)\mathcal{S}_{\alpha_\mathrm{m}}^{-1} \Big(\frac{1}{2}I_{d} - \mathcal{K}_{\alpha_\mathrm{m}}\Big)\Big[\gamma^{\textbf{ext}}_{0}\mathrm{U}_{\mathrm{e}}\Big].
	\end{equation}
	We denote by $\mathbb{A}^{\textbf{ext}} := \mathbb{H}_{\alpha_\mathrm{m}} + \Big(\frac{1}{2}I_{d} - \mathcal{K}_{\alpha_\mathrm{m}}^{*}\Big)\mathcal{S}_{\alpha_\mathrm{m}}^{-1} \Big(\frac{1}{2}I_{d} - \mathcal{K}_{\alpha_\mathrm{m}}\Big)$, which is known as Steklov-Poincare heat operator.
	Again from the first boundary-integral equation of (\ref{2ndIE}), we get 
	\begin{equation}\label{a1111}
		\gamma^{\textbf{int}}_{1}\mathrm{U}_{\mathrm{i}} = \mathcal{S}_{\alpha_\mathrm{p}}^{-1} \Big(\frac{1}{2}I_{d} + \mathcal{K}_{\alpha_\mathrm{p}}\Big)\Big[\gamma^{\textbf{int}}_{0}\mathrm{U}_{\mathrm{i}}\Big] - \frac{\omega \cdot \boldsymbol{\Im}(\varepsilon_\mathrm{p})}{2\pi\gamma_{\mathrm{p}}}\mathcal{S}_{\alpha_\mathrm{p}}^{-1} \gamma^{\textbf{int}}_{0}\mathcal{V}\Big[ |\mathrm{E}|^{2} \Big](\mathrm{x},t).
	\end{equation}
	Therefore, from the transmission condition of (\ref{eq:heat1}) and replacing the equations (\ref{a3}) and (\ref{a1111}) in that, we get
	\begin{align}\label{aa4}
		\nonumber\gamma_{c}\gamma^{\textbf{int}}_{1}\mathrm{U}_{\mathrm{i}} = \gamma_{m}\gamma^{\textbf{ext}}_{1}\mathrm{U}_{\mathrm{e}} \\
		\gamma_{c}\mathcal{S}_{\alpha_\mathrm{p}}^{-1} \Big(\frac{1}{2}I_{d} + \mathcal{K}_{\alpha_\mathrm{p}}\Big)\Big[\gamma^{\textbf{int}}_{0}\mathrm{U}_{\mathrm{i}}\Big] - \gamma_{\mathrm{c}}\frac{\omega \cdot \boldsymbol{\Im}(\varepsilon_\mathrm{p})}{2\pi\gamma_{\mathrm{p}}}\mathcal{S}_{\alpha_\mathrm{p}}^{-1} \gamma^{\textbf{int}}_{0}\mathcal{V}\Big[ |\mathrm{E}|^{2} \Big]= -\gamma_{m}\mathbb{A}^{\textbf{ext}}\Big[\gamma^{\textbf{ext}}_{0}\mathrm{U}_{\mathrm{e}}\Big] 
	\end{align}
	Now from the transmission condition $\gamma^{\textbf{int}}_{0}\mathrm{U}_{\mathrm{i}} = \gamma^{\textbf{ext}}_{0}\mathrm{U}_{\mathrm{e}}$, we obtain from the equation (\ref{aa4})
	\begin{align}\label{a4}
		\nonumber \gamma_{c}\mathcal{S}_{\alpha_\mathrm{p}}^{-1} \Big(\frac{1}{2}I_{d} + \mathcal{K}_{\alpha_\mathrm{p}}\Big)\Big[\gamma^{\textbf{int}}_{0}\mathrm{U}_{\mathrm{i}}\Big] = -\gamma_{m}\mathbb{A}^{\textbf{ext}}\Big[\gamma^{\textbf{int}}_{0}\mathrm{U}_{\mathrm{i}}\Big] +
		\frac{\omega \cdot \boldsymbol{\Im}(\varepsilon_\mathrm{p})}{2\pi\gamma_{\mathrm{p}}}\gamma_{c}\mathcal{S}_{\alpha}^{-1} \gamma^{\textbf{int}}_{0}\mathcal{V}\Big[ |\mathrm{E}|^{2} \Big], \\ \text{equivalently,} \ \  \
		\gamma^{\textbf{int}}_{0}\mathrm{U}_{\mathrm{i}} = -\frac{\gamma_{\mathrm{m}}}{\gamma_{\mathrm{p}}} \Big(\frac{1}{2}I_{d} + \mathcal{K}_{\alpha_\mathrm{p}}\Big)^{-1}\mathcal{S}_{\alpha_\mathrm{p}} \mathbb{A}^{\textbf{ext}}\Big[\gamma^{\textbf{int}}_{0}\mathrm{U}_{\mathrm{i}}\Big] + \frac{\omega \cdot \boldsymbol{\Im}(\varepsilon_\mathrm{p})}{2\pi\gamma_{\mathrm{p}}}\Big(\frac{1}{2}I_{d} + \mathcal{K}_{\alpha_\mathrm{p}}\Big)^{-1}\gamma^{\textbf{int}}_{0}\mathcal{V}\Big[ |\mathrm{E}|^{2} \Big].
	\end{align}
	\subsection{Proof of Lemma \ref{4.1}}
	We start with recalling the definition of the fundamental solution for heat equation
	\begin{equation}\label{defs}
		\Phi(\mathrm{y},\tau; \mathrm{v},\mathrm{s}):=\  \begin{cases}
			\dfrac{\alpha}{4\pi(\tau-\mathrm{s})}\textbf{exp}\Big(-\dfrac{\alpha|\mathrm{y}-\mathrm{v}|^2}{4(\tau-\mathrm{s})}\Big) , \quad \tau > \mathrm{s} \\
			0 ,\quad \text{otherwise}
		\end{cases},
	\end{equation}
	where we set $\alpha := \alpha_\mathrm{p} = \frac{\rho_{\mathrm{p}}C_{\mathrm{p}}}{\gamma_{\mathrm{p}}}$, and the corresponding double layer heat operator $\mathcal{K}$
	\begin{align}\label{maineq}
		\mathcal{K}[f](\mathrm{y},\tau) := \frac{1}{\alpha}\displaystyle \int_{0}^{\tau}\int_{\partial\Omega}\partial_{\nu_\mathrm{y}}\Phi(\mathrm{y},\tau; \mathrm{v},\mathrm{s}) f(\mathrm{v},\mathrm{s}) d\sigma_{\mathrm{v}}d\mathrm{s}.
	\end{align}
	We see that $
	\partial_{\nu_\mathrm{y}}\Phi(\mathrm{y},\tau; \mathrm{v},\mathrm{s}) := -\frac{\alpha (\mathrm{v}-\mathrm{y})\cdot\nu_\mathrm{y}}{2(\tau-\mathrm{s})}\Phi(\mathrm{y},\tau; \mathrm{v},\mathrm{s}).$
	Furthermore, we recall the fundamental solution corresponding to the Laplace equation $\mathbb{G}^{(0)}(\mathrm{y},\mathrm{v}) := -\frac{1}{2\pi}\ln(|\mathrm{y}-\mathrm{v}|)$ and then $
	\partial_{\nu_\mathrm{y}}\mathbb{G}^{(0)}(\mathrm{y},\mathrm{v}) := -\frac{(\mathrm{v}-\mathrm{y})\cdot \nu_\mathrm{y}}{2\pi}\frac{1}{|\mathrm{y}-\mathrm{v}|^2}.
	$
	Hence, we derive the following representation from  (\ref{maineq}) applied to $f:=f_{\xi,\mathrm{t}, \mathrm{z}}( \mathrm{y}, \mathrm{\tau})=\Phi^{\textbf{e}}(\xi,\mathrm{t};\mathrm{z},\cdot)$, which is constant with respect to the variable $\mathrm{y}$,
	\begin{align}
		\nonumber
		\mathcal{K}\Big[\Phi^{\textbf{e}}(\xi,\mathrm{t};\mathrm{z},\cdot)\Big](\mathrm{y},\tau) &= \frac{1}{\alpha}\displaystyle \int_{0}^{\tau}\int_{\partial\Omega} -\frac{\alpha (\mathrm{v}-\mathrm{y})\cdot\nu_\mathrm{y}}{2(\tau-\mathrm{s})}\Phi(\mathrm{y},\tau; \mathrm{v},\mathrm{s})\Phi^{\textbf{e}}(\xi,\mathrm{t};\mathrm{z},\mathrm{s})d\sigma_{\mathrm{v}}d\mathrm{s} \\ \nonumber
		&=\frac{1}{\alpha}\displaystyle \int_{0}^{t}\int_{\partial\Omega}-\frac{\alpha (\mathrm{v}-\mathrm{y})\cdot\nu_\mathrm{y}}{2(\tau-\mathrm{s})} \dfrac{\alpha}{4\pi(\tau-\mathrm{s})}\textbf{exp}\Big(-\dfrac{\alpha|\mathrm{y}-\mathrm{v}|^2}{4(\tau-\mathrm{s})}\Big)\Phi^{\textbf{e}}(\xi,\mathrm{t};\mathrm{z},\mathrm{s})d\sigma_{\mathrm{v}}d\mathrm{s} \\ &=\displaystyle\int_{\partial\Omega} -\dfrac{(\mathrm{v}-\mathrm{y})\cdot \nu_\mathrm{y}}{2\pi}\frac{1}{|\mathrm{y}-\mathrm{v}|^2}\Bigg[\displaystyle\int_{0}^{\tau}\dfrac{\alpha |\mathrm{y}-\mathrm{v}|^2}{4(\tau-\mathrm{s})^2}\textbf{exp}\Big(-\dfrac{\alpha|\mathrm{y}-\mathrm{v}|^2}{4(\tau-\mathrm{s})}\Big)\Phi^{\textbf{e}}(\xi,\mathrm{t};\mathrm{z},\mathrm{s})d\mathrm{s}\Bigg]d\sigma_{\mathrm{v}}
	\end{align}
	Now, let us define
	\begin{equation}\label{defvarphi}
		\varphi(\mathrm{v},\mathrm{y},\mathrm{t}, \tau) := \displaystyle\int_{0}^{\tau}\dfrac{\alpha |\mathrm{y}-\mathrm{v}| ^2}{4(\mathrm{s}-\tau)^2}\textbf{exp}\Big(-\dfrac{\alpha|\mathrm{y}-\mathrm{v}| ^2}{4(\mathrm{s}-\tau)}\Big)\Phi^{\textbf{e}}(\xi,\mathrm{t};\mathrm{z},\mathrm{s})d\mathrm{s}.
	\end{equation}
	We use the change of variable $\mathrm{m} :=\frac{\sqrt{\alpha}|\mathrm{y}-\mathrm{v}|}{2\sqrt{\tau-\mathrm{s}}}$, then it follows that 
	$
	\mathrm{s} = \tau - \frac{\alpha|\mathrm{y}-\mathrm{v}|^2}{4\mathrm{m}^2}$ and then $d\mathrm{s} = \frac{1}{2}\alpha|\mathrm{y}-\mathrm{v}|^2\mathrm{m}^{-3}.
	$
	Consequently, we obtain from (\ref{defvarphi}) 
	\begin{align}\label{finalvarphi}
		\nonumber
		\varphi(\mathrm{v},\mathrm{y},\mathrm{t}, \tau) &= \displaystyle\int_{\frac{\sqrt{\alpha}|\mathrm{y}-\mathrm{v}|}{2\sqrt{\tau}}}^{\infty}\textbf{exp}(-\mathrm{m}^2)\Phi^{\textbf{e}}\Big(\xi,\mathrm{t}; \mathrm{z},\tau - \dfrac{\alpha|\mathrm{y}-\mathrm{v}|^2}{4\mathrm{m}^2}\Big)\frac{1}{(\tau-\mathrm{s})}\mathrm{m}^2\frac{1}{2}\alpha|\mathrm{y}-\mathrm{v}|^2\mathrm{m}^{-3}d\mathrm{m} \\ \nonumber&=\displaystyle\int_{\frac{\sqrt{\alpha}|\mathrm{y}-\mathrm{v}|}{2\sqrt{\tau}}}^{\infty}\textbf{exp}(-\mathrm{m}^2)\Phi^{\textbf{e}}\Big(\xi,\mathrm{t}; \mathrm{z},\tau - \dfrac{\alpha|\mathrm{y}-\mathrm{v}|^2}{4\mathrm{m}^2}\Big)\frac{\alpha|\mathrm{y}-\mathrm{v}|^2}{2(\tau-\mathrm{s})}\mathrm{m}^{-1}d\mathrm{m} \\ \nonumber&= \displaystyle\int_{\frac{\sqrt{\alpha}|\mathrm{y}-\mathrm{v}|}{2\sqrt{\tau}}}^{\infty}\textbf{exp}(-\mathrm{m}^2)\Phi^{\textbf{e}}\Big(\xi,\mathrm{t}; \mathrm{z},\tau - \dfrac{\alpha|\mathrm{y}-\mathrm{v}|^2}{4\mathrm{m}^2}\Big)2\mathrm{m}^2\mathrm{m}^{-1}d\mathrm{m} \\
		&= \displaystyle\int_{\frac{\alpha|\mathrm{y}-\mathrm{v}|^2}{4\tau}}^{\infty}\textbf{exp}(-\mathrm{m}^2)\Phi^{\textbf{e}}\Big(\xi,\mathrm{t}; \mathrm{z},\tau - \dfrac{\alpha|\mathrm{y}-\mathrm{v}|^2}{4\mathrm{m}^2}\Big)d(\mathrm{m}^2).
	\end{align}
	We rewrite (\ref{finalvarphi}) as follows
	\begin{align}\label{finalvarphi1}
		\nonumber
		\varphi(\mathrm{v},\mathrm{y},\mathrm{t}, \tau) &= \displaystyle\int_{\frac{\alpha|\mathrm{y}-\mathrm{v}|^2}{4\tau}}^{\infty}\textbf{exp}(-\mathrm{m}^2)\Bigg[\Phi^{\textbf{e}}\Big(\xi,\mathrm{t}; \mathrm{z},\tau - \dfrac{\alpha|\mathrm{y}-\mathrm{v}|^2}{4\mathrm{m}^2}\Big) - \Phi^{\textbf{e}}\Big(\xi,\mathrm{t};\mathrm{z}, \tau \Big)\Bigg]d(\mathrm{m}^2) \\ 
		&+ \Phi^{\textbf{e}}\Big(\xi,\mathrm{t};\mathrm{z}, \tau \Big)\displaystyle\int_{\frac{\alpha|\mathrm{y}-\mathrm{v}|^2}{4\tau}}^{\infty}\textbf{exp}(-\mathrm{m}^2) d(\mathrm{m}^2).
	\end{align}
	From the identity $\displaystyle\int_{0}^\infty \textbf{exp}(-\mathrm{m}^2)d(\mathrm{m}^2) = 1$, we derive the following
	\begin{align}\label{finalvarphi1}
		\varphi(\mathrm{v},\mathrm{y},\mathrm{t}, \tau) - \Phi^{\textbf{e}}(\xi,\mathrm{t};\mathrm{z}, \tau ) \nonumber &= \displaystyle\int_{\frac{\alpha|\mathrm{y}-\mathrm{v}|^2}{4\tau}}^{\infty}\textbf{exp}(-\mathrm{m}^2)\Bigg[\Phi^{\textbf{e}}\Big(\xi,\mathrm{t}; \mathrm{z},\tau - \dfrac{\alpha|\mathrm{y}-\mathrm{v}|^2}{4\mathrm{m}^2}\Big) - \Phi^{\textbf{e}}\Big(\xi,\mathrm{t};\mathrm{z}, \tau \Big)\Bigg]d(\mathrm{m}^2)
		\nonumber
		\\ &- \Phi^{\textbf{e}}(\xi,\mathrm{t};\mathrm{z}, \tau )\displaystyle\int_{0}^{\frac{\alpha|\mathrm{y}-\mathrm{v}|^2}{4\tau}}\textbf{exp}(-\mathrm{m}^2) d(\mathrm{m}^2).
	\end{align}
	Next, we do the following computation
	\begin{align}\label{esti1}
		\Phi^{\textbf{e}}(\xi,\mathrm{t};\mathrm{z}, \tau ) = \displaystyle\int_{\mathrm{t}}^\tau \partial_{s}\Phi^{\textbf{e}}(\xi,\mathrm{t};\mathrm{z}, \mathrm{s} )d\mathrm{s} \nonumber &\lesssim \Vert 1\Vert_{\mathrm{H}^{\frac{1}{4}}(\mathrm{t},\tau)}\Vert \partial_{s}\Phi^{\textbf{e}}(\xi,\mathrm{t};\mathrm{z}, \cdot )\Vert_{\mathrm{H}^{-\frac{1}{4}}(\mathrm{t},\tau)}
		\\ &\lesssim \tau^{\frac{1}{2}} \Vert \partial_{s}\Phi^{\textbf{e}}(\xi,\mathrm{t};\mathrm{z}, \cdot )\Vert_{\mathrm{H}^{-\frac{1}{4}}(0,\tau)}.
	\end{align}
	In a similar way as before, we write the following expression
	\begin{equation}
		\Phi^{\textbf{e}}\Big(\xi,\mathrm{t}; \mathrm{z},\tau - \dfrac{\alpha|\mathrm{y}-\mathrm{v}|^2}{4\mathrm{m}^2}\Big) - \Phi^{\textbf{e}}\Big(\xi,\mathrm{t};\mathrm{z}, \tau \Big) = \displaystyle\int_{\tau}^{\tau - \frac{\alpha|\mathrm{y}-\mathrm{v}|^2}{4\mathrm{m}^2}}\partial_{s}\Phi^{\textbf{e}}\Big(\xi,\mathrm{t};\mathrm{z}, \mathrm{s}\Big)d\mathrm{s}.
	\end{equation}
	\noindent
	We observe that $\delta\ge \frac{\sqrt{\alpha}|\mathrm{y}-\mathrm{v}|}{2\sqrt{\tau}}$ which implies $\tau-\frac{\alpha|\mathrm{y}-\mathrm{v}|^2}{4\mathrm{m}^2} \ge 0$.
	Hence we get 
	\begin{align}\label{esti2}
		\Phi^{\textbf{e}}\Big(\xi,\mathrm{t}; \mathrm{z},\tau - \dfrac{\alpha|\mathrm{y}-\mathrm{v}|^2}{4\mathrm{m}^2}\Big) - \Phi^{\textbf{e}}\Big(\xi,\mathrm{t};\mathrm{z}, \tau \Big) = \mathcal{O} \Big(\frac{\sqrt{\alpha}|\mathrm{y}-\mathrm{v}|}{2\mathrm{m}}\Vert \partial_{s}\Phi^{\textbf{e}}(\xi,\mathrm{t};\mathrm{z}, \cdot )\Vert_{\mathrm{H}^{-\frac{1}{4}}(0,\tau)}\Big).
	\end{align}
	Therefore, we plug (\ref{esti1}) and (\ref{esti2}) in (\ref{finalvarphi1}) to obtain
	\begin{align}\label{varphi2}
		\nonumber
		\varphi(\mathrm{v},\mathrm{y},\mathrm{t}, \tau) - \Phi^{\textbf{e}}(\xi,\mathrm{t};\mathrm{z}, \tau ) &= \mathcal{O}\Bigg(\displaystyle \int_{\frac{\alpha|\mathrm{y}-\mathrm{v}|^2}{4\tau}}^{\infty}\frac{\textbf{exp}(-\mathrm{m}^2)}{\mathrm{m}}d(\mathrm{m}^2)\frac{\sqrt{\alpha}|\mathrm{y}-\mathrm{v}|}{2}\Vert \partial_{s}\Phi^{\textbf{e}}(\xi,\mathrm{t};\mathrm{z}, \cdot )\Vert_{\mathrm{H}^{-\frac{1}{4}}(0,\tau)}\Bigg)\\ &+ \mathcal{O} \Bigg(\displaystyle\int_{0}^{\frac{\alpha|\mathrm{y}-\mathrm{v}|^2}{4\tau}}\textbf{exp}(-\mathrm{m}^2) d(\mathrm{m}^2)\tau^{\frac{1}{2}}\Vert \partial_{s}\Phi^{\textbf{e}}(\xi,\mathrm{t};\mathrm{z}, \cdot )\Vert_{\mathrm{H}^{-\frac{1}{4}}(0,\tau)}\Bigg).
	\end{align}
	Now, we look into the integral $\displaystyle \int_{\frac{\alpha|\mathrm{y}-\mathrm{v}|^2}{4\tau}}^{\infty}\frac{\textbf{exp}(-\mathrm{m}^2)}{\mathrm{m}}d(\mathrm{m}^2)$ and substitute $\omega = \mathrm{m}^2$, then we deduce that 
	\begin{equation*}
		\displaystyle \int_{\frac{\alpha|\mathrm{y}-\mathrm{v}|^2}{4t}}^{\infty}\omega^{\frac{1}{2}-1}\textbf{exp}(-\omega)d\omega = \Gamma \Big(\frac{1}{2},\frac{\alpha|\mathrm{y}-\mathrm{v}|^2}{4\tau}\Big).    
	\end{equation*}
	Moreover, it is well known that $\Gamma\Big(\frac{1}{2},\frac{\alpha|\mathrm{y}-\mathrm{v}|^2}{4\tau}\Big) = \sqrt{\pi} \ \textbf{erfc}\Big(\frac{\sqrt{\alpha}|\mathrm{y}-\mathrm{v}|}{2\sqrt{\tau}}\Big)$, where $\textbf{erfc}(\mathrm{z})$ is the complementary error function and $\textbf{erfc}(\mathrm{z}) = 1- \textbf{erf}(\mathrm{z}).$ The Taylor series expansion 
	$\textbf{erf}(\mathrm{z})= {\sqrt{\pi}}\Big(\mathrm{z}-\frac{\mathrm{z}^3}{3}+\frac{\mathrm{z}^5}{10}-...\Big)$, which holds for every $\mathrm{z} \in \mathbb{C}$, implies that
	\begin{align}\label{tay1}
		\Gamma\Big(\frac{1}{2},\frac{\alpha|\mathrm{y}-\mathrm{v}|^2}{4\tau}\Big) = \sqrt{\pi} \ \textbf{erfc}\Big(\frac{\sqrt{\alpha}|\mathrm{y}-\mathrm{v}|}{2\sqrt{\tau}}\Big) = \sqrt{\pi}\Bigg(1-\frac{2}{\sqrt{\pi}}\Big(\mathrm{\frac{\sqrt{\alpha}|\mathrm{y}-\mathrm{v}|}{2\sqrt{\tau}}}-\frac{\mathrm{(\frac{\sqrt{\alpha}|\mathrm{y}-\mathrm{v}|}{2\sqrt{\tau}}})^3}{3}+\frac{\mathrm{(\frac{\sqrt{\alpha}|\mathrm{y}-\mathrm{v}|}{2\sqrt{\tau}}})^5}{10}-...\Big)\Bigg).
	\end{align}
	We also note that 
	\begin{align}\label{tay2}
		\displaystyle\int_{0}^{\frac{\alpha|\mathrm{y}-\mathrm{v}|^2}{4\tau}}\textbf{exp}(-\mathrm{m}^2) d(\mathrm{m}^2) = 1- \textbf{exp}\Big(\frac{\alpha|\mathrm{y}-\mathrm{v}|^2}{4\tau}\Big) = 1- \Big[1-\frac{\alpha|\mathrm{y}-\mathrm{v}|^2}{4\tau} + \frac{1}{2!}\Big(\frac{\alpha|\mathrm{y}-\mathrm{v}|^2}{4\tau}\Big)^2-... \Big]. 
	\end{align}
	As a next step, we plug (\ref{tay1}) and (\ref{tay2}) into (\ref{varphi2}) to derive
	\begin{align}
		\nonumber
		\varphi(\mathrm{v},\mathrm{y},\mathrm{t}, \tau) - \Phi^{\textbf{e}}(\xi,\mathrm{t};\mathrm{z}, \tau ) &= \mathcal{O}\Bigg(\sqrt{\pi}\frac{\sqrt{\alpha}|\mathrm{y}-\mathrm{v}|}{2}\Vert \partial_{s}\Phi^{\textbf{e}}(\xi,\mathrm{t};\mathrm{z}, \cdot )\Vert_{\mathrm{H}^{-\frac{1}{4}}(0,\tau)}\Bigg)\\ \nonumber&+ \mathcal{O}\Bigg(\frac{\sqrt{\alpha}|\mathrm{y}-\mathrm{v}|}{2}\frac{\sqrt{\alpha}|\mathrm{y}-\mathrm{v}|}{\sqrt{\tau}}\Vert \partial_{s}\Phi^{\textbf{e}}(\xi,\mathrm{t};\mathrm{z}, \cdot )\Vert_{\mathrm{H}^{-\frac{1}{4}}(0,\tau)}\Bigg)\\&+ \mathcal{O} \Bigg(\frac{\alpha|\mathrm{y}-\mathrm{v}|^2}{4\tau}\tau^{\frac{1}{2}}\Vert \partial_{s}\Phi^{\textbf{e}}(\xi,\mathrm{t};\mathrm{z}, \cdot )\Vert_{\mathrm{H}^{-\frac{1}{4}}(0,\tau)}\Bigg),  
	\end{align}
	and then we obtain
	\begin{align}
		\nonumber
		\varphi(\mathrm{v},\mathrm{y},\mathrm{t}, \tau) - \Phi^{\textbf{e}}(\xi,\mathrm{t};\mathrm{z}, \tau ) &= \mathcal{O}\Bigg(\sqrt{\pi}\frac{\sqrt{\alpha}|\mathrm{y}-\mathrm{v}|}{2}\Vert \partial_{s}\Phi^{\textbf{e}}(\xi,\mathrm{t};\mathrm{z}, \cdot )\Vert_{\mathrm{H}^{-\frac{1}{4}}(0,\tau)}\Bigg)\\ \nonumber&+ \mathcal{O}\Bigg(\frac{\sqrt{\alpha}|\mathrm{y}-\mathrm{v}|}{2\sqrt{\tau}}\sqrt{\alpha}|\mathrm{y}-\mathrm{v}|\Vert \partial_{s}\Phi^{\textbf{e}}(\xi,\mathrm{t};\mathrm{z}, \cdot )\Vert_{\mathrm{H}^{-\frac{1}{4}}(0,\tau)}\Bigg)\\&+ \mathcal{O}\Bigg(\frac{\sqrt{\alpha}|\mathrm{y}-\mathrm{v}|}{2\sqrt{\tau}}\frac{\sqrt{\alpha}|\mathrm{y}-\mathrm{v}|}{2}\Vert \partial_{s}\Phi^{\textbf{e}}(\xi,\mathrm{t};\mathrm{z}, \cdot )\Vert_{\mathrm{H}^{-\frac{1}{4}}(0,\tau)}\Bigg) .  
	\end{align}
	Consequently, because of the boundness of the term $\frac{\sqrt{\alpha}|\mathrm{y}-\mathrm{v}|}{2\sqrt{\tau}}$ we deduce that
	\begin{align}
		\varphi(\mathrm{v},\mathrm{y},\mathrm{t}, \tau) - \Phi^{\textbf{e}}(\xi,\mathrm{t};\mathrm{z}, \tau ) = \mathcal{O}\Bigg(\sqrt{\alpha}|\mathrm{y}-\mathrm{v}| \ \Vert \partial_{s}\Phi^{\textbf{e}}(\xi,\mathrm{t};\mathrm{z}, \cdot )\Vert_{\mathrm{H}^{-\frac{1}{4}}(0,\tau)}\Bigg).
	\end{align}
	This completes the proof.
	
	\subsection{Proof of Lemma \ref{5.1}}
	
	The goal is to show that the operator defined thought the expression
	\begin{equation}\label{T}
		\mathbb{T}\Big[ \mathrm{H}\Big](\mathrm{x}) = \int_{\Omega} \nabla \nabla\cdot \Big(\mathbb{G}^{(\mathrm{k})}-\mathbb{G}^{(0)}\Big)(\mathrm{x},\mathrm{y})\mathrm{H}(\mathrm{y})d\mathrm{y}.
	\end{equation}
	is a bounded operator from $\mathbb{L}^{2}(\Omega) \rightarrow \mathbb{H}^{2}(\Omega).$
	
	\noindent From the definition of Hankel function, we have the following asymptotic expansion
	\begin{equation}\label{expansion2}
		\partial_{\mathrm{ij}}\Big(\mathbb{G}^{(\mathrm{k})}-\mathbb{G}^{(0)}\Big)(\mathrm{x},\mathrm{y}) = 2a \delta_{\mathrm{ij}}\log|x-y| + 2\mathrm{a} \dfrac{(\mathrm{x}-\mathrm{y})_{\mathrm{i}}(\mathrm{x}-\mathrm{y})_{\mathrm{j}}}{|\mathrm{x}-\mathrm{y}|^{2}} + 2\mathrm{b} \delta_{\mathrm{ij}} + \mathcal{O}(|\mathrm{x}-\mathrm{y}|),\quad \mathrm{x}\sim \mathrm{y},
	\end{equation}
	where $\delta_{\mathrm{ij}}$ is the Kronecker delta. Here a and b are constants.
	\bigbreak
	\noindent
	We see that the leading term of the expansion (\ref{expansion2}), when injected into (\ref{T}), is the kernel of the volume potential operator $\mathcal{V}$ which is bounded from $\mathbb{L}^2(\Omega)$ to $ \mathbb{H}^{2}(\Omega)$. The other terms of the expansion (\ref{expansion2}) define smoother kernels. Therefore they define smoother operators as compared to $\mathcal{V}$. This implies that the operator $\mathbb{T}$ is bounded $\mathbb{L}^2(\Omega)$ to $\mathbb{H}^{2}(\Omega)$.
	
	\subsection{Modelling the Electric Permittivity for Lorentzian Particles}\label{Lorentz}
	The Lorentz model (also known as Drude-Lorentz oscillator model) describes an electron as damped harmonic oscillator driven by a strong force. 
	The purpose of this model is to determine, using Newton's Second Law, the motion of an electron of the particle (and its polarization as it is proportional to it) when it comes in contact with the incidence electric field. The polarization moment satisfies the following equation
	\begin{align}\label{Lorentz-model-equation}
		\frac{\mathrm{d}^2P}{\mathrm{d}\mathrm{t}^2} + \gamma \frac{\mathrm{d}P}{\mathrm{d}\mathrm{t}} + \omega_0^2P =  \varepsilon_\infty \omega^2_p u^\mathrm{i}
	\end{align}
	where $\omega_0^2$ is the undamped resonance frequency, $\gamma$ is the electric damping parameter, $\omega_p$ the electric plasmonic frequency and $\varepsilon_\infty$ is the permittivity of the vacuum. The above equation provides a basis for determining the polarization and the dielectric susceptibility (i.e. the equivalent dielectric constant). Indeed, consider the driving force to be an oscillating incidence electric field, which is given by the time-harmonic acoustic plane wave $ u^{\mathrm{i}}(\mathrm{x},\mathrm{t}):= \mathrm{e}^{i(\mathrm{k}\mathrm{x}\cdot \mathrm{d}-\omega \mathrm{t})}$
	where, $\mathrm{k} = \frac{\omega}{\mathrm{c_0}}$ is the wave number, $\omega$ the frequency of incidence, $\mathrm{c_0}$ the speed of light, in the vacuum, and $\mathrm{d}$ the direction of incidence. We look for the polarization of the form $P(x, t):=P(x, \omega) e^{-i \omega t}$, then we obtain its expression, from (\ref{Lorentz-model-equation}), of the form $
	P(x, t)=\frac{\epsilon_\infty \omega_p}{-\omega^2-i\gamma\omega+\omega_0^2}u^{\mathrm{i}}(\mathrm{x},\mathrm{t}).$
	On the other hand, we know that the polarization is related to the incident field through the electric susceptibility $\chi_e$ as $P(x, t)=\epsilon_\infty \chi_e u^{\mathrm{i}}(\mathrm{x},\mathrm{t})$. As the electric susceptibility is defined through the relation $\frac{\epsilon_p}{\epsilon_\infty}=1+\chi_e$, we derive the following expression of the permittivity $\varepsilon_\mathrm{p}(\omega)$  
	\begin{align}\label{Lorentz-model-elect-permi-Appendix}
		\varepsilon_\mathrm{p}(\omega) = \varepsilon_\infty\Bigg[1+\dfrac{\omega_\mathrm{p}^2}{\omega_0^2-\omega^2-i\gamma\omega} \Bigg].   
	\end{align}
	For more information and details about the model (\ref{Lorentz-model-elect-permi-Appendix}), we refer to \cite{A-B-B, lmodel2, lmodel1}.
	\vskip 0.1in
	\noindent
	Based on this model, we discuss the regimes, in terms of the incident frequency $\omega$, where the nanoparticle behaves as  a plasmonic nanoparticle or a as dielectric nanoparticle.
	\subsubsection{Modeling with Plasmonic Nanoparticles}\label{plasmonics}
	The purpose of this section is to use the Lorentz model to generate the condition on the permittivity used in Section \ref{sec4}, namely $|1-\alpha\lambda_{\mathrm{n_0}}^{(3)}| = \mathcal{O}(\delta^\mathrm{h}),$ where $\alpha = \frac{1}{\varepsilon_\mathrm{p}(\omega)}-\frac{1}{\varepsilon_\mathrm{m}}$ and $\mathrm{h}<1$. Here, $\lambda^{(3)}_\mathrm{n_0}$'s are the eigenvalues of the Magnetization operator $\mathbb{M}: \nabla \mathrm{H}_{\text{arm}}\rightarrow \nabla \mathrm{H}_{\text{arm}}$ and $\varepsilon_\mathrm{p}(\omega), \varepsilon_\mathrm{m}$ are the corresponding permittivity of the particle and host medium respectively. We state the following lemma.
	\begin{lemma}
		If we consider the electric damping parameter $\gamma \sim \delta^\mathrm{h}$ and the frequency $\omega$ of the incidence wave as $\omega^2 = \omega_0^2 + \omega_\mathrm{p}^2\dfrac{\varepsilon_\mathrm{m}+\lambda_\mathrm{n_0}^{(3)}}{\lambda_\mathrm{n_0}^{(3)}(1-\varepsilon_\infty^{-1}\varepsilon_\mathrm{m})+\varepsilon_\mathrm{m}} + O(\delta^\mathrm{h})$, we have the required equality $|1-\alpha\lambda_{\mathrm{n_0}}^{(3)}| = \mathcal{O}(\delta^\mathrm{h}).$
	\end{lemma}
	\noindent
	Before justifying this lemma, let us mention that under the conditions on $\gamma$ and $\omega$, the permittivity of the nanoparticle has a negative real part (as $\omega^2>\omega_0^2$). This negativity property is the one that is usually used in the literature to design plasmonic nanoparticles. In our analysis, the choice of the incident frequency of the form presented in the previous lemma is key.
	\bigbreak
	\noindent
	\textbf{Proof.} From the definition of $\alpha$ and the fact that both $\varepsilon_\mathrm{p}$ and $\varepsilon_\mathrm{m}$ are moderate in terms of the parameter $\delta$, then we have
	\begin{align}
		\nonumber
		1-\alpha\lambda_\mathrm{n_0}^{(3)} = \mathcal{O}(\delta^\mathrm{h}) 
		\nonumber
		\iff \varepsilon_\mathrm{p} = \lambda_\mathrm{n_0}^{(3)} - \dfrac{\Big(\lambda_\mathrm{n_0}^{(3)}\Big)^2}{\varepsilon_\mathrm{m}+ \lambda_\mathrm{n_0}^{(3)}} + \mathcal{O}(\delta^\mathrm{h})
	\end{align}
	Then substituting the value of $\varepsilon_\mathrm{p}$, derived from the Lorentz model, we deduce
	\begin{align}
		1-\alpha\lambda_\mathrm{n_0}^{(3)} = \mathcal{O}(\delta^\mathrm{h}) \iff \dfrac{\omega_\mathrm{p}^2}{\omega_0^2-\omega^2-i\gamma\omega} = -1 + \varepsilon_\infty^{-1}\lambda_\mathrm{n_0}^{(3)} - \dfrac{\varepsilon_\infty^{-1}(\lambda_\mathrm{n_0}^{(3)})^2}{\varepsilon_\mathrm{m}+ \lambda_\mathrm{n_0}^{(3)}} + \mathcal{O}(\delta^\mathrm{h}).
	\end{align}
	Now, let us denote $\beta:= -1 + \varepsilon_\infty^{-1}\lambda_\mathrm{n_0}^{(3)} - \dfrac{\varepsilon_\infty^{-1}\Big(\lambda_\mathrm{n_0}^{(3)}\Big)^2}{\varepsilon_\mathrm{m}+ \lambda_\mathrm{n_0}^{(3)}}$ and considering $\gamma = \mathcal{O}(\delta^\mathrm{h})$ we obtain
	\begin{align}
		1-\alpha\lambda_\mathrm{n_0}^{(3)} = \mathcal{O}(\delta^\mathrm{h})  \iff \omega^2 = \omega_0^2-\omega_\mathrm{p}^2\cdot\beta^{-1} + \mathcal{O}(\delta^\mathrm{h}).
	\end{align}
	Consequently, with the calculation $\beta^{-1} = \dfrac{\varepsilon_\mathrm{m}+\lambda_\mathrm{n_0}^{(3)}}{\lambda_\mathrm{n_0}^{(3)}(\varepsilon_\infty^{-1}\varepsilon_\mathrm{m}-1)-\varepsilon_\mathrm{m}}$ we derive the needed frequency to acquired the desired equality
	\begin{align}
		1-\alpha\lambda_\mathrm{n_0}^{(3)} = \mathcal{O}(\delta^\mathrm{h}) \iff \omega^2 = \omega_0^2 - \omega_\mathrm{p}^2\dfrac{\varepsilon_\mathrm{m}+\lambda_\mathrm{n_0}^{(3)}}{\lambda_\mathrm{n_0}^{(3)}(\varepsilon_\infty^{-1}\varepsilon_\mathrm{m}-1)-\varepsilon_\mathrm{m}} + \mathcal{O}(\delta^\mathrm{h})
	\end{align}
	which completes the proof.
	\bigbreak
	\noindent
	\subsubsection{Modelling with Dielectric Nanoparticles}\label{dielectric1}
	The purpose of this section is to show under which conditions, on the incident frequency $\omega$ and the damping parameter $\gamma$, the nanoparticle behaves as a dielectric nanoparticle. Namely, the real part, and also the imaginary part, of its permittivity is large with a small ratio between the imaginary and the real parts. Such scales infer the nanoparticle to enjoy high contrast of the relative index of refraction and small relative absorption. These last properties are those usually used in the literature to describe the dielectric nanoparticles. In our analysis, we need precise choices of the incident frequency $\omega$ and damping frequency $\gamma$. 
	\bigbreak
	\noindent
	Recall the eigenvalues $\overline{\lambda}_\mathrm{n_0}^{(\boldsymbol{\ell})}$ of the Newtonian operator stated on the nanoparticle $\Omega$ in the 2D-case. We used the assumption that $\overline{\lambda}_\mathrm{n_0}^{(\boldsymbol{\ell})} \sim \delta^2 \log \delta$, with $\delta$ as the diameter of the $\Omega$. This assumption is valid for $n_0=1$, i.e. the first eigenvalue. We also recall the contrast parameter $\bm{\tau_{\mathrm{p}}}=\frac{\varepsilon_\mathrm{p}}{\varepsilon_\mathrm{\infty}}-1$. 
	\begin{lemma}
		If the frequency of the incidence wave $\omega$ is chosen close to the undamped resonance frequency $\omega^2_0$ with $\omega_0^2-\omega^2 \sim \delta^2|\log\delta|\Big(\overline{\lambda}_\mathrm{n_0}^{(\boldsymbol{\ell})}\mu_\mathrm{m}\omega_0^2\Big)$ and $\gamma\omega \sim \delta^2|\log\delta|^{1-\mathrm{h}-\mathrm{s}}\Big(\overline{\lambda}_\mathrm{n_0}^{(\boldsymbol{\ell})}\mu_\mathrm{m}\omega_0^2\Big)^2$, then we obtain $\boldsymbol{\Re}\Big(\bm{\tau_{\mathrm{p}}}\Big)\sim \delta^{-2}|\ln{\delta}|^{-1}\Big(\overline{\lambda}_\mathrm{n_0}^{(\boldsymbol{\ell})}\mu_\mathrm{m}\omega_0^2\Big)^{-1}$ and $\boldsymbol{\Im}\Big(\bm{\tau_\mathrm{p}}\Big) \sim \delta^{-2}|\log{\delta}|^{-1-\mathrm{h}-\mathrm{s}}\Big(\overline{\lambda}_\mathrm{n_0}^{(\boldsymbol{\ell})}\mu_\mathrm{m}\omega_0^2\Big)^{-1}$, where $\mathrm{h}<1$ and $\mathrm{s}>0$. 
		\bigbreak
		\noindent
		With these choices, we see that $\frac{\omega^2}{\omega_0^2}\sim 1$. In addition, we derive the needed relation $|1-\omega^2\mu_\mathrm{m}\varepsilon_\mathrm{p}\lambda_\mathrm{n_0}^{(\boldsymbol{\ell})}| \sim |\log\delta|^{-\mathrm{h}}.$
	\end{lemma}
	\noindent
	{\bf{Proof.}} From the Lorentz's model, we see that 
	\begin{align}
		\nonumber
		\bm{\tau_\mathrm{p}} &\sim \dfrac{\omega_\mathrm{p}^2}{\omega_0^2-\omega^2-i\gamma\omega} 
		\\ \nonumber &\sim \dfrac{\omega_\mathrm{p}^2(\omega_0^2-\omega^2)}{(\omega_0^2-\omega^2)^2 + (\gamma\omega)^2} +i \dfrac{\gamma\omega}{(\omega_0^2-\omega^2)^2 + (\gamma\omega)^2}
	\end{align}
	Furthermore, as we assume $\gamma\omega = o\Big(\omega_0^2-\omega^2\Big)$, then we obtain
	\begin{align}
		\bm{\tau_\mathrm{p}} &\sim \frac{\omega_\mathrm{p}^2}{\omega_0^2-\omega^2} + i \frac{\gamma\omega}{(\omega_0^2-\omega^2)^2}.
	\end{align}
	Therefore, we derive the scales $$\boldsymbol{\Re}\Big(\bm{\tau_{\mathrm{p}}}\Big)\sim \delta^{-2}|\ln{\delta}|^{-1}\Big(\overline{\lambda}_\mathrm{n_0}^{(\boldsymbol{\ell})}\mu_\mathrm{m}\omega_0^2\Big)^{-1} \mbox{ and }~  \boldsymbol{\Im}\Big(\bm{\tau_\mathrm{p}}\Big) \sim \delta^{-2}|\log{\delta}|^{-1-\mathrm{h}-\mathrm{s}}\Big(\overline{\lambda}_\mathrm{n_0}^{(\boldsymbol{\ell})}\mu_\mathrm{m}\omega_0^2\Big)^{-1}, \mbox{ where } \mathrm{h}<1.$$

\end{document}